\documentclass[3p]{elsarticle}

\usepackage{graphics}
\usepackage{graphicx}
 \usepackage{epsfig}

\usepackage{amssymb}
\usepackage{stmaryrd}
 \usepackage{amsthm}
  \usepackage{amsmath}
\usepackage{tikz-cd}

 \usepackage{float}
 \usepackage{newfloat}
\usepackage{comment}
\usepackage{ccicons}
\usepackage{epstopdf}

\usepackage{multirow}
\usepackage{color}
\usepackage{hhline}
\usepackage{booktabs}
\usepackage{rotating} 
\usepackage{hyperref}
\usepackage{bm}
\newtheorem{definition}{Definition}[section]
\newtheorem{theorem}{Theorem}[section]
\newtheorem{corollary}[theorem]{Corollary}
\newtheorem{proposition}[theorem]{Proposition}
\newtheorem{remark}{Remark}[section]

\newtheorem{lemma}{Lemma}[section]

\definecolor{nverde}{RGB}{0,61,0} 
\definecolor{cr1}{RGB}{200,0,0}
\definecolor{cr2}{RGB}{0,0,200}
\definecolor{cr12}{RGB}{100,0,100}

\newcommand{\wh}{\tilde{\mathbf{w}}_h} 
\newcommand{\vh}{\mathbf{v}_h} 
\newcommand{\ph}{\tilde{p}_h} 
\newcommand{\Zh}{\tilde{Z}_h} 
\newcommand{\Ah}{\tilde{\mathbf{A}}_h} 
\newcommand{\Yh}{Y_h} 
\newcommand{\qh}{q_h} 
\newcommand{\Bh}{\mathbf{B}_h} 
\newcommand{\Eh}{\tilde{\mathbf{E}}_h} 
\newcommand{\dx}{d \mathbf{x}} 
\newcommand{\halb}{\frac{1}{2}}
\newcommand{\n}{\mathbf{n}} 
\renewcommand{\v}{\mathbf{v}} 
\newcommand{\w}{\mathbf{w}} 
\newcommand{\A}{\mathbf{A}} 
\newcommand{\E}{\mathbf{E}} 
\newcommand{\B}{\mathbf{B}} 
\newcommand{\nph}{n+\frac{1}{2}}
\newcommand{\x}{\mathbf{x}} 

\newcommand{\CURL}{\textbf{CURL}}

\newcommand{\dpar}[2]{\frac{\partial #1}{\partial #2}}
\newcommand{\btau}{\bm \tau}
\newcommand{\N}{\mathbb N}
\newcommand{\R}{\mathbb R}

\DeclareFloatingEnvironment[fileext=lop]{Diagram}

\graphicspath{ {./} }

\journal{Journal of Computational Physics}

\allowdisplaybreaks

\begin{document}

\begin{frontmatter}

\title{A simple and general framework for the construction of exactly div-curl-grad compatible discontinuous Galerkin finite element schemes on unstructured simplex meshes}  

\author[UniZH]{R. Abgrall}
\ead{remi.abgrall@math.uzh.ch}

\author[UniTN]{M. Dumbser}
\ead{michael.dumbser@unitn.it}
\cortext[cor1]{Corresponding author}

\author[CEA]{P.-H. Maire}
\ead{pierre-henri.maire@cea.fr}

\address[UniZH]{Institute of Mathematics, University of Z\"urich, Winterthurerstrasse 190, CH-8057 Z\"urich, Switzerland} 

\address[UniTN]{Laboratory of Applied Mathematics, DICAM, University of Trento, via Mesiano 77, 38123 Trento, Italy} 

\address[CEA]{CEA CESTA, 33116 Le Barp, France}

\begin{abstract}
We introduce a new family of discontinuous Galerkin (DG) finite element schemes for the discretization of first order systems of hyperbolic partial differential equations (PDE) on unstructured simplex meshes in two and three space dimensions that respect the two basic vector calculus identities exactly also at the discrete level, namely that the curl of the gradient is zero and that the divergence of the curl is zero. The key ingredient here is the construction of two compatible discrete nabla operators, a primary one and a dual one, both defined on general unstructured simplex meshes in multiple space dimensions. Our new schemes extend existing cell-centered finite volume methods based on corner fluxes to arbitrary high order of accuracy in space. 
An important feature of our new method is the fact that only two different discrete function spaces are needed to represent the numerical solution, and the choice of the appropriate function space for each variable is related to the origin and nature of the underlying PDE. The first class of variables is discretized at the aid of a \textit{discontinuous} Galerkin approach, where the numerical solution is represented via piecewise polynomials of degree $N$ and which are allowed to jump across element interfaces. This set of variables is related to those PDE which are mere consequences of the definitions, derived from some abstract scalar and vector potentials, and for which involutions like the divergence-free or the curl-free property must hold if satisfied by the initial data. The second class of variables is discretized via classical \textit{continuous} Lagrange finite elements of approximation degree $M=N+1$ and is related to those PDE which can be derived as the Euler-Lagrange equations of an underlying variational principle. 

The primary nabla operator takes as input the data from the FEM space and returns data in the DG space, while the dual nabla operator takes as input the data from the DG space and produces output in the FEM space. The two discrete nabla operators satisfy a discrete Schwarz theorem on the symmetry of discrete second derivatives. From there, both discrete vector calculus identities follow automatically.    

We apply our new family of schemes to three hyperbolic systems with involutions: the system of linear acoustics, in which the velocity field must remain curl-free and the vacuum Maxwell equations, in which the divergence of the magnetic field and of the electric field must remain zero. In our approach, only the magnetic field will remain exactly divergence free. As a third model we study the Maxwell-GLM system of Munz \textit{et al.} \cite{MunzCleaning}, which contains a unique mixture of curl-curl and div-grad operators and in which the magnetic field may be either curl-free or divergence-free, depending on the choice of the initial data. In all cases we prove that the proposed schemes are exactly total energy conservative and thus nonlinearly stable in the $L^2$ norm. We finally apply our method to the incompressible Euler equations, which are of totally different nature than the previous PDE, but which shows that the coupling between the new DG schemes and classical standard DG methods is straightforward.  
\end{abstract}


\end{frontmatter}


\section{Introduction}\label{sec:intro}

Hyperbolic equations with stationary differential constraints (involutions) are well-known for more than one century now. The most prominent system of equations is of course the system of the Maxwell equations of electromagnetics \cite{MaxwellOrg}, where the magnetic field must remain divergence-free for all times if it was initially divergence-free. This mathematical statement is equivalent to the physical statement that there exist no magnetic monopoles. 

The first compatible discretization of the Maxwell equations in the time domain is the so-called Yee scheme \cite{Yee66}, which is based on suitably staggered meshes for the definition of the discrete electric and magnetic field in combination with compatible differential operators that guarantee that the discrete magnetic field remains always divergence-free if it was initially divergence-free. After the Yee scheme a lot of numerical methods have been developed which are able to satisfy the classical vector calculus identities that the divergence of a curl is zero and that the curl of a gradient is zero exactly also at the discrete level. Without pretending completeness, we refer the reader to the well-known mimetic finite differences \cite{HymanShashkov1997,Margolin2000,Lipnikov2014}, compatible finite volume schemes \cite{BalsaraCED,BalsaraKaeppeli,HazraBalsara} and compatible finite elements \cite{Nedelec1,Nedelec2,Hiptmair,Arnold2006,Monk,Alonso2015,CPSo2016,ddrham,BBZ22,ZAR,Zampa1} and references therein. Compatible schemes that preserve the curl-free property of a vector field exactly on the discrete level are still quite rare, but important steps forward in this direction have been made for example in \cite{JeltschTorrilhon2006,BalsaraCurlFree,SIGPR,Fambri20,SIST,Perrier3,Dhaouadi2023NSK,CurlFreeTwoFluid}. Very recently, notable progress in the field of exactly constraint-preserving finite volume and discontinuous Galerkin schemes has been also achieved in \cite{Sidilkover2025,Barsukow2024} and \cite{Perrier1,Perrier2}.  

Most of the above schemes achieve discrete compatibility with the curl and divergence involutions at the aid of suitable mesh staggering and / or the use of well-chosen function spaces in which the numerical solution is sought. And this will also be the case in the new scheme presented in this paper. 

The two basic vector calculus identities  
\begin{equation}
	\nabla \times \nabla Z = 0, \qquad \qquad 
	\nabla \cdot \nabla \times \mathbf{A} = 0,
\end{equation}
with $Z$ a scalar field and $\mathbf{A}$ a vector field 
are an immediate consequence of the Schwarz theorem on the symmetry of second partial derivatives
\begin{equation}
	\partial_j \partial_k Z   = \partial_k \partial_j Z,
	\qquad  \qquad 
	\partial_j \partial_k A_m = \partial_k \partial_j A_m.
	\label{eqn.cont.schwarz} 
\end{equation} 
This can be easily seen in Einstein index notation, assuming summation over two repeated indices, and using the classical fully antisymmetric Levi-Civita tensor 
$\epsilon_{ijk}$ for which it is well known that its contraction with a symmetric tensor is zero ($\epsilon_{ijk} T_{jk}=0$ for $T_{jk}=T_{kj}$ since $\epsilon_{ikj}=-\epsilon_{ijk}$): 
\begin{equation}
	\nabla \times \nabla Z = \epsilon_{ijk} \, \partial_j \partial_k Z = 0, 
	\qquad \qquad 
	\nabla \cdot \nabla \times \mathbf{A} = \partial_i ( \epsilon_{ijk} \partial_j A_k) = \epsilon_{ijk} \,\partial_i \partial_j  A_k = 0,
\end{equation}
due to the symmetry of second partial derivatives \eqref{eqn.cont.schwarz} and because the contraction of the antisymmetric Levi-Civita tensor with a symmetric tensor is zero.   

Our new method will mimic the above exactly also on the discrete level and is based on two main ingredients: two suitably chosen function spaces and two discrete nabla operators, a primary one and a dual one that together satisfy a discrete version of the Schwarz theorem. 

In contrast to the DG scheme introduced in \cite{Perrier1}, due to the appropriate staggering of function spaces, our method needs no numerical flux at all for the DG method and in contrast to the compatible DG schemes on staggered Cartesian meshes proposed in \cite{BoscheriCompatibleDG}  it does not require any ad hoc correction of the curl and divergence errors that arise due to the jump terms. 

The first function space is given by the discontinuous Galerkin (DG) space of piecewise polynomials of degree $N$, which are allowed to jump from one element to the other. The second function space is the one given by continuous Lagrange finite elements of approximation degree $M=N+1$. The method has thus also a weak relation to $P_NP_M$ schemes \cite{Dumbser2008}, but without reconstruction. The key idea behind our new scheme is the fact that the DG space of piecewise polynomials of degree $N$ is naturally obtained when applying the  gradient operator to the continuous finite element space of degree $M=N+1$, hence the discrete derivative applied to the discrete solution in the FEM space is actually the exact derivative applied to the discrete solution given in the continuous finite element space. The primary nabla operator is therefore applied to the FEM space and produces a discrete derivative in the DG space, while the dual nabla operator is applied to the discrete solution in the DG space and provides an output in the FEM space. The dual nabla operator is necessarily weak, since the DG polynomials are allowed to jump across element interfaces. {The chosen DG space $\mathcal{U}_h^N$ is discontinuous and of one degree less than the chosen continuous finite element space $\mathcal{W}_h^{N+1}$. However, the DG space is still rich enough to be able to represent element-wise derivatives of functions in $\mathcal{W}_h^{N+1}$ exactly, which is an essential key property that is needed in our new schemes. } 

In this paper we consider three different types of first order hyperbolic prototype systems with involutions. For simplicity of notation, we assume that the sound speed, the light speed and the cleaning speeds are unity. One of the key features of the numerical method presented in this paper is the fact that the choice of the proper function space in which the discrete solution is sought is based on the underlying variational principle from which the governing equations can be derived. In the framework of SHTC systems of Godunov and Romenski \cite{God1961,Rom1998,GodRom2003} the governing PDE come in pairs, where the first equation of each pair is a mere consequence of the definitions from some abstract scalar and vector potentials and the other equation of the pair is the Euler-Lagrange equation of the underlying variational principle. Usually, it is the first equation in the pair which is endowed with curl or divergence involutions, as the associated quantity is defined via the gradient of a scalar potential and/or the curl of a vector potential. In this paper, we use two different function spaces for the numerical discretization. The first approximation space $\mathcal{U}_h^N$ is the DG space with piecewise polynomials of degree $N$, while the second approximation space $\mathcal{W}_h^{N+1}$ is the space of continuous Lagrange finite elements (FEM) of degree $M=N+1$. We will then systematically define all quantities of the first class, with evolution equations that are mere consequences of the definitions of abstract scalar and vector potentials, in the DG space, while the quantities of the second class, whose evolution stems from the Euler-Lagrange equations of an underlying variational principle, will be discretized in the continuous finite element space. 

\paragraph{The system of linear acoustics}

Assuming unit density and unit sound speed the governing equations for linear acoustics read  
\begin{eqnarray}
	\frac{\partial \mathbf{v}}{\partial t} + \nabla p & = &  0,        \label{eqn.acoustics.v} \\
	\frac{\partial p}{\partial t} + \nabla \cdot \mathbf{v} & = &  0,  \label{eqn.acoustics.p}    
\end{eqnarray}
and satisfy the involution $\nabla \times \mathbf{v}=0$. The pressure and the velocity can be formally derived from a scalar potential $Z$ as follows: 
\begin{equation}
    \mathbf{v} = \nabla Z, \qquad \textnormal{ and } \qquad 
    p = -\partial_t Z.
    \label{eqn.acoustics.Z}
\end{equation}
Differentiating the first identity in \eqref{eqn.acoustics.Z} with respect to time and inserting the second relation immediately leads to the velocity equation \eqref{eqn.acoustics.v}. Instead, the pressure equation can be obtained as the usual Euler-Lagrange equation 
\begin{equation}
    -\frac{\partial}{ \partial t} \frac{\partial \Lambda}{\partial (\partial_t Z) } - 
    \frac{\partial}{ \partial x_k} \frac{\partial \Lambda}{\partial (\partial_k Z) } = 0,
    \label{eqn.acoustics.el}
\end{equation}
associated to the Lagrangian
\begin{equation}
    \Lambda = \halb \left[ \left( \partial_t Z \right)^2 - \left( \nabla Z \right)^2  \right] = \halb \left( p^2 - \mathbf{v}^2 \right). 
\end{equation}
For this reason in the numerical scheme presented later the velocity field, which is a mere consequence of the definitions, is discretized by the DG scheme, while the pressure, which is evolved via an Euler-Lagrange equation, is discretized via the continuous finite element method. 

\paragraph{The vacuum Maxwell equations}

The vacuum Maxwell equations for unit light speed read 
\begin{eqnarray}
	\frac{\partial \mathbf{B}}{\partial t} + \nabla \times \mathbf{E} & = &  0,  \label{eqn.maxwell.B} \\
	\frac{\partial \mathbf{E}}{\partial t} - \nabla \times \mathbf{B} & = &  0,  \label{eqn.maxwell.E}
\end{eqnarray}
and satisfy the involutions $\nabla \cdot \mathbf{B}=0$ and  $\nabla \cdot \mathbf{E}=0$. The electric and the magnetic field can be obtained from a vector potential $\mathbf{A}$ as 
\begin{equation}
    \mathbf{B} = \nabla \times \mathbf{A}, \qquad \textnormal{ and } \qquad 
    \mathbf{E} = -\partial_t \mathbf{A}. 
    \label{eqn.maxwell.A}
\end{equation}
Applying the time derivative to the first equation in \eqref{eqn.maxwell.A} and then using the definition given in the second equation directly leads to the Faraday law \eqref{eqn.maxwell.B}, as a mere consequence of the definitions given in \eqref{eqn.maxwell.A}. 
The second Maxwell equation \eqref{eqn.maxwell.E} can be obtained from the classical Euler-Lagrange equations
\begin{equation}
    -\frac{\partial}{ \partial t} \frac{\partial \Lambda}{\partial (\partial_t \textbf{A}) } - 
    \frac{\partial}{ \partial x_k} \frac{\partial \Lambda}{\partial (\partial_k \textbf{A}) } = 0,
    \label{eqn.maxwell.el}
\end{equation}
associated with the Lagrangian
\begin{equation}
    \Lambda = \halb \left[ \left( \partial_t \mathbf{A} \right)^2 - \left( \nabla \times \mathbf{A} \right)^2  \right] = \halb \left( \mathbf{E}^2 - \mathbf{B}^2 \right). 
    \label{eqn.maxwell.lambda}
\end{equation}
This is why later  in the numerical scheme the magnetic field (first PDE of the pair, i.e. a trivial consequence of the definitions) will be discretized by the DG scheme, while the electric field (second PDE of the pair, i.e. stemming from the Euler-Lagrange equations of the underlying variational principle) will be discretized via the continuous finite element method. 

\paragraph{The Maxwell-GLM system}

The Maxwell-GLM system introduced by Munz \textit{et al.} in \cite{MunzCleaning} reads  
\begin{eqnarray}
	\frac{\partial \mathbf{B}}{\partial t} + \nabla \times \mathbf{E} + \nabla p & = &  0,  \label{eqn.maxmunz.B} \\
    \frac{\partial p}{\partial t} + \nabla \cdot \mathbf{B} & = & 0, \label{eqn.maxmunz.phi} \\ 
    \frac{\partial \mathbf{E}}{\partial t} - \nabla \times \mathbf{B} + \nabla q & = &  0,  \label{eqn.maxmunz.E} \\
    \frac{\partial q}{\partial t} + \nabla \cdot \mathbf{E} & = & 0. \label{eqn.maxmunz.psi}    
\end{eqnarray}
\begin{itemize}
    \item For Maxwell-type initial data $p=q=0$ it satisfies the involutions 
    $\nabla \cdot \mathbf{B}=0$ and  $\nabla \cdot \mathbf{E}=0$, while 
    \item for acoustic-type initial data $\mathbf{E}=0$ and $q=0$ it satisfies 
    $\E=0$ and the involution $\nabla \times \mathbf{B} = 0$. 
\end{itemize}
 The field variables can be obtained from a vector potential $\mathbf{A}$ and from a scalar potential $Z$ as 
\begin{equation}
    \mathbf{B} = \nabla \times \mathbf{A} + \nabla Z, \qquad 
    \mathbf{E} = -\partial_t \mathbf{A}, \qquad  
    q = \nabla \cdot \mathbf{A}, \qquad 
    p = -\partial_t Z. \qquad  
    \label{eqn.maxmunz.AZ}
\end{equation}
Applying the time derivative to the first equation in \eqref{eqn.maxmunz.AZ} and then using the definitions given in the second and the fourth equation directly leads to the PDE \eqref{eqn.maxmunz.B}, as a mere consequence of the definitions given in \eqref{eqn.maxmunz.AZ}. Likewise, application of the time derivative to the third equation and using the definition given in the second equation of \eqref{eqn.maxmunz.AZ} leads to \eqref{eqn.maxmunz.psi}. 
The time evolution equations for the electric field \eqref{eqn.maxmunz.E} and for the scalar $p$, see \eqref{eqn.maxmunz.phi}, can be obtained from the classical Euler-Lagrange equations
\begin{eqnarray}
    -\frac{\partial}{ \partial t} \frac{\partial \Lambda}{\partial (\partial_t \textbf{A}) } - 
    \frac{\partial}{ \partial x_k} \frac{\partial \Lambda}{\partial (\partial_k \textbf{A}) } & = & 0, 
    \label{eqn.maxmunz.el.A} \\ 
    -\frac{\partial}{ \partial t} \frac{\partial \Lambda}{\partial (\partial_t Z) } - 
    \frac{\partial}{ \partial x_k} \frac{\partial \Lambda}{\partial (\partial_k Z) } & = & 0,
    \label{eqn.maxmunz.el.Z}
\end{eqnarray}
associated with the Lagrangian
\begin{equation}
    \Lambda = \halb \left[ \left( \partial_t \mathbf{A} \right)^2 - \left( \nabla \times \mathbf{A} \right)^2 - \left( \nabla \cdot \mathbf{A} \right)^2 + \left( \partial_t Z \right)^2 - \nabla Z \cdot \nabla Z \right], 
    \label{eqn.maxmunz.lambda}
\end{equation}
as recently discovered in \cite{MaxwellGLM}. 
Following the same universal reasoning as outlined above, in the numerical scheme the magnetic field $\mathbf{B}$ and the scalar $q$ (variables of the first class, with evolution equations that are a mere consequence of the definitions) will be discretized by the DG scheme, while the electric field $\mathbf{E}$ and the scalar $p$ (variables of the second class, with evolution equations given by the Euler-Lagrange equations) will be discretized via the continuous finite element method. 

\paragraph{Incompressible Euler equations} 
In addition, we also consider the nonlinear hyperbolic-elliptic system of the incompressible Euler equations with constant fluid density $\rho$, 
\begin{eqnarray}
	\frac{\partial \rho  \mathbf{v}}{\partial t} + \nabla \cdot \left( \rho \mathbf{v} \otimes \mathbf{v} \right) + \nabla p & = &  0,  \label{eqn.euler.v} \\
	 \nabla \cdot \mathbf{v} & = &  0,  \label{eqn.euler.p}
\end{eqnarray}
which is not directly related to the above equations and in which the divergence-free constraint of the velocity field is \textit{not} an involution, but the equation that is needed to determine the pressure in \eqref{eqn.euler.v}. For the ease of notation we also introduce the flux tensor containing the nonlinear convective terms $\mathbf{f} = \rho \mathbf{v} \otimes \mathbf{v}$. {Considering the incompressible Euler equations in this work allows to show that it is indeed possible to combine the new DG schemes proposed in this paper with classical DG discretizations for nonlinear time dependent PDE.}

\section{Numerical method}\label{sec:method}

The problems we consider are defined over a domain $\Omega\subset \mathbb R^d$, with $d \in \left\{2, 3 \right\}$ the number of space dimensions.
From now on, we assume that $\Omega$ has a polyhedric boundary.   This simplification is by no means essential. A triangulation / tetrahedrization of $\Omega$ is obtained  by a family of simplex elements $T_k$, $k=1, \dots, N_e$. The triangulation is assumed to be conformal. We denote by $\mathcal{E}_h$ the set of internal edges/faces of $\mathcal{T}_h$, and by $\mathcal{F}_h$ those contained in $\partial \Omega$.  
The mesh is assumed to be shape regular, $h_T$ represents the diameter of the element $T=T_k$. Similarly, if $e\in \mathcal{E}_h\cup \mathcal{F}_h$, $h_e$ represents its diameter.

 Throughout this paper, we follow Ciarlet's definition \cite{ciarlet,ErnGuermond} of a finite element approximation: we have a set of degrees of freedom $\Sigma_T$ of linear forms acting on the set $\mathbb P_N(T)$ of polynomials of degree $N$  defined on $T$ such that the linear mapping
 $$q\in \mathbb P_N\mapsto \big (\sigma_1(q), \ldots, \sigma_{|\Sigma_T|}(q)\big )$$
 is one-to-one. The space $\mathbb P_N$ is spanned by basis functions that are interpolant at the degrees of freedom in $T$. 
  We have in mind  Lagrange interpolations of degree $N$ where the degrees of freedom are associated to points in $T$. However, this is not essential.
  For the sake of simplicity, if $\{\sigma_j\}_{\Sigma_T}$ is the set of degrees of freedom in $T$, we will denote them simply by their index $j$ when there is no ambiguity. 


The DG space is denoted by $\mathcal{U}^N_h$. Its precise definition depends on the problem under consideration, as well as the space dimension. In all cases, 
$$\mathcal{U}^N_h=\bigoplus\limits_{T_k\in \mathcal{T}}\{ u_h \text{ is polynomial of degree }N\text{ in }T_k\text{ with coefficients in }\mathbb R^m\}$$ where $m$ is the number of components of the vectors needed to set up the PDE. There is no continuity requirement across the faces or edges of the elements $T_k$. We denote by ${\phi}_c$ the associated basis or ansatz functions and by $u_c$ the corresponding degrees of freedom, so in $T_k$,
$$u_h = \sum_{c} \phi_c(\x) u_c, \quad u_c \in \mathbb R^m.$$ In principle, any basis is fine. In practice, we use a classical nodal Lagrange FEM basis of degree $N$.

Instead, the continuous FEM space is denoted by $\mathcal{W}^{N+1}_h$. Again, its precise definition depends on the problem under consideration, as well as the space dimension. We define, {for $M\in \N$}
$${\mathcal{W}^{M}_h=\mathcal{U}^{M}_h\cap \big (C^0(\Omega)\big )^m}.$$
Here global continuity is enforced across the faces /edges of the elements,  which is why we need the mesh to be conformal. The associated basis functions  are denoted by ${\psi}_p$. In practice, they are chosen as classical nodal Lagrange FEM basis of degree $M=N+1$, {i.e. we consider $\mathcal{W}^{N+1}_h$}. The generic data representation with basis functions $\psi_p$ and degrees of freedom $w_p$ reads 
$$\tilde{w}_h = \sum_{p} \psi_p(\x) \, w_p, \quad w_p \in \mathbb R^m.$$

{At this point we stress again that the judicious choice of the discrete ansatz spaces for the different kind of variables arising in the PDE systems treated in this manuscript is actually one of the major key ingredients of the present paper: 
i) those PDE that are \textbf{mere consequences of the definitions} and which can be simply derived from the definitions of some abstract scalar or vector potentials are discretized via a \textbf{DG scheme} in the discontinuous ansatz space $\mathcal{U}^N_h$, while 
ii) the PDE that originate from a \textbf{variational principle} via the \textbf{Euler-Lagrange equations} are discretized in the continuous finite element space, i.e. in $\mathcal{W}^{N+1}_h$. Note that this statement does not hold for the incompressible Euler equations, which do not fall in the general class of SHTC systems treated in this paper.} 

{We also would like to underline that for all linear PDE considered in this paper we compute all integrals exactly, since the integration of piecewise polynomials is easily feasible on simplex elements. For nonlinear PDE we employ quadrature rules of sufficient accuracy. In the proofs of all theorems below, we assume exact integration. } 

\subsection{Pair of discretely compatible nabla operators}
\label{sec.disc.nabla}

Making use of the two approximation spaces $\mathcal{W}^{N+1}_h$ and $\mathcal{U}^N_h$ we now introduce a pair of discretely compatible nabla operators. In the following we will use the subscripts $a$ and $c$ for the DG degrees of freedom in a cell (triangle / tetrahedron $T_k$) and the subscripts $p$ and $q$ for the nodal degrees of freedom of the continuous FEM. Furthermore, the tensor indices of  physical fields are denoted by $i$, $j$, $k$, $l$ and $m$, respectively. 

\subsubsection{Discrete primary nabla operator} 

We define the rank three stiffness tensor as 
$$ \mathbb K = \left\{ \mathbf{K}_{cp} \right\} = \int \limits_{\Omega} \phi_c \nabla \psi_p d\mathbf{x}, $$  
or, equivalently, in Einstein index notation, 
$$
K_{cpm} = \int \limits_{\Omega} \phi_c \partial_m \psi_p d\mathbf{x}.
$$
Note that in a DG scheme, the $\phi_c$ are allowed to jump across the elements. To get the discrete primary nabla operator we must multiply with the inverse of the element-local mass matrix of the DG scheme, which is easy to compute in the DG setting (when using the Dubiner basis \cite{Dubiner} for $\phi_c$, it can even be made diagonal),      
%
%
$$ D_{ac} = \int \limits_{\Omega} \phi_a \phi_c d\mathbf{x}  $$
i.e. 
\begin{equation} 
\nabla_c^p = D_{ca}^{-1} \mathbf{K}_{ap}, 
\qquad \textnormal{or, in index notation,} \qquad  
\left(\partial_m\right)_c^p = D_{ca}^{-1} K_{apm}.
\label{eqn.primary.nabla}
\end{equation} 
Let now $Z$ be a function defined on $\Omega$ with values in $\mathbb R$ and its gradient $\v = \nabla Z$ with values in $\mathbb R^d$. The approximation of $Z$ in $\mathcal{W}^{N+1}_h$ is denoted by $\Zh$, with degrees of freedom ${Z}_q$ and representation 
\begin{equation} 	
\Zh = \sum_q \psi_q(\mathbf{x}) {Z}_q:= \psi_q Z_q \qquad \in \mathcal{W}^{N+1}_h,
\label{eqn.Z.ansatz}
\end{equation} 	
while the discrete gradient is denoted by $\vh \in \mathcal{U}^N_h$ and representation
\begin{equation} 	
\vh = \sum_a \phi_a(\mathbf{x}) \v_a := \phi_a \v_a = \phi_a v_{ma} \qquad \in \mathcal{U}^N_h,
\label{eqn.v.ansatz}
\end{equation} 	
where the index $m$ in $v_{ma}$ refers to the components of the velocity field and the index $a$ to the degrees of freedom. 
The degrees of freedom of the discrete gradient $\nabla_h$ applied to $\tilde{Z}_h$ then simply read 
\begin{equation}
\v_c = \nabla_c^p Z_p = \sum_{p} D_{ca}^{-1} \mathbf{K}_{ap} Z_p,  
\qquad 
\textnormal{ or, in index notation, } 
\qquad 
v_{mc} = \left( \partial_m \right)_c^p Z_p = D_{ca}^{-1} K_{apm} Z_p.    
\end{equation}
We stress that thanks to the judicious choice of the discrete function spaces for $\vh$ and $\Zh$, the primary discrete nabla operator is \textit{exact}, i.e. $\vh = \nabla_h \Zh = \phi_c \nabla_c^p Z_p = \nabla \Zh$, with $\nabla$ the exact nabla operator.  

For the sake of completeness, we also report the application of the discrete primary nabla operator to a vector field $\Ah = \psi_q \mathbf{A}_{q} = \psi_q A_{kq} \in \mathcal{W}^{N+1}_h$, with $k$ the tensor index of the physical field and $q$ the index for the degrees of freedom:
\begin{equation}
 \nabla_c^p \mathbf{A}_{p} = \sum_{p} D_{ca}^{-1} \mathbf{K}_{ap} \mathbf{A}_{p}, 
 \qquad \textnormal{ or, in index notation } \qquad 
 \left( \partial_m \right)_c^p A_{kp} = D_{ca}^{-1} K_{apm} A_{kp}.  
\end{equation}

For $N=0$ we have $\phi_c=1$, the inverse mass matrix is one divided by the cell volume and $\nabla \psi_p$ is the gradient of continuous P1 Lagrange FEM
and we obtain the discrete primary nabla operator used in the Lagrangian schemes \cite{Despres2005,Despres2009,Maire2007,Maire2009,LMR2016,Maire2020,Boscheri_hyperelast_22,Boscheri24_mhd,HTCLagrange}. 

\subsubsection{Discrete dual nabla operator} 

The discrete dual nabla operator can be obtained via the integration by parts formula
$$
 \int \limits_\Omega \psi_p \nabla \phi_c \, d\x = 
 \int \limits_{\Omega \backslash (\mathcal{E}_h \cup \mathcal{F}_h)} \psi_p \nabla \phi_c \, d\x + \int \limits_{\mathcal{E}_h \cup \mathcal{F}_h} \psi_p \left( \phi_c^+ - \phi_c^- \right) \, \n \, dS 
 = 
 - \int \limits_\Omega \nabla \psi_p  \, \phi_c \, d\x + \int \limits_{\partial \Omega}  \psi_p \phi_c \, \n \, dS.
$$
Assuming suitable boundary conditions on $\partial \Omega$ that allow to drop the boundary contributions on $\partial \Omega$, this simplifies to 
$$
 \int \limits_\Omega \psi_p \nabla \phi_c d\x = 
 \int \limits_{\Omega \backslash \mathcal{E}_h} \psi_p \nabla \phi_c \, d\x + \int \limits_{\mathcal{E}_h} \psi_p \left( \phi_c^+ - \phi_c^- \right) \, \n \, d\x 
 = 
 - \int \limits_\Omega \nabla \psi_p  \, \phi_c \, d\x. 
$$

We then have a discrete dual nabla operator that reads 

\begin{equation} 
\tilde{\nabla}_p^c = -\mathbf{K}_{cp}, 
\qquad \textnormal{or, in index notation,} \qquad  
(\tilde{\partial}_m)_p^c = - K_{cpm},
\label{eqn.dual.nabla}
\end{equation} 
and which deliberately does not yet contain the multiplication by the inverse of the global FEM mass matrix 
\begin{equation}
    M_{pq} = \int \limits_{\Omega} \psi_p \psi_q \dx.     
\end{equation}
In the case $N=0$ we have $\phi_c=1$, hence we recover the discrete dual nabla operator as in \cite{Despres2005,Despres2009,Maire2007,Maire2009,LMR2016,Maire2020,Boscheri_hyperelast_22,Boscheri24_mhd,HTCLagrange}.      

Let now $Y$ be a function defined on $\Omega$ with values in $\mathbb R$ and its gradient $\w = \nabla Y$ with values in $\mathbb R^d$. The approximation of $Y$ in $\mathcal{U}^N_h$ is denoted by $\Yh$, with degrees of freedom ${Y}_a$ and representation 
\begin{equation} 	
\Yh = \sum_a \phi_a(\mathbf{x}) {Y}_a:= \phi_a Y_a \qquad \in \mathcal{U}^N_h,
\label{eqn.Y.ansatz}
\end{equation} 	
while the discrete gradient is denoted by $\wh \in \mathcal{W}^{N+1}_h$ and with representation
\begin{equation} 	
\wh = \sum_p \psi_p(\mathbf{x}) \w_p := \psi_p \w_p = \psi_p w_{mp} \qquad \in \mathcal{W}^{N+1}_h,
\label{eqn.w.ansatz}
\end{equation} 	
where the index $m$ in $w_{mp}$ refers to the spatial components of the vector field and the index $p$ to the associated degrees of freedom of the discrete solution $\wh$. The discrete dual gradient $\tilde{\nabla}_h$ applied to $\Yh$ reads 
\begin{equation}
M_{pq} \w_q = \tilde{\nabla}_p^c Y_c = -\sum_{c} \mathbf{K}_{cp} Y_c,  
\qquad 
\textnormal{ or, in index notation, }  
\qquad 
M_{pq} w_{mq} = ( \tilde{\partial}_m )_p^c = - K_{cpm} Y_c.    
\end{equation}
Again for completeness, we also provide the discrete dual nabla operator applied to a vector field $\mathbf{A}_h = \phi_a A_{ka} \in \mathcal{U}^N_h$, with $k$ the tensor index of the physical field and $a$ the index for the degrees of freedom:
\begin{equation}
 \tilde{\nabla}_p^c \mathbf{A}_{c} = -\sum_{c} \mathbf{K}_{cp} \mathbf{A}_{c}, 
 \qquad \textnormal{ or, in index notation } \qquad 
 ( \tilde{\partial}_m )_p^c A_{kc} =  - K_{cpm} A_{kc}.  
\end{equation}

{
\subsubsection{Basic properties} 
We first present two basic properties of the primary discrete nabla operator introduced in this paper. Since the continuous finite element space $\mathcal{W}_h^{N+1}$ is globally continuous, it is obvious that all \textit{tangential derivatives}, i.e. the tangential components of the gradient, are \textit{continuous} across elements. Furthermore, the DG space $\mathcal{U}_h^N$ is rich enough to represent all first derivatives of functions represented in the finite element space $\mathcal{W}_h^{N+1}$ exactly. We therefore have the following 
\begin{proposition}
	\label{prop.td} 
	For any discrete scalar field $\Zh \in \mathcal{W}_h^{N+1}$ and any pair of elements $T_b$ and $T_c$  the jump of the tangential derivatives vanishes across the common edge / face $\partial T_{bc} = T_b \cap T_c$, i.e.  
	\begin{equation}
		 \left( \nabla \Zh^+ - \nabla \Zh^- \right) \cdot \mathbf{t} = 0, 
 	\end{equation} 
	with $\mathbf{t}$ all vectors tangential to $\partial T_{bc}$ and $\Zh^+$, $\Zh^-$ the boundary-extrapolated values of $\Zh$ seen from elements $T_b$ and $T_c$, respectively.  
\end{proposition} 
\begin{proof}
	Due to the \underline{global continuity} of $\Zh \in \mathcal{W}_h^{N+1}$ one has clearly  
	$$ \Zh^+ = \Zh^-, \quad \textnormal{ or, equivalently } \quad \Zh^+ - \Zh^- = 0, \qquad \forall \mathbf{x} \in \partial T_{bc}.$$
	It is therefore obvious that for all vectors $\mathbf{t}$ tangential to $\partial T_{bc}$ one has 
	$$\frac{\partial \Zh^+}{\partial \mathbf{t}} = \frac{\partial \Zh^-}{\partial \mathbf{t}}, \qquad \textnormal{ or, equivalently } \qquad 
	\left( \nabla \Zh^+ - \nabla \Zh^- \right) \cdot \mathbf{t} = 0, \qquad \forall \mathbf{x} \in \partial T_{bc}.
	$$ 
\end{proof}
As direct \textit{consequence} of Proposition \ref{prop.td} we also immediately have the property that the \textit{normal component} of the curl of a vector field $\Ah \in \mathcal{W}_h^{N+1}$ is \textit{continuous} across elements. 
\begin{proposition}
	\label{prop.nd} 
	For any discrete vector field $\Ah \in \mathcal{W}_h^{N+1}$ and any pair of elements $T_b$ and $T_c$  the jump of the normal component of the curl vanishes across the common edge / face $\partial T_{bc} = T_b \cap T_c$, i.e.  
	\begin{equation}
		\left( \nabla \times \Ah^+ - \nabla \times \Ah^- \right) \cdot \mathbf{n} = 0,
	\end{equation} 
	with $\mathbf{n}$ the unit normal vector of $\partial T_{bc}$ pointing from $T_c$ to $T_b$ and $\Ah^+$, $\Ah^-$ the boundary-extrapolated values of $\Ah$ seen from elements $T_b$ and $T_c$, respectively.  
\end{proposition} 
\begin{proof}
	The proof becomes easier when using the Einstein index notation, i.e. here we will denote the components of the vector field $\Ah$ by $\tilde{A}_p$ and drop the subscript $h$ to ease notation. \\ Furthermore, the components of the tangential vectors $\mathbf{t}$ will be denoted by $t_k$ and those of the unit normal vector $\mathbf{n}$ by $n_k$. 
	Due to the \underline{global continuity} of $\Ah \in \mathcal{W}_h^{N+1}$ one has clearly componentwise   
	$$ \tilde{A}_p^+ = \tilde{A}_p^-, \quad \textnormal{ or } \quad \tilde{A}_p^+ - \tilde{A}_p^- = 0, \qquad \forall \mathbf{x} \in \partial T_{bc}.$$
	It is therefore obvious that for all vectors $\mathbf{t}$ tangential to $\partial T_{bc}$ one has componentwise 
	$$	\left( \partial_j \tilde{A}_p^+ - \partial_j \tilde{A}_p^- \right) t_j = 0, \qquad \forall \mathbf{x} \in \partial T_{bc}.
	$$ 
	Hence, in general, the only non-vanishing component of the jump in the gradient of $\Ah$ at the element boundary $\partial T_{bc}$ is the jump of the derivative in the normal direction, i.e. 
	$$ \left( \partial_j  \tilde{A}_p^+ - \partial_j  \tilde{A}_p^- \right) = \partial_j \left( \tilde{A}_p^+ - \tilde{A}_p^- \right) = \partial_m \left( \tilde{A}_p^+ - \tilde{A}_p^- \right) n_m \, n_j. $$
	Using the above result and the usual fully antisymmetric Levi-Civita tensor $\epsilon_{ijk}$ in Einstein index notation it is now easy to see that 
	$$ 
	\left( \nabla \times \Ah^+ - \nabla \times \Ah^- \right) \cdot \mathbf{n} = 
	\epsilon_{ijk} \left( \partial_j \tilde{A}_k^+ - \partial_j \tilde{A}_k^- \right) n_i = 
	\epsilon_{ijk} \left( \partial_m \tilde{A}_k^+ - \partial_m \tilde{A}_k^- \right) n_m n_i n_j = 0, 
	$$ 
	since the contraction of the anti-symmetric Levi-Civita tensor with the symmetric tensor $n_i n_j$ vanishes, i.e. $\epsilon_{ijk} n_i n_j = 0$. 
\end{proof}
}

\subsubsection{Discrete vector calculus identities} 

We now prove that the two basic vector calculus identities hold exactly on the discrete level thanks to a discrete analogue of the Schwarz theorem. 

{
\begin{lemma}
	\label{lemma.ds}
 The two discrete nabla operators \eqref{eqn.primary.nabla} and \eqref{eqn.dual.nabla} satisfy the following discrete analogue of the Schwarz theorem
 \begin{equation}
 	  (\tilde{\partial_j})_p^c  \left(\partial_k \right)_c^q Z_q = (\tilde{\partial_k})_p^c  \left(\partial_j\right)_c^q Z_q.
 	  \label{eqn.ds.scalar} 
 \end{equation}
\end{lemma}
\begin{proof}
    To ease notation and for the sake of convenience in the following we use the Einstein index notation. From the definition of the two discrete nabla operators \eqref{eqn.primary.nabla} and \eqref{eqn.dual.nabla} we have 
\begin{equation} 
	-K_{cpj}  D_{ca}^{-1} K_{aqk} Z_q = (\tilde{\partial}_j)_p^c \, (\partial_k)_c^q Z_q = - \int \limits_{\Omega} \nabla \psi_p \nabla \Zh \dx = \int \limits_{\Omega \backslash \mathcal{E}_h} \psi_p \partial_j \partial_k \Zh \dx 
	+ \int \limits_{\mathcal{E}_h} \psi_p  \left( \partial_k \Zh^+  - \partial_k \Zh^- \right)  n_j \, dS.
	\label{eqn.auxeq}
\end{equation}
Here, $\mathcal{E}_h$ is the skeleton of the mesh composed of the union of all internal edges / faces. 
Since $\Zh$ is smooth inside each simplex $T_k$, the classical continuous Schwarz theorem obviously holds, $\partial_j \partial_k \Zh = \partial_k \partial_j \Zh$, and due to the global continuity of $\Zh$ across element edges / faces $\nabla \Zh$ is allowed to jump across $\mathcal{E}_h$ only in the normal direction and not in the tangential directions, hence also the jump term 
$$ \left( \partial_k \Zh^+  - \partial_k \Zh^- \right)  n_j = \left( \partial_m \Zh^+ - \partial_m \Zh^- \right) n_m n_k n_j = \left( \partial_j \Zh^+  - \partial_j \Zh^- \right)  n_k $$ is obviously symmetric in the indices $j$ and $k$, see also the proof of Proposition \ref{prop.nd}. 
The symmetry in the indices $j$ and $k$ in \eqref{eqn.auxeq} is now clear, i.e. we have the sought discrete analogue of the Schwarz theorem that reads
\begin{eqnarray} 
	&& (\tilde{\partial}_j)_p^c \, (\partial_k)_c^q Z_q = -K_{cpj}  D_{ca}^{-1} K_{aqk} Z_q = 
	\int \limits_{\Omega \backslash \mathcal{E}_h} \psi_p \partial_j \partial_k \Zh \dx 
	+ \int \limits_{\mathcal{E}_h} \psi_p  \left( \partial_k \Zh^+  - \partial_k \Zh^- \right)  n_j \, dS \nonumber \\ 
	&& =   
	\int \limits_{\Omega \backslash \mathcal{E}_h} \psi_p \partial_k \partial_j \Zh \dx 
	+ \int \limits_{\mathcal{E}_h} \psi_p  \left( \partial_j \Zh^+  - \partial_j \Zh^- \right)  n_k \, dS 
	= -K_{cpk}  D_{ca}^{-1} K_{aqj} Z_q = (\tilde{\partial}_k)_p^c \, (\partial_j)_c^q Z_q. 
	\label{eqn.discrete.schwarz}
\end{eqnarray}     	 
\end{proof}
As an immediate consequence of the above Lemma we have the following 
\begin{corollary}
	\label{corr.ds}
	The two discrete nabla operators \eqref{eqn.primary.nabla} and \eqref{eqn.dual.nabla} also satisfy a discrete analogue of the Schwarz theorem for discrete vector fields 
	\begin{equation}
		(\tilde{\partial_j})_p^c  \left(\partial_k \right)_c^q A_{mq} = (\tilde{\partial_k})_p^c  \left(\partial_j\right)_c^q A_{mq}.
		\label{eqn.ds.vector} 
	\end{equation}
\end{corollary}
\begin{proof}
	 Replacing $Z_q$ by $A_{mq}$ in the proof of the previous Lemma \eqref{lemma.ds} leads immediately to the sought result. 
\end{proof}
}

\begin{theorem}
Given the two discrete nabla operators \eqref{eqn.primary.nabla} and \eqref{eqn.dual.nabla}, 
a discrete scalar potential $\Zh \in \mathcal{W}^{N+1}_h$ with degrees of freedom $Z_q$ and a discrete vector potential $\Ah \in \mathcal{W}^{N+1}_h$ with degrees of freedom $\A_q$ and representation
\begin{equation}
    \Zh = \psi_q  Z_q \quad \in \mathcal{W}^{N+1}_h, \qquad 
    \Ah = \psi_q \A_q \quad \in \mathcal{W}^{N+1}_h, 
    \label{eqn.disc.potential}
\end{equation}
the following discrete vector calculus identities hold:
\begin{equation}
   { \tilde{\nabla}_p^c \times \left( \nabla_c^q Z_q \right) = 0,}
    \label{eqn.disc.rotgrad}
\end{equation}
and
\begin{equation}
   { \tilde{\nabla}_p^c \cdot 
   \left( \nabla_c^q \times \mathbf{A}_q \right)} = 0.
    \label{eqn.disc.divrot}
\end{equation}
\end{theorem}
\begin{proof}
    In the following we use the Einstein index notation. We also make use of the classical fully antisymmetric Levi-Civita tensor $\epsilon_{ijk}$, which is well-known and whose definition does not need to be repeated here. 
    We first prove the identity \eqref{eqn.disc.rotgrad}: 
    \begin{equation}
        \tilde{\nabla}_p^c \times \nabla_c^q Z_q = 
        \epsilon_{ijk} \, (\tilde{\partial}_j)_p^c \, (\partial_k)_c^q Z_q 
        = -\epsilon_{ijk} \, K_{cpj} \, D_{ca}^{-1} \, K_{aqk} \, Z_q = 0, 
        \label{eqn.disc.rotgrad2}
    \end{equation}    
    where we have simply used the definitions of the two compatible nabla operators \eqref{eqn.primary.nabla} and \eqref{eqn.dual.nabla} in index notation, and the fact that the contraction of the antisymmetric Levi-Civita symbol with a symmetric tensor is zero. Clearly $K_{cpj} D_{ca}^{-1} K_{aqk}  = K_{cpk}  D_{ca}^{-1} K_{aqj} Z_q$ is symmetric in the indices $j$ and $k$, thanks to the discrete Schwarz theorem \eqref{eqn.ds.scalar} of the previous Lemma \ref{lemma.ds}, i.e. $(\tilde{\partial}_j)_p^c \, (\partial_k)_c^q
    = (\tilde{\partial}_k)_p^c \, (\partial_j)_c^q$ and obviously one has $\epsilon_{ijk} T_{jk} = 0$ if $T_{jk}=T_{kj}$ since the contraction of an anti-symmetric tensor with a symmetric one vanishes.

    \bigskip 

    We now prove the second identity \eqref{eqn.disc.divrot}: 
    \begin{equation}
        \tilde{\nabla}_p^c \cdot \nabla_c^q \times \mathbf{A}_q = 
        (\tilde{\partial}_i)_p^c \, \epsilon_{ijk}  (\partial_j)_c^q A_{kq}
        = - \epsilon_{ijk} K_{cpi} \, D_{ca}^{-1} \, K_{aqj} A_{kq} = 0,
        \label{eqn.disc.divrot2}
    \end{equation}    
    for the same reasons as before, i.e. simply using 
    \eqref{eqn.primary.nabla} and \eqref{eqn.dual.nabla} in index notation, the discrete Schwarz theorem for vector fields \eqref{eqn.ds.vector} and knowing that the contraction of the antisymmetric Levi-Civita symbol with a symmetric tensor is zero. 
\end{proof}
{
\paragraph{Proper choice of initial data} 
We emphasize that the discrete solution in the DG space $\mathcal{U}_h^N$ \textit{cannot} be initialized via simple element-local $L^2$ projection. Also in this case the guiding principle follows the same philosophy as the one that dictates the choice of ansatz spaces used throughout this paper. Since the DG space is used for those variables which are defined via the gradient or the curl of some abstract scalar and vector potentials, the initial condition in the DG space $\mathcal{U}_h^N$ must necessarily be also properly computed as the discrete gradient or discrete curl of the corresponding discrete potential in the continuous finite element space $\mathcal{W}_h^{N+1}$. }

\subsection{Compatible semi-implicit DG scheme for acoustics}

In this section we present our new approach for the system of linear acoustics \eqref{eqn.acoustics.v} and \eqref{eqn.acoustics.p}. Following the universal guiding principle outlined above, according to which the quantities whose evolution equations are the Euler-Lagrange equations are sought in the FEM space $\mathcal{W}^{N+1}_h$ and the quantities based on mere consequences of the definitions are located in the DG space $\mathcal{U}^N_h$, the discrete solution of the velocity is sought in the DG space $\mathcal{U}^N_h$, while the discrete pressure field is sought in the continuous FEM space $\mathcal{W}^{N+1}_h$: 
\begin{equation} 	
\vh = \sum_j \phi_j(\mathbf{x}) {\mathbf{v}}_j \quad \in \mathcal{U}^N_h, \qquad \qquad 
\ph = \sum_j \psi_j(\mathbf{x}) {p}_j \quad \in \mathcal{W}^{N+1}_h.
\label{eqn.acoustics.ansatz}
\end{equation} 	
We have $N_{dof}^v$ vectorial degrees of freedom for the velocity and $N_{dof}^p$ pressure degrees of freedom. 
One could therefore call the presented scheme a staggered method, but with the staggering not realized in the mesh, but rather in the functional spaces used for the representation of the discrete solution.  

We first illustrate the approach in the semi-discrete setting, before discussing also the time discretization. Multiplication of the velocity equation  \eqref{eqn.acoustics.v} by a test function $\phi_i \in \mathcal{U}^N_h$ and integrating over $T_k$, for all $T_k$ the weak formulation of the velocity equation becomes,
\begin{equation}	
	\int \limits_{T_k} \phi_i \partial_t \vh \, \dx + 
	\int \limits_{T_k} \phi_i \nabla \ph \, \dx = 0.  
    \label{eqn.acoustics.vh}
\end{equation} 
Since the pressure field $\ph \in \mathcal{W}^{N+1}_h$ is \textit{continuous} across element boundaries, nothing more  needs to be done for Eqn. \eqref{eqn.acoustics.vh}.

 Likewise, the velocity field $\vh \in \mathcal{U}^N_h$ is \textit{discontinuous} across element boundaries and therefore  \eqref{eqn.acoustics.p} 
 needs to be interpreted in the sense of distributions, i.e. for any test function $\psi_i \in \mathcal{W}^{N+1}_h$
\begin{equation}
    	\int \limits_{\Omega} {\psi}_i \partial_t \ph \, \dx -  
	\int \limits_{\Omega} \nabla {\psi}_i \cdot \vh \, \dx + 
     \int \limits_{\partial \Omega} {\psi}_i \mathbf{v} \cdot \n \, d \gamma = 0,  
    \label{eqn.acoustics.ph.global}
\end{equation}
{where $\mathbf{v} \cdot \n$ on $\partial \Omega$ is the imposed physical boundary condition of the original continuous problem. Since the DG space explicitly admits discontinuities at element boundaries, this choice is indeed possible.}
Summing the discrete velocity equation over all triangles immediately leads to 
\begin{equation}
	\int \limits_{\Omega} \phi_i \partial_t \vh \, \dx + 
	\int \limits_{\Omega} \phi_i \nabla \ph \, \dx = 0.  
    \label{eqn.acoustics.vh.global}    
\end{equation}
{\paragraph{Proper choice of discrete initial data} 
	In order to guarantee compatible discrete initial data for the DG method, the 
	discrete initial velocity field $\v_h(\x,0) \in \mathcal{U}^N_h$ must be computed as 
	the discrete gradient of a discrete	scalar potential 
	$\tilde{Z}_h = \psi_q Z_q \in \mathcal{W}^{N+1}_h$ defined in the continuous finite element space.  The degrees of freedom of the initial velocity field are then simply given by $\v^0_c = \nabla_c^q Z_q$, with the primary nabla operator $\nabla_c^q$ defined above.  
	We stress that simple elementwise $L^2$ projection will in general lead to incompatible 
	discrete initial data for the DG method. Concerning the initial data of the discrete pressure, which is defined in the continuous finite element space $\mathcal{W}^{N+1}_h$, standard nodewise evaluation of the initial data is sufficient.}
In the following we now prove global total energy conservation of the semi-discrete scheme for the boundary condition $\mathbf{v} \cdot \n = 0$ on $\partial \Omega$: 
\begin{theorem}
    For impermeable boundary conditions $\mathbf{v} \cdot \n = 0$ on $\partial \Omega$ the scheme \eqref{eqn.acoustics.vh} and \eqref{eqn.acoustics.ph.global} conserves global total energy in the sense 
    \begin{equation}	
	   \partial_t \int \limits_{\Omega}  \frac{1}{2} \left( \vh^2 + \ph^2 \right) \, \dx = 0.   
    \end{equation} 
\end{theorem}
\begin{proof}
    Since $\phi_i \in \mathcal{U}^N_h$ and $\psi_i \in \mathcal{W}^{N+1}_h$ we can also use $\vh$ and $\ph$ as test functions in \eqref{eqn.acoustics.vh.global} and \eqref{eqn.acoustics.ph.global}, thus obtaining
\begin{equation}
	\int \limits_{\Omega} \vh \cdot \partial_t \vh \, \dx + 
	\int \limits_{\Omega} \vh \cdot \nabla \ph \, \dx = 0, 
    \label{eqn.acoustics.vh.vh}    
\end{equation}
and 
\begin{equation}
    	\int \limits_{\Omega} \ph \partial_t \ph \, \dx -  
	\int \limits_{\Omega} \nabla \ph  \cdot \vh \, \dx = 0,  
    \label{eqn.acoustics.ph.ph}
\end{equation}    
where we have also used the boundary condition $\mathbf{v} \cdot \n = 0$ on $\partial \Omega$. 
Summing up the above equations \eqref{eqn.acoustics.vh.vh} and \eqref{eqn.acoustics.ph.ph} yields
\begin{equation} 
\int \limits_{\Omega} \vh \cdot \partial_t \vh \, \dx 
+ \int \limits_{\Omega} \ph \partial_t \ph \, \dx
+ \int \limits_{\Omega} \vh \cdot \nabla \ph \, \dx 
- \int \limits_{\Omega} \nabla \ph  \cdot \vh \, \dx =  
\partial_t \int \limits_{\Omega}  \frac{1}{2} \left( \vh^2 + \ph^2 \right) \, \dx = 0. 
\end{equation}
\end{proof}
Using a classical Crank-Nicholson time discretization with the usual notation 
\begin{equation} 
\vh^{\nph} = \halb \left( \vh^{n} + \vh^{n+1} \right), 
\qquad \textnormal{ and } \qquad 
\ph^{\nph} = \halb \left( \ph^{n} + \ph^{n+1} \right) 
\end{equation}
and assuming again $\mathbf{v} \cdot \n = 0$ on $\partial \Omega$, let us show that we also get a discrete energy equality. Indeed, 
\begin{equation}
	\int \limits_{\Omega} \phi_i \left( \vh^{n+1} - \vh^n \right) \, \dx + 
	\Delta t \int \limits_{\Omega} \phi_i \nabla \ph^{\nph} \, \dx = 0,  
    \label{eqn.acoustics.vh.fd}    
\end{equation}
\begin{equation}
    	\int \limits_{\Omega} {\psi}_i \left( \ph^{n+1} - \ph^n \right) \, \dx -  
	\Delta t \int \limits_{\Omega} \nabla {\psi}_i \cdot \vh^{\nph} \, \dx = 0.  
    \label{eqn.acoustics.ph.fd}
\end{equation}
i.e. our new compatible DG scheme becomes: 
\label{acoustic:CN}
\begin{equation}	
	\int \limits_{T_k} \phi_i  \, \frac{{\mathbf{v}}^{n+1}-{\mathbf{v}}^{n}}{\Delta t}\; d \mathbf{x}  + 
	\int \limits_{T_k} \phi_i  \, \tilde{p}^{n+\frac{1}{2}} \; d \mathbf{x} = 0  
	\label{eqn.dg.v} 
\end{equation} 
\begin{equation}	
	\int \limits_{\Omega} {\psi}_i \frac{\tilde{p}^{n+1}-\tilde{p}^{n}}{\Delta t}\; d \mathbf{x}  -   
	\int \limits_{\Omega} \nabla {\psi}_i \cdot {\mathbf{v}}^{n+\frac{1}{2}} \; d \mathbf{x} 
	= 0.  
	\label{eqn.fem.p} 
\end{equation} 
Multiplying \eqref{eqn.dg.v} by ${\mathbf{v}}_j^{n+\frac{1}{2}}$ and \eqref{eqn.fem.p} by 
$ \tilde{p}_j^{n+\frac{1}{2}}$, summing up, and summing the outcomes of \eqref{eqn.dg.v} and \eqref{eqn.fem.p}, the spatial terms disappear, hence
$${\mathbf{v}}^{n+\frac{1}{2}}\cdot \big ( {\mathbf{v}}^{n+1}-{\mathbf{v}}^{n}\big ) +
\tilde{p}^{n+\frac{1}{2}} \big ( \tilde{p}^{n+1}-\tilde{p}^{n}\big )=
\frac{1}{2} \bigg ( \big ({\mathbf{v}}^{n+1}\big )^2+\big ( \tilde{p}^{n+1}\big )^2\bigg )-\frac{1}{2}\bigg ( \big ({\mathbf{v}}^{n}\big )^2+\big ( \tilde{p}^{n}\big )^2\bigg ) = 0, $$
{which is the fully-discrete conservation law of total energy. Note that the particular choice of ansatz spaces $\mathcal{U}_h^N$ and $\mathcal{W}_h^{N+1}$ is not needed for total energy conservation, but it is instead needed for the exact preservation of the involutions, in this case the curl-free property of the velocity field.} 

\bigskip

From a practical point of view we now insert the ansatz \eqref{eqn.acoustics.ansatz} into the relations above and obtain 
\begin{equation}	
	\sum_j\int \limits_{\Omega} \phi_i \phi_j  d \mathbf{x} \, \left({\mathbf{v}}_j^{n+1}-{\mathbf{v}}_j^{n}\right)  
	+ \Delta t \, \sum_j\bigg ( \int \limits_{\Omega} \phi_i \nabla \psi_j d \mathbf{x}\bigg ) \, {p}_j^{n+\frac{1}{2}} = 0,  
	\label{eqn.dg.v.disc} 
\end{equation} 
\begin{equation}	
	\sum_j\int \limits_{\Omega} {\psi}_i \psi_j d \mathbf{x} \, \left( {p}_j^{n+1}-{p}_j^{n} \right)  -   
	\Delta t \, \sum_j\bigg (\int \limits_{\Omega} \nabla {\psi}_i \phi_j d \mathbf{x} \bigg )\cdot {\mathbf{v}}_j^{n+\frac{1}{2}} 
	= 0.  
	\label{eqn.fem.p.disc} 
\end{equation} 
To ease notation we now introduce the following vectors $\mathbf{v}$ and $\mathbf{p}$ that contains the degrees of freedom  of the velocity vector and the pressure, and  matrices: the mass matrix of the DG scheme 
\begin{equation}
   \mathbb D=(D_{ij}), \qquad D_{ij} = \int \limits_{\Omega} \phi_i \phi_j  d \mathbf{x},
    \label{eqn.mass.DG}
\end{equation}
the global mass matrix of the continuous finite element method 
\begin{equation}
   \mathbb M=(M_{ij}), \qquad  M_{ij} = \int \limits_{\Omega} \psi_i \psi_j  d \mathbf{x},
    \label{eqn.mass.FEM}
\end{equation}
and the stiffness tensor
\begin{equation}
   \mathbb K=(K_{ij}), \qquad  {K}_{ij} = \int \limits_{\Omega} \phi_i \nabla \psi_j  d \mathbf{x}\in \mathbb R^d, 1\leq i\leq N_{dof}^p, 1\leq j\leq N_{dof}^v.
\end{equation}
The matrix $\mathbb K$ is a $N_{dof}^p\times N_{dof}^v$ matrix which "coefficients" are vectors of $\mathbb R^d$. It maps the space of pressure onto the space of velocity. The matrix $\mathbb K^T$ is the matrix
$$\mathbb K^T=(K_{ji}^T)_{1\leq j\leq N_{dof}^v, 1\leq i\leq N_{dof}^p}.$$ If $\langle~.~,~.~\rangle_{v}$ is the Euclidian norm for the velocity degrees of freedom,
$$\langle \mathbf{v}, \mathbf{v}\rangle_v=\sum_{i=1}^{N_{Dof}^v}\Vert \mathbf{v}_j\Vert^2$$ and 
$\langle~.~,~.~\rangle_{p}$ is the Euclidian norm for the pressure degrees of freedom
$$\langle \mathbf{p},\mathbf{p}\rangle_{p}=\sum_{i=1}^{N_{dof}^p} p_j^2,$$  we note that
\begin{equation*}
\begin{split}\langle  \mathbb K \mathbf{p}, \mathbf{v} \rangle_v+\langle \mathbf{p}, \mathbb K^T \mathbf{v}\rangle_p&=
\sum_T \sum_{\sigma^v_j\in T}\mathbf{v}_j\cdot\int_T \phi_i \nabla p\, d\mathbf{x}-\sum_{\sigma_j^p} p_j \int_\Omega \nabla\psi_i\cdot \mathbf{v}\, d\mathbf{x}\\
&\quad=
\sum_T \int_T \mathbf{v}\cdot \nabla p\; d\mathbf{x}-\int_\Omega \nabla p \cdot \mathbf{v}\; d\mathbf{x}+
\int_{\partial\Omega} p\mathbf{v}\cdot \n \, d\gamma\\&\qquad=\int_{\partial\Omega} p\mathbf{v}\cdot \n \, d\gamma,
\end{split}
\end{equation*}
i.e. for the $\mathbf{v}\cdot \n=0$ boundary condition, the operators are skew symmetric:
\begin{equation}\label{eq:skew}\langle  \mathbb K \mathbf{p}, \mathbf{v} \rangle_v=-\langle \mathbf{p}, \mathbb K^T \mathbf{v}\rangle_p\end{equation}
because the boundary integrals vanish in that case.
The scheme then becomes 
\begin{equation}	
\mathbb D \big ( \mathbf{v}^{n+1}-\mathbf{v}^n\big ) +\Delta t \mathbb K^T \mathbf{p}^{n+1/2}=0
	\label{eqn.dg.v2} 
\end{equation} 
\begin{equation}	
\mathbb M \big (\mathbf{p}^{n+1}-\mathbf{p}^n\big ) -\Delta t \mathbb K \;\bullet  \mathbf{v}^{n+1/2}=0
	\label{eqn.fem.p2} 
\end{equation} 
where we define 
$$(\mathbb K \bullet \mathbf{v})_i=\sum_{j=1}^{N_{dof}^v} K_{ji}^T\cdot \mathbf{v}_i.$$
Since the mass matrix $D_{ij}$ is based on the discontinuous DG basis and test functions it is block-diagonal and can therefore be easily inverted. Via the choice of the orthogonal Dubiner basis functions on simplex elements \cite{Dubiner} it can be even made globally diagonal. Eqn. \eqref{eqn.dg.v2} therefore can be rewritten as
\begin{equation}	
\mathbf{v}^{n+1}=\mathbf{v}^n -\Delta t \mathbb D^{-1}\mathbb K^T\mathbf{p}^{n+1/2}.
	\label{eqn.dg.v3} 
\end{equation} 
We therefore can easily use the Schur complement to reduce the above system with the unknowns $\mathbf{v}_j^{n+1}$ and $p_j^{n+1}$  to a much simpler system with only the scalar pressure as unknown. More precisely, we formally substitute \eqref{eqn.dg.v3} into \eqref{eqn.fem.p2} and obtain the following discrete pressure wave equation 
\begin{equation} 
\bigg (\mathbb M -\frac{\Delta t^2}{4}\mathbb K\mathbb D^{-1} \mathbb K^T\bigg )\mathbf{p}^{n+1}
=
\mathbb M \mathbf{p}^n +\Delta t\mathbb K \mathbf{v}^n+\frac{\Delta t^2}{4}\mathbb K\mathbb D^{-1}\mathbb K^T \mathbf{p}^n.
    \label{eqn.acoustics.pwave}    
\end{equation} 
The resulting system is clearly symmetric. It is symmetric positive definite in the space of solutions satisfying $\mathbf{v}\cdot \n=0$ since for any $\mathbf{p}$, using \eqref{eq:skew},
$$\langle \mathbf{p}, \bigg (\mathbb M-\frac{\Delta^2}{4} \mathbb K\mathbb D^{-1} \mathbb K^T\bigg )\mathbf{p}\rangle_p=\langle \mathbf{p}, \mathbb M \mathbf{p}\rangle_p
+\frac{\Delta t^2}{4}
\langle \mathbb K^T \mathbf{p}, \mathbb D^{-1} \mathbb K^T\mathbf{p}\rangle_v\geq \langle \mathbf{p}, \mathbb M \mathbf{p}\rangle_p.$$

Hence for its solution we employ a matrix-free conjugate gradient method \cite{cgmethod}. Once the degrees of freedom ${p}_j^{n+1}$ of the new pressure field have been obtained, the degrees of freedom of the discrete velocity field can be easily updated via \eqref{eqn.dg.v2}. 


\subsection{Compatible semi-implicit DG scheme for the vacuum Maxwell equations}

For the vacuum Maxwell equations  \eqref{eqn.maxwell.B} and \eqref{eqn.maxwell.E} we seek the discrete magnetic field in the DG space $\mathcal{U}^N_h$, since the evolution of $\B$ is a mere consequence of the definitions of $\E$ and $\B$ from the underlying vector potential $\A$, while the discrete electric field is defined in the continuous FEM space $\mathcal{W}^{N+1}_h$, since the evolution of $\E$ is given by the Euler-Lagrange equations of an underlying variational principle. Hence, we have:  
\begin{equation} 	
\Bh = \sum_j \phi_j(\mathbf{x}) \mathbf{B}_j  \in \mathcal{U}^N_h, \qquad \qquad 
\Eh = \sum_j \psi_j(\mathbf{x}) \mathbf{E}_j  \in \mathcal{W}^{N+1}_h.
\label{eqn.maxwell.ansatz}
\end{equation} 	
Again, for the sake of simplicity, we first present the method for the semi-discrete case. Multiplication of the induction equation  \eqref{eqn.maxwell.B} by a test function $\phi_i \in \mathcal{U}^N_h$ and integrating over $T_k$, for all $T_k$ the weak formulation of the {induction} equation becomes 
\begin{equation}	
	\int \limits_{T_k} \phi_i \partial_t \Bh \, \dx + 
	\int \limits_{T_k} \phi_i \nabla \times \Eh \, \dx = 0.  
    \label{eqn.maxwell.Bh}
\end{equation} 
Since $\Bh$ is discontinuous across element boundaries, similar to the acoustics case, we need to interpret the curl of the discrete magnetic field in the sense of distributions: for any test function {$\psi_i \in \mathcal{W}^{N+1}_h$}, we write for \eqref{eqn.maxwell.E}
\begin{equation}	
	\int \limits_{\Omega} {\psi}_i \partial_t \Eh \, \dx +  
	\int \limits_{\Omega} \nabla {\psi}_i \times \Bh \, \dx 
    - \int \limits_{\partial \Omega} {\psi}_i \, \mathbf{B} \times \n \, \dx = 0.  
    \label{eqn.maxwell.Eh0}
\end{equation} 

{Again, here $\mathbf{B} \times \n $ is the original boundary condition of the underlying continuous problem on $\partial \Omega$.}
The problem is hence described by the following two global identities, where $\phi_i\in \mathcal{U}^N_h$ and $\psi_i\in \mathcal{W}^{N+1}_h$ are arbitrary: 
\begin{equation}	
	\int \limits_{\Omega} \phi_i \partial_t \Bh \, \dx + 
	\int \limits_{\Omega} \phi_i \nabla \times \Eh \, \dx = 0,  
    \label{eqn.maxwell.Bh.global}
\end{equation} 
and
\begin{equation}	
	\int \limits_{\Omega} {\psi}_i \partial_t \Eh \, \dx +  
	\int \limits_{\Omega} \nabla {\psi}_i \times \Bh \, \dx 
    - \int \limits_{\partial \Omega} {\psi}_i \, \mathbf{B} \times \n \, \dx = 0.  
    \label{eqn.maxwell.Eh.global}
\end{equation} 
It is now easy to prove global total energy conservation. 
\begin{theorem}
    In the case of the boundary condition  $\mathbf{B} \times \n = 0$ on $\partial \Omega$ the scheme \eqref{eqn.maxwell.Bh} and \eqref{eqn.maxwell.Eh.global} conserves global total energy in the sense  
    \begin{equation}
	   \partial_t \int \limits_{\Omega}  \frac{1}{2} \left( \Bh^2 + \Eh^2 \right) \, \dx = 0.   
    \end{equation}    
\end{theorem}
\begin{proof}
    Using $\Bh$ as test function in \eqref{eqn.maxwell.Bh.global} and $\Eh$ as test function in 
    \eqref{eqn.maxwell.Eh.global}, and making use of the boundary condition 
    $\mathbf{B} \times \n = 0$ on $\partial \Omega$ we get: 
    \begin{equation}	
    	\int \limits_{\Omega} \Bh \cdot \partial_t \Bh \, \dx + 
    	\int \limits_{\Omega} \Bh \cdot \nabla \times \Eh \, \dx = 0,  
        \label{eqn.maxwell.Bh.Bh}
    \end{equation}     
    and
    \begin{equation}	
    	\int \limits_{\Omega} \Eh \cdot \partial_t \Eh \, \dx -  
    	\int \limits_{\Omega} \Bh \cdot \nabla \times \Eh  \dx  = 0.  
        \label{eqn.maxwell.Eh.Eh}
    \end{equation}     
    Summing up the above equations yields the sought result 
    \begin{equation}
    	\int \limits_{\Omega} \Bh \cdot \partial_t \Bh \, \dx + 
        \int \limits_{\Omega} \Eh \cdot \partial_t \Eh \, \dx + 
    	\int \limits_{\Omega} \Bh \cdot \nabla \times \Eh \, \dx - 
        \int \limits_{\Omega} \Bh \cdot \nabla \times \Eh  \dx = 
        \partial_t \int \limits_{\Omega} \halb \left( \Bh^2 + \Eh^2 \right) \, \dx = 0.          
    \end{equation}
        \end{proof}

Again we can proceed with a classical Crank-Nicolson discretization in time 
\begin{equation} 
\Bh^{\nph} = \halb \left( \Bh^{n} + \Bh^{n+1} \right), 
\qquad \textnormal{ and } \qquad 
\Eh^{\nph} = \halb \left( \Eh^{n} + \Eh^{n+1} \right), 
\end{equation}
which leads to 
\begin{equation}	
	\int \limits_{\Omega} \phi_i \left( \Bh^{n+1} - \Bh^n \right)  d \mathbf{x}  + 
	\Delta t \int \limits_{\Omega} \phi_i \nabla \times \Eh^{\nph} d \mathbf{x} = 0  
	\label{eqn.dg.B} 
\end{equation} 
and 
\begin{equation}	
	\int \limits_{\Omega} {\psi}_i \left( \Eh^{n+1} - \Eh^n \right) d \mathbf{x} +  
	\Delta t \int \limits_{\Omega} \nabla {\psi}_i \times \Bh^{\nph} d \mathbf{x} 
	= 0  
	\label{eqn.fem.E} 
\end{equation} 
Global total energy conservation of the fully discrete scheme is again completely obvious and analogous to the semi-discrete case when using $\Bh^{\nph}$ and $\Bh^{\nph}$ as test functions.  
With the ansatz \eqref{eqn.maxwell.ansatz} we obtain 
\begin{equation}	
	\sum_j\int \limits_{\Omega} \phi_i \phi_j  d \mathbf{x} \, \left( {\mathbf{B}}_j^{n+1}-{\mathbf{B}}_j^{n} \right)  + 
	\sum_j\Delta t \int \limits_{\Omega} \phi_i \nabla \psi_j d \mathbf{x} \times {\mathbf{E}}_j^{n+\frac{1}{2}} = 0  
	\label{eqn.dg.B1} 
\end{equation} 
and 
\begin{equation}	
	\sum_j\int \limits_{\Omega} {\psi}_i \psi_j d \mathbf{x} \, \left({\mathbf{E}}_j^{n+1}-\tilde{\mathbf{E}}_j^{n} \right) +  
	\sum_j\Delta t \int \limits_{\Omega} \nabla {\psi}_i \phi_j d \mathbf{x} \times  {\mathbf{B}}_j^{n+\frac{1}{2}} 
	= 0  
	\label{eqn.fem.E1} 
\end{equation} 
We can introduce the DG mass matrix $\mathbb D$, the FEM mass matrix $\mathbb M$ and the stiffness matrix $\mathbb K$ as before and the matrix $\mathbb P$ defined by:
$$\mathbb P_{ij}=\int \limits_{\Omega} \nabla {\psi}_i \phi_j d \mathbf{x},$$ and rewrite \eqref{eqn.dg.B1} and \eqref{eqn.fem.E1} as
$$\mathbb D\big ( \mathbf{B}^{n+1}-\mathbf{B}^n\big )+\Delta t \mathbb K\times \bigg( \frac{\mathbf{E}^{n+1}+\mathbf{E}^n}{2}\bigg )=0$$
and
$$\mathbb M\big ( \mathbf{E}^{n+1}-\mathbf{E}^n\big )-\Delta t \mathbb P\times \bigg( \frac{\mathbf{B}^{n+1}+\mathbf{B}^n}{2}\bigg )=0.$$
Here we have used a slight abuse of notations, namely if $\mathbb H$ is a matrix with elements $H_{ij}$ are vectors of $\mathbb R^d$, and if $\mathbf{X}$ is a vector with components $X_i$ also  vectors of $\mathbb R^d$, we define
$${\mathbb H}\times \mathbf{X}=\bigg ( \sum_j H_{ij}\times X_j\bigg ).$$
Introducing the scalar products
$$\langle \mathbf{B},\mathbf{B}\rangle_B=\sum_{j=1}^{N_{dof}^B}\Vert \mathbf{B}_j\Vert^2, \qquad
\langle\mathbf{E},\mathbf{E}\rangle_E=\sum_{j=1}^{N_{dof}^E}\Vert \mathbf{E}_j\Vert^2, $$
we see easily that
$$\langle \mathbf{B},\mathbb K\times\mathbf{E}\rangle_B-\langle \mathbf{E},\mathbb P\times\mathbf{B}\rangle_E=\int_{\partial \Omega}\mathbf{E}\cdot (\mathbf{B}\times \n\big )\ d\gamma,$$
so that if $\mathbf{B}\times \n=0$ on $\partial \Omega$, 
\begin{equation}\label{skew:B}
\langle \mathbf{B},\mathbb K\times\mathbf{E}\rangle_B-\langle \mathbf{E},\mathbb P\times\mathbf{B}\rangle_E=0.
\end{equation}
Also here we make use of the Schur complement, since the DG mass matrix $D_{ij}$ is trivial to invert, i.e. we can insert \eqref{eqn.dg.B} into \eqref{eqn.fem.E} to obtain a discrete vector wave equation for the electric field:
\begin{equation} 
\bigg ( \mathbb M +\frac{\Delta t^2}{4} \mathbb P\times \mathbb D^{-1}\mathbb K\times \bigg ) \mathbf{E}^{n+1}=
\mathbb M \mathbf{E}^n-\Delta t \mathbb P\times \mathbf{B}^{n}-\frac{\Delta t^2}{4}\mathbb P\times\big ( \mathbb D^{-1}\mathbb K\times \mathbf{E}^n\big )
\end{equation}
Again, we note  thanks to \eqref{skew:B}, that,
\begin{equation*}
\begin{split}\langle \mathbf{E}, \bigg ( \mathbb M +\frac{\Delta t^2}{4} \mathbb P\times \big (\mathbb D^{-1}\mathbb K\times \big )\bigg ) \mathbf{E}\rangle_E&=
\langle \mathbf{E},\mathbb M\mathbf{E}\rangle_E+\frac{\Delta t^2}{4}\langle \mathbf{E}, \mathbb P\times \big (\mathbb D^{-1}\mathbb K\times\mathbf{E}\big )\rangle_E\\
&=\langle \mathbf{E},\mathbb M\mathbf{E}\rangle_E+\frac{\Delta t^2}{4}\langle \mathbb K\times\mathbf{E},  \big (\mathbb D^{-1}\mathbb K\times\mathbf{E}\big )\rangle_E\\
&\geq \langle \mathbf{E},\mathbb M\mathbf{E}\rangle_E
\end{split}
\end{equation*}
so that the matrix $\mathbb M +\frac{\Delta t^2}{4} \mathbb P\times\big ( \mathbb D^{-1}\mathbb K\times\big )$ is positive definite.
Once the electric field has been obtained, the magnetic field can be calculated from \eqref{eqn.dg.B1}.
{\paragraph{Proper choice of discrete initial data} 
	In order to guarantee compatible discrete initial data for the DG method, the 
	discrete initial magnetic field $\B_h(\x,0) \in \mathcal{U}^N_h$ must be computed as 
	the discrete curl of a discrete	vector potential 
	$\tilde{\A}_h = \psi_q \A_q \in \mathcal{W}^{N+1}_h$ defined in the continuous finite element space.  The degrees of freedom of the initial magnetic field are then simply given by $\B^0_c = \nabla_c^q \times \A_q$, with the usual primary nabla operator $\nabla_c^q$ introduced before.  
	We emphasize again that simple elementwise $L^2$ projection will in general lead to incompatible 
	discrete initial data for the DG method. Concerning the initial data of the discrete electric field, which is defined in the continuous finite element space $\mathcal{W}^{N+1}_h$, standard nodewise evaluation of the initial data is sufficient.}

\subsection{Compatible semi-implicit DG scheme for the Maxwell-GLM system}

After the systems of linear acoustics and the vacuum Maxwell equations we can tackle also more complex PDEs as the one given by the Maxwell-GLM system \eqref{eqn.maxmunz.B} - \eqref{eqn.maxmunz.psi}. It is very interesting to note that the Maxwell-GLM system actually admits \textit{both} involutions on $\B$: 
\begin{itemize}
    \item for Maxwell-compatible initial data 
    \begin{equation}
        p(\x,0)=q(\x,0)=0, \qquad \textnormal{ and } \qquad 
        \nabla \cdot \B(\x,0) = 0, 
        \label{eqn.mm.maxwellic}
    \end{equation} 
    the magnetic field remains divergence-free for all times, $\nabla \cdot \B = 0$, if it was divergence-free initially, while  
    \item for acoustics-compatible initial data 
    \begin{equation}
        \E(\x,0)=0, \quad q(\x,0)=0, \qquad \textnormal{ and } \qquad 
        \nabla \times \B(\x,0) = 0, 
        \label{eqn.mm.acousticsic}
    \end{equation}     
    it remains curl-free for all times, $\nabla \times \B = 0$, since it initially was assumed to be curl-free.
\end{itemize}
 
The new structure-preserving DG schemes presented in this paper will naturally satisfy both involutions, depending on the chosen initial data. Following again the universal guidelines detailed above, the discrete magnetic field $\Bh$ and the discrete cleaning scalar $\qh$ are sought in the DG space $\mathcal{U}^N_h$, since the PDE for $\B$ and $q$ are mere consequences of the definitions from the abstract potentials $\A$ and $Z$, while the quantities $\Eh$ and $\ph$ are sought in the FEM space $\mathcal{W}^{N+1}_h$, since the evolution equations of $\E$ and $p$ are the \textit{Euler-Lagrange} equations of an underlying variational principle. Therefore, we choose     
\begin{equation} 	
\Bh = \sum_j \phi_j(\mathbf{x}) {\mathbf{B}}_j \, \in \mathcal{U}^N_h, \quad  
\ph = \sum_j \psi_j(\mathbf{x}) {p}_j \, \in \mathcal{W}^{N+1}_h, \quad
\Eh = \sum_j \psi_j(\mathbf{x}) {\mathbf{E}}_j \, \in \mathcal{W}^{N+1}_h, \quad  
\qh = \sum_j \phi_j(\mathbf{x}) {q}_j \, \in \mathcal{U}^N_h, 
\label{eqn.maxmunz.ansatz}
\end{equation} 	
The four discrete quantities $\Bh$, $\ph$, $\Eh$ and $\qh$ are again staggered within carefully chosen function spaces that respect the overarching framework of SHTC systems at the discrete level. 
To facilitate the reader, we first illustrate the approach again in the semi-discrete setting, before discussing also the time discretization. 

Multiplication of equation  \eqref{eqn.maxmunz.B} by a test function $\phi_i \in \mathcal{U}^N_h$ and integrating over $T_k$, for all $T_k$ the weak formulation of the magnetic field reads 
\begin{equation}	
	\int \limits_{T_k} \phi_i \partial_t \Bh \, \dx
    + \int \limits_{T_k} \phi_i \nabla \times \Eh \, \dx  
    + \int \limits_{T_k} \phi_i \nabla \ph \, \dx = 0.  
    \label{eqn.maxmunz.Bh}
\end{equation} 
Since the discrete electric field $\Eh \in \mathcal{W}^{N+1}_h$ and the scalar field $\ph \in \mathcal{W}^{N+1}_h$ are both \textit{continuous} across element boundaries, nothing more  needs to be done for eqn. \eqref{eqn.maxmunz.Bh}.

The discrete magnetic field $\Bh \in \mathcal{U}^N_h$ is \textit{discontinuous} across element boundaries and therefore \eqref{eqn.maxmunz.phi} needs to be interpreted in the sense of distributions, i.e.
\begin{equation}
    	\int \limits_{\Omega} {\psi}_i \partial_t \ph \, \dx -  
	\int \limits_{\Omega} \nabla {\psi}_i \cdot \Bh \, \dx + 
     \int \limits_{\partial \Omega} {\psi}_i \B \cdot \n \, d \gamma = 0.  
    \label{eqn.maxmunz.ph.global}
\end{equation}

Likewise, multiplication of equation  \eqref{eqn.maxmunz.psi} by a test function $\phi_i \in \mathcal{U}^N_h$ and integrating over $T_k$, for all $T_k$ the weak formulation of the $q$ field reads 
\begin{equation}	
	\int \limits_{T_k} \phi_i \partial_t \qh \, \dx
    + \int \limits_{T_k} \phi_i \nabla \cdot \Eh \, \dx  
    = 0.  
    \label{eqn.maxmunz.qh}
\end{equation} 
Since the discrete electric field $\Eh \in \mathcal{W}^{N+1}_h$ is \textit{continuous} across element boundaries, nothing more  needs to be done for eqn. \eqref{eqn.maxmunz.qh}.
%
The discrete magnetic field $\Bh \in \mathcal{U}^N_h$ and also $\qh \in \mathcal{U}^N_h$ are \textit{discontinuous} across element boundaries and therefore \eqref{eqn.maxmunz.E} 
needs to be interpreted in the sense of distributions:
\begin{equation}	
	\int \limits_{\Omega} {\psi}_i \partial_t \Eh \, \dx +  
	\int \limits_{\Omega} \nabla {\psi}_i \times \Bh \, \dx - 
    \int \limits_{\Omega} \nabla {\psi}_i \, \qh \, \dx -
	\int \limits_{\partial \Omega}  {\psi}_i  \B \times \n \, d \gamma + 
    \int \limits_{\partial \Omega}  {\psi}_i \, q \n \, d \gamma 
    = 0.  
    \label{eqn.maxmunz.Eh.global}
\end{equation}

In the following we now prove global total energy conservation of the semi-discrete scheme for the boundary conditions $\B \cdot \n = 0$, $\B \times \n = 0$, $q=0$ on $\partial \Omega$: 
\begin{theorem}
    For the boundary conditions $\B \cdot \n = 0$, $\B \times \n = 0$, $q=0$ on $\partial \Omega$ the scheme \eqref{eqn.maxmunz.Bh}, \eqref{eqn.maxmunz.ph.global}, \eqref{eqn.maxmunz.Eh.global}  and \eqref{eqn.maxmunz.qh} conserves global total energy in the sense 
    \begin{equation}	
	   \partial_t \int \limits_{\Omega}  \frac{1}{2} \left( \Bh^2 + \Eh^2 + \ph^2 + \qh^2 \right) \, \dx = 0.   
    \end{equation} 
\end{theorem}
\begin{proof}
    Since $\phi_i \in \mathcal{U}^N_h$ and $\psi_i \in \mathcal{W}^{N+1}_h$ we can also use $\Bh$, $\qh$, $\Eh$ and $\ph$ as test functions in \eqref{eqn.maxmunz.Bh} - \eqref{eqn.maxmunz.qh}, thus obtaining
\begin{equation}
	\int \limits_{\Omega} \Bh \cdot \partial_t \Bh \, \dx
    + \int \limits_{\Omega} \Bh \cdot \nabla \times \Eh \, \dx  
    + \int \limits_{\Omega} \Bh \cdot \nabla \ph \, \dx = 0,  
    \label{eqn.maxmunz.Bh.Bh}    
\end{equation}
\begin{equation}
    	\int \limits_{\Omega} \ph \partial_t \ph \, \dx -  
	\int \limits_{\Omega} \nabla \ph  \cdot \Bh \, \dx = 0,  
    \label{eqn.maxmunz.ph.ph}
\end{equation}    
\begin{equation}
	\int \limits_{\Omega} \Eh \cdot \partial_t \Eh \, \dx
    - \int \limits_{\Omega} \Bh \cdot \nabla \times \Eh \, \dx  
    - \int \limits_{\Omega} \Eh \cdot \nabla \qh \, \dx = 0,  
    \label{eqn.maxmunz.Eh.Eh}    
\end{equation}
\begin{equation}
    	\int \limits_{\Omega} \qh \partial_t \qh \, \dx +  
	\int \limits_{\Omega} \nabla \qh  \cdot \Eh \, \dx = 0,  
    \label{eqn.maxmunz.qh.qh}
\end{equation}   
where we have also used the boundary conditions $\B \cdot \n = 0$, $\B \times \n = 0$, $q=0$ on $\partial \Omega$. 
Summing up the above equations \eqref{eqn.acoustics.vh.vh} and \eqref{eqn.acoustics.ph.ph} yields
\begin{equation} 
\int \limits_{\Omega} \Bh \cdot \partial_t \Bh \, \dx 
+ \int \limits_{\Omega} \ph \partial_t \ph \, \dx
+ \int \limits_{\Omega} \Eh \cdot \partial_t \Eh \, \dx 
+ \int \limits_{\Omega} \qh  \partial_t \qh \, \dx =  
\partial_t \int \limits_{\Omega}  \frac{1}{2} \left( \Bh^2 + \Eh^2 + \ph^2 + \qh^2 \right) \, \dx = 0. 
\end{equation}
\end{proof}
Using a classical Crank-Nicholson time discretization with the usual notation 
\begin{equation} 
\Bh^{\nph} = \halb \left( \Bh^{n} + \Bh^{n+1} \right), 
\quad 
\ph^{\nph} = \halb \left( \ph^{n} + \ph^{n+1} \right), 
\quad
\Eh^{\nph} = \halb \left( \Eh^{n} + \Eh^{n+1} \right), 
\quad 
\qh^{\nph} = \halb \left( \qh^{n} + \qh^{n+1} \right) 
\end{equation}
and assuming again $\mathbf{B} \cdot \n = 0$, $\B \times \n = 0$, $q=0$ on $\partial \Omega$, let us show that we also get a discrete energy equality. Indeed, 
\begin{equation}
	\int \limits_{\Omega} \phi_i \left( \Bh^{n+1} - \Bh^n \right) \, \dx + 
	\Delta t \int \limits_{\Omega} \phi_i \nabla \times \Eh^{\nph} \, \dx +
    \Delta t \int \limits_{\Omega} \phi_i \nabla \ph^{\nph} \, \dx = 0,  
    \label{eqn.maxmunz.Bh.fd}    
\end{equation}
\begin{equation}
    	\int \limits_{\Omega} {\psi}_i \left( \ph^{n+1} - \ph^n \right) \, \dx -  
	\Delta t \int \limits_{\Omega} \nabla {\psi}_i \cdot \Bh^{\nph} \, \dx = 0,  
    \label{eqn.maxmunz.ph.fd}
\end{equation}
\begin{equation}
	\int \limits_{\Omega} \psi_i \left( \Eh^{n+1} - \Eh^n \right) \, \dx + 
	\Delta t \int \limits_{\Omega} \nabla \psi_i  \times \Bh^{\nph} \, \dx -
    \Delta t \int \limits_{\Omega} \nabla \psi_i \qh^{\nph} \, \dx = 0,  
    \label{eqn.maxmunz.Eh.fd}    
\end{equation}
\begin{equation}
    	\int \limits_{\Omega} {\phi}_i \left( \qh^{n+1} - \qh^n \right) \, \dx +  
	\Delta t \int \limits_{\Omega} {\phi}_i  \nabla \cdot \Eh^{\nph} \, \dx = 0,  
    \label{eqn.maxmunz.qh.fd}
\end{equation}
Multiplying \eqref{eqn.maxmunz.Bh.fd} by ${\mathbf{B}}_j^{n+\frac{1}{2}}$, \eqref{eqn.maxmunz.ph.fd} by $ {p}_j^{n+\frac{1}{2}}$, 
\eqref{eqn.maxmunz.Eh.fd} by ${\mathbf{E}}_j^{n+\frac{1}{2}}$ and \eqref{eqn.fem.p} by $ {q}_j^{n+\frac{1}{2}}$
and summing up over all degrees of freedom and all equations, the spatial terms disappear as before. We then obtain 
\begin{eqnarray}
\int \limits_{\Omega} 
\left( 
\Bh^{n+\frac{1}{2}}\cdot \big ( \Bh^{n+1}-\Bh^{n}\big ) +
\ph^{n+\frac{1}{2}} \big ( \ph^{n+1}-\ph^{n}\big )
+
\Eh^{n+\frac{1}{2}}\cdot \big ( \Eh^{n+1}-\Eh^{n}\big ) +
\qh^{n+\frac{1}{2}} \big ( \qh^{n+1}-\qh^{n}\big ) \right) \, \dx 
= && \nonumber \\ 
\int \limits_{\Omega} 
\left( 
\frac{1}{2} \bigg( \big (\Bh^{n+1}\big )^2+\big ( \ph^{n+1} \big )^2 + \big (\Eh^{n+1}\big )^2+\big (\qh^{n+1}\big )^2 \bigg) - 
\frac{1}{2} \bigg( \big (\Bh^{n}\big )^2+\big ( \ph^{n} \big )^2 + \big (\Eh^{n}\big )^2+\big (\qh^{n}\big )^2 \bigg)
\right) \, \dx = 0, 
&& 
\end{eqnarray} 
with is the fully-discrete conservation of total energy.

\bigskip

From a practical point of view we can again make use of the Schur complement, since the global mass matrix of the DG scheme is block diagonal (or even diagonal for orthogonal basis functions). Inserting \eqref{eqn.maxmunz.Bh.fd} and \eqref{eqn.maxmunz.qh.fd} into \eqref{eqn.maxmunz.Eh.fd} and \eqref{eqn.maxmunz.ph.fd} and making use of the discrete vector calculus identities we obtain one discrete scalar wave equation for the field $\ph$ and one vector wave equation for the field $\Eh$:
\begin{equation} 
\bigg (\mathbb M + \frac{\Delta t^2}{4}\mathbb K\mathbb D^{-1} \mathbb K^T\bigg )\mathbf{p}^{n+1}
=
\mathbb M \mathbf{p}^n 
+\Delta t\mathbb K \mathbf{B}^n
-\frac{\Delta t^2}{4}\mathbb K\mathbb D^{-1}\mathbb K^T \mathbf{p}^n
    \label{eqn.maxmunz.pwave}    
\end{equation} 
\begin{equation} 
 \mathbb M \, \mathbf{E}^{n+1} 
 - \frac{\Delta t^2}{4} \mathbb K \times \mathbb D^{-1}\mathbb K^T \times \mathbf{E}^{n+1}
 + \frac{\Delta t^2}{4} \mathbb K \mathbb D^{-1}\mathbb K^T \cdot \mathbf{E}^{n+1}
 = \mathbf{b}^n 
\label{eqn.maxmunz.Ewave}   
\end{equation}
with the right hand side vector
\begin{equation}
\mathbf{b}^n  = 
    \mathbb M \mathbf{E}^n - \Delta t \mathbb K \times \mathbf{B}^{n}
    +\frac{\Delta t^2}{4}\mathbb K \times \mathbb D^{-1}\mathbb K^T \times \mathbf{E}^n 
    + \Delta t \mathbb K \mathbf{q}^{n}
    - \frac{\Delta t^2}{4} \mathbb K \mathbb D^{-1}\mathbb K^T \cdot \mathbf{E}^{n}
\end{equation}
Once the new electric field $\Eh^{n+1}$ and the new scalar $\ph^{n+1}$ are known, the magnetic field $\Bh$ and the scalar $\qh$ can be readily updated via \eqref{eqn.maxmunz.Bh.fd}
and \eqref{eqn.maxmunz.qh.fd}. 
{\paragraph{Proper choice of discrete initial data} 
	Depending on the test problem, for the Maxwell-GLM system the 
	discrete initial magnetic field $\B_h(\x,0) \in \mathcal{U}^N_h$ must be either computed as 
	the discrete curl of a discrete	vector potential 
	$\tilde{\A}_h = \psi_q \A_q \in \mathcal{W}^{N+1}_h$ in the case of Maxwell-compatible initial data (divergence-free magnetic field), or as the discrete gradient of a scalar potential 
	$\tilde{Z}_h = \psi_q Z_q \in \mathcal{W}^{N+1}_h$ 
	in the case of acoustics-compatible initial data (curl-free magnetic field). The degrees of freedom of the initial magnetic field are then simply given by $\B^0_c = \nabla_c^q \times \A_q$ for Maxwell-compatible initial data and as $\B^0_c = \nabla_c^q  Z_q$ for acoustics-compatible initial data.}

\subsection{Application to the incompressible Euler equations via a simple projection method}
The divergence-free constraint in the incompressible Euler equations is \textit{not} an involution of the PDE but it is rather the remainder of the pressure equation of the compressible Euler equations in the limit when the Mach number tends to zero, see e.g. \cite{KlaMaj,KlaMaj82,KleinMach}. Nevertheless, we try to apply our new scheme also to the set of PDE \eqref{eqn.euler.v} and \eqref{eqn.euler.p} {by making use of a classical projection approach, see e.g. \cite{HW65,PS72,Pat80,Kan86,Guer06,Tavelli2015,Tavelli2016,busto2018projection} and related work on semi-implicit pressure-based schemes on staggered meshes \cite{Casulli1990,Casulli1999,CasulliWalters2000,KleinMach,RoMu99,PM05,DC16,SIGPR}. However, in this paper we employ a suitable staggering of the discrete velocity and pressure fields in the \textit{approximation spaces} $\mathcal{U}_h^N$ and $\mathcal{W}_h^{N+1}$, and \textit{not} via staggering of the discrete velocity and pressure on the \textit{mesh}, as used in classical semi-implicit projection schemes.} The ansatz for the velocity and for the pressure is again given by \eqref{eqn.acoustics.ansatz}, i.e. $\vh \in \mathcal{U}^N_h$ and $\ph \in \mathcal{W}^{N+1}_h$. In the following we assume that the density $\rho$ is a global constant in space and time. 
Since $\vh \in \mathcal{U}^N_h$ is discretized via the DG scheme, we can now simply combine our new compatible discretization of the pressure gradient with a standard DG discretization of the nonlinear convective terms. 
Therefore, multiplication of the velocity equation with a test function $\phi_i \in \mathcal{U}^N_h$ integrating the integral containing the nonlinear convective flux by parts and introducing a numerical flux of the form 
$\hat{\mathbf{f}} = \hat{\mathbf{f}}\left(\vh^-,\vh^+\right) $  together with the same discretization of the pressure gradient as used in acoustics leads to 
\begin{equation}	
	\rho \int \limits_{T_k} \phi_i \partial_t \vh \, \dx - 
    \int \limits_{T_k} \nabla \phi_i \cdot \mathbf{f} \, \dx  + 
    \int \limits_{\partial T_k} \phi_i \, \hat{\mathbf{f}}\left(\vh^-,\vh^+\right) \cdot \n \, d \gamma  
    + 
	\int \limits_{T_k} \phi_i \nabla \ph \, \dx = 0.  
    \label{eqn.euler.vh}
\end{equation} 
In analogy to the acoustics system, also in the case of the incompressible Euler equations the discrete pressure equation \eqref{eqn.euler.p} is interpreted in the sense of distributions:
\begin{equation}
\int_\Omega\nabla \psi_i\cdot \mathbf{v}\; d\mathbf{x}=\int_{\partial\Omega}\psi_i\mathbf{v}\cdot\n\; d\gamma
\label{eqn.euler.ph}
\end{equation}
%
%
%
As usual in the context of projection schemes, we first compute an intermediate velocity field denoted by $\vh^*$, which is not yet divergence-free and which only takes into account the nonlinear convective terms: 

\begin{equation}	
	\rho \int \limits_{T_k} \phi_i \vh^* \, \dx = \rho \int \limits_{T_k} \phi_i \vh^n \, \dx 
    +  
    \Delta t \int \limits_{T_k} \nabla \phi_i \cdot \mathbf{f} \, \dx  - 
    \Delta t \int \limits_{\partial T_k} \phi_i \, \hat{\mathbf{f}}\left(\vh^-,\vh^+\right) \cdot \n \, d \gamma,
    \label{eqn.vstar}
\end{equation} 
where we use a simple numerical flux of Ducros type: 
\begin{equation}
 \hat{\mathbf{f}}\left(\vh^-,\vh^+\right) \cdot \n = \rho v_n  \halb \left( \vh^+ + \vh^- \right) -  \halb \rho |v_n|  \left( \vh^+ - \vh^- \right), 
\end{equation}
with the averaged normal velocity
\begin{equation}
    v_n = \halb \left( \vh^+ + \vh^- \right) \cdot \n 
\end{equation}
After summing over all elements $T_k$ the remaining terms in the momentum equation then become
\begin{equation}
    \rho \int \limits_{\Omega} \phi_i \vh^{n+1} \, \dx =  \rho \int \limits_{\Omega} \phi_i \vh^* \, \dx - \Delta t  \int \limits_{\Omega} \phi_i \nabla \ph^{n+1} \, \dx.  
\label{eqn.euler.vupdate}
\end{equation}
Assuming impermeable wall boundary conditions $\mathbf{v} \cdot \n = 0$ on $\partial \Omega$, the discrete divergence-free condition for the velocity field \eqref{eqn.euler.ph} reads  
\begin{equation}	
	\int \limits_{\Omega}  \nabla {\psi}_i \cdot \vh^{n+1} \, d\mathbf{x} = 0.  
    \label{eqn.euler.divv}
\end{equation} 
As in the previous cases, we simply substitute the new velocity $\vh^{n+1}$ into the discrete divergence-free condition \eqref{eqn.euler.divv} to obtain the discrete pressure Poisson equation for $\ph^{n+1}$, 
and which is given by  
\begin{equation}
-\rho \Delta t^2 \mathbb K\mathbb D \mathbb K^T \mathbf{p}^{n+1} = 
- \rho \Delta t\mathbb K \mathbf{v}^*
    \label{eqn.euler.pressure}
\end{equation}
Once the new pressure is known from \eqref{eqn.euler.pressure}, the new velocity field can be obtained via \eqref{eqn.euler.vupdate}. 
{ 
From our chosen semi-implicit discretization it is obvious that the following discrete energy inequality holds:  
\begin{proposition}
	The semi-implicit scheme \eqref{eqn.euler.vupdate} and \eqref{eqn.euler.divv} is energy stable in the sense 
\begin{equation}
	\int \limits_{\Omega} \halb \rho \left( \vh^{n+1} \right)^2 \, \dx \leq  
	\int \limits_{\Omega} \halb \rho \left( \vh^{*} \right)^2 \, \dx.
\end{equation}	
\end{proposition} 
\begin{proof}
	Using $\vh^{n+1}$ as test function in \eqref{eqn.euler.vupdate} and $-\Delta t \, \ph^{n+1}$ as test function in \eqref{eqn.euler.divv} and summing up all equations leads to 
\begin{equation}
	\rho \int \limits_{\Omega} \vh^{n+1} \cdot \vh^{n+1} \, \dx -  \rho \int \limits_{\Omega} \vh^{n+1} \cdot \vh^* \, \dx + \Delta t  \int \limits_{\Omega} \vh^{n+1} \cdot \nabla \ph^{n+1} \, \dx - \Delta t \int \limits_{\Omega}  \nabla \ph^{n+1} \cdot \vh^{n+1} \, d\mathbf{x} = 0.  
\label{eqn.euler.diss}
\end{equation}	
It is obvious that the last two terms in \eqref{eqn.euler.diss} cancel. Then, splitting the first term in two equal contributions and adding and subtracting $\halb \rho \int \limits_{\Omega} \vh^{*} \cdot \vh^{*} \, \dx $  one arrives at 
\begin{equation}
	\halb \rho \int \limits_{\Omega} \vh^{n+1} \cdot \vh^{n+1} \, \dx + 
	\halb \rho \int \limits_{\Omega} \vh^{n+1} \cdot \vh^{n+1} \, \dx - 
	\halb \rho \int \limits_{\Omega} 2 \vh^{n+1} \cdot \vh^* \, \dx +
	\halb \rho \int \limits_{\Omega} \vh^{*} \cdot \vh^{*} \, \dx  - 
	\halb \rho \int \limits_{\Omega} \vh^{*} \cdot \vh^{*} \, \dx
    = 0. 
	\label{eqn.euler.diss2}
\end{equation}	
Collecting the square contained in the second, third and fourth expression and rearranging terms yields 
\begin{equation}
	 \int \limits_{\Omega} \halb \rho \vh^{n+1} \cdot \vh^{n+1} \, \dx - 
	 \int \limits_{\Omega} \halb \rho \vh^{*} \cdot \vh^{*} \, \dx
	 + 
	 \int \limits_{\Omega} \halb \rho \left( \vh^{n+1} - \vh^{*} \right)^2 \dx 
	= 0. 
	\label{eqn.euler.diss3}
\end{equation}	
Since obviously $\int \limits_{\Omega} \halb \rho \left( \vh^{n+1} - \vh^{*} \right)^2 \dx \geq 0$ we obtain the sought result 
\begin{equation}
	\int \limits_{\Omega} \halb \rho \left( \vh^{n+1} \right)^2 \, \dx \leq  
	\int \limits_{\Omega} \halb \rho \left( \vh^{*} \right)^2 \, \dx.
\end{equation}
\end{proof}
}

For a very recent semi-implicit and asymptotic preserving  discretization of the compressible and incompressible Euler equations based on compatible finite element spaces, the reader is referred to \cite{Zampa2}. {More research on the incompressible Euler and Navier-Stokes equations using the new schemes introduced in this paper will be left to future work. Here, we only want to provide a simple proof of concept, since going into more detail is clearly out of scope of the present paper.}
{\paragraph{Proper choice of discrete initial data} 
	For the incompressible Euler system, to guarantee an initially divergence-free velocity field the 
	discrete velocity field $\v_h \in \mathcal{U}^N_h$ must be initialized as the discrete curl of a discrete	vector potential  $\tilde{\A}_h = \psi_q \A_q \in \mathcal{W}^{N+1}_h$. The degrees of freedom of the initial velocity field are simply given by $\v^0_c = \nabla_c^q \times \A_q$.}



\section{Analysis of the compatibility conditions}
\label{sec.comp} 



In this section, we show that the scheme for the equations of linear acoustics, as well as the one for the vacuum Maxwell equations conserve the curl (resp. the divergence) of the velocity (resp. the magnetic field).
We do it for the semi discrete scheme. It is obvious that the same proof also applies to the fully discrete scheme with Crank-Nicolson time discretization.  
{This section is complementary to Section \ref{sec.disc.nabla}, where a different notation was used and where the main properties of the two discrete differential operators introduced in this paper have also been studied.}

\subsection{Discrete curl}
We are given $\mathbf{w}$ such that on each element $T$,
$$\mathbf{w}=\sum_{\sigma_T\in T} \mathbf{w}_{\sigma_T}\varphi_{\sigma_T}$$
and we would like to define a curl operator. We note that $\varphi_{\sigma_T}\in \mathbb P_N(T)$, and we consider the finite element space spaned by the $\{\psi_{\sigma_p}\}$ where the $\sigma_p$ are the Lagrange points on the elements, 
$\psi_{\sigma_p}\in \mathbb P_{N+1}(T)$ and is globally continuous. We will assume that $\psi$ vanishes on $\partial \Omega$, i.e. we consider only the internal degrees of freedom (none on the boundary).
If $\mathbf{w}$ were regular enough, we would have
$$\int_\Omega \psi_{\sigma_p}\text{curl }\mathbf{w}\; \dx = \textcolor{red}{-} \int_\Omega \nabla\psi_{\sigma_p}\times \mathbf{w}\; \dx.$$
{ Since
	$$\int_\Omega \psi_{\sigma_p}\text{curl }\mathbf{w}\; \dx=\sum_{T}\int_T \psi_{\sigma_p}\text{curl }\mathbf{w}\; \dx=
	\sum_T \bigg ( \int_T\text{curl }\big ( \psi_{\sigma_p}\mathbf{w}\big ) \; \dx-\int_T \nabla\psi_{\sigma_p}\times \mathbf{w}\;\dx\bigg )$$
	this is equal to 
	$ -\int_\Omega \nabla\psi_{\sigma_p}\times \mathbf{w}\; \dx$ if $$\sum_T \int_T\text{curl }\big ( \psi_{\sigma_p}\mathbf{w}\big ) \; \dx=0.$$If $\mathbf{w}$ is regular enough, in 2D we have 
	$$\sum_T \int_{\partial T} \psi_{\sigma_p}\mathbf{w}\cdot d\btau =0$$ where $d\btau$ is the oriented tangent  vector on $T$. Since $\psi_{\sigma_p}$ and $\mathbf{w}$ are continuous across the boundary of $T$,  this term cancels.
	\begin{remark}[About 3D]
		In 3D, though it is not the purpose of the paper, we have for a regular enough vector field $\mathbf{f}$ and a simplex $K$
		$$\int_K\text{curl }\mathbf{f}\;\dx=\int_{\partial K}\mathbf{f}\times \mathbf{n} \, dS$$ and we would have the same relation by assuming the continuity of $\mathbf{f}$ across the faces of the simplex.
\end{remark}}

This leads to the following natural definition of the $\CURL$ operator that is defined at the $\sigma_p$, so this will defined an object in $\mathcal {W}_h$.

{
	\begin{definition}[Discrete curl of an element of $\mathcal{U}_h^N$]
		Let $\Omega$ be an open set of $\R^d$ with $d=2,3$. Let  $\mathbf{w}\in \mathcal{U}_h^N$.  Its $\mathcal{W}_h^{N+1}$ curl defined at the $\mathcal{W}_h^{N+1}$ Lagrange points $\sigma_p$ by
		\begin{equation}
			\label{CURL}
			\mathbb M\begin{pmatrix} \CURL_{\sigma_p}\mathbf{w}\end{pmatrix}_{\sigma_p}= \textcolor{red}{-} \begin{pmatrix}
				\int_\Omega \nabla\psi_{\sigma_p}\times \mathbf{w}\; \dx.\end{pmatrix}_{\sigma_p}
		\end{equation} where $\mathbb M$ is the global mass matrix of the continuous FEM.
	\end{definition}

	Let us show that  for the acoustic problem, the discrete curl of the velocity does not evolve in time.
	\begin{proposition}\label{prop:curl}
		Let be $(\mathbf{v}_h, \tilde{p}_h)\in \mathcal{U}_h^N\times \mathcal{W}_h^{N+1}$ be a solution of the semi discrete scheme \eqref{eqn.acoustics.v}-\eqref{eqn.acoustics.p} with the perioic boundary conditions or the following ones 
		$$\nabla \tilde{p}_h\cdot \btau=0 \text{ in } 2D \text{ or }\nabla\tilde{p}_h\times \mathbf{n}=0 \text{ in }3D$$
		where $\btau$ is a tangent vector to $\partial \Omega$ in 2D and $\mathbf{n}$ is a normal to $\partial \Omega$ in 3D.
		Then, for all $\sigma_p$, we have
		$$\dfrac{d}{dt}\mathbf{CURL}_{\sigma_p}\mathbf{v}_h=0.$$
	\end{proposition}
	\begin{proof}
		If $\psi_{p}\in \mathbb P_{N+1}$, then ${\nabla \psi_p}_{\vert T}\in \big (\mathbb P_{N}\big )^d$. so we can use this as test function and we can write
		\begin{equation*}
			\int_\Omega \nabla \psi_p\times \dfrac{d\vh}{dt}\; \dx=\sum_{T}\int_T\nabla \psi_p\times \dfrac{d\vh}{dt}\; \dx\\ =
			-\sum_T \int_T \nabla \psi_p\times \nabla \tilde{p}_h \; \dx.
		\end{equation*}
		Then we use the identity
		$$\text{curl}\big ( \varphi \mathbf{w}\big )=\varphi\big (\text{curl}\,\mathbf{w}\big )+\nabla\varphi\times \mathbf{w}$$
		to see that
		$$\int_T  \nabla \psi_p\times \nabla \tilde{p}_h \; \dx=\int_T\text{ curl}\big (\psi_p \nabla \tilde{p}_h\big ) \dx -\int_T \psi_p \big ( \text{curl}\nabla \tilde{p}_h) \dx $$
		Since in $T$, $\tilde{p}_h$ and $\psi_p$ are smooth, $\text{curl}\nabla \tilde{p}_h=0$ and then, using Stokes formula and the boundary conditions on $\nabla\tilde{p}_h$, 
		\begin{itemize}
			\item In 2D,
			$$\int_T  \nabla \psi_p\times \nabla \tilde{p}_h \; \dx=\int_{\partial T}\psi_p \nabla \tilde{p}_h\cdot \mathbf{\tau} \;d\gamma.$$
			Let $e$ be any edge of $T$, with $\btau$ the tangent vector that is compatible with the orientation of $T$. The edge $e$ is the edge of another triangle $T'$: here we use that  the mesh is conformal. The tangent vector $\btau'$ is $-\btau$. Since the pressure is  continuous accros $e$, taking into account the boundary condition, the total contribution of $\nabla \psi_p\cdot \mathbf{\tau}$ over all internal edge vanishes by lemma \ref{zob}:
			$$\sum_T \int_T \nabla \psi_p\times \nabla \tilde{p}_h \; \dx=0.$$
			\item In 3D, 
			we get
			$$\int_T  \nabla \psi_p\times \nabla \tilde{p}_h \; \dx=\int_{\partial T}\psi_p \nabla \tilde{p}_h\times \mathbf{n} \;dS.$$
			From lemma \ref{zob}, since $\tilde{p}_h$ is continuous across any face  $f=T_1\cap T_2$, thanks to the conformality of the triangulation, using that the normal to $f$ seen from $T_1$ is the opposite of that seen from $T_2$, using the boundary conditions, we see that 
			$$\sum_T \int_T \nabla \psi_p\times \nabla \tilde{p}_h \; \dx=0.$$
		\end{itemize}
		This shows that $$\dfrac{d}{dt}\CURL_{\sigma_p} \mathbf{v}=0.$$
	\end{proof}
	
	The following lemma is very classical, we provide its proof for the sake of completeness.
	\begin{lemma}\label{zob}
		Let  $g\in \mathcal{W}_h^M$.
		\begin{itemize}
			\item in 2D:  Let $e$ be any  internal edge of the triangulation, and $T_1$, $T_2$ the two triangles that share $e$. Let  $\btau$ is the tangent vector along $e$. Then
			$$\nabla g_{\vert T_1}\cdot \btau=\nabla g_{\vert T_2}\cdot \btau$$ on $e$
			\item In 3D.  Let $f$ be any internal face of the triangulation and $T_1$, $T_2$ the two simplex that share $f$. If $\mathbf{n}$ is a normal of $f$, then
			$$\nabla g_{\vert T_1}\times  \mathbf{n}=\nabla g_{\vert T_2}\cdot \mathbf{n}$$ on $f$.
		\end{itemize}
	\end{lemma}
	\begin{proof} $ $
		\begin{itemize}
			\item 2D case:
			Let $\mathbf{x}$ a point in the interior of $e$, and $h\in \R$ small enough so that $\mathbf{x}+h\btau\in e$.
			Then we have in $T_1$,
			$$g_{\vert T_1}(\mathbf{x}+h\btau)-g_{\vert T_1}(\mathbf{x})=h \nabla g_{\vert T_1}(\mathbf{x})\cdot \btau) +o(h)$$
			i.e.
			$$\nabla g_{\vert T_1}(\mathbf{x})\cdot \btau =\lim\limits_{h\rightarrow 0} \dfrac{g_{\vert T_1}(\mathbf{x}+h\btau)-g_{\vert T_1}(\mathbf{x})}{h}$$
			on one hand, and on the other hand
			$$\nabla g_{\vert T_2}(\mathbf{x})\cdot \btau =\lim\limits_{h\rightarrow 0} \dfrac{g_{\vert T_2}(\mathbf{x}+h\btau)-g_{\vert T_2}(\mathbf{x})}{h}.$$
			But on $e$, by continuity, $g_{\vert T_1}=g_{\vert T_2}$, so that we get the result.
			\item 3D case:
			We write $\nabla g=\dpar{g}{n}\mathbf{n}+\dpar{g}{\tau_1}\btau_1+\dpar{g}{\tau_2}\btau_2$ where $\mathbf{n}\cdot\btau_1=0$ and $\btau_2=\mathbf{n}\times \btau_1$. Then we see that 
			$g(\mathbf{x}+\varepsilon\btau_1)=g(\mathbf{x})+\dpar{g}{\tau_1}(\mathbf{x})\btau_1\varepsilon+o(\varepsilon)$ and since $\btau_1$ depends only on the face, as well as $\btau_2$, and $g$ is continuous at $\mathbf{x}$, then $\dpar{g}{\tau_1}(\mathbf{x})\btau_1$ will be continuous across the face. Since
			$$\nabla g\times \mathbf{n}=\dpar{g}{\tau_1}\btau_1\times \mathbf{n}+\dpar{g}{\tau_2}\btau_2\times \mathbf{n},$$ we get the result.
		\end{itemize}
	\end{proof}
	\begin{proposition}
		Let be $(\mathbf{v}_h, \tilde{p}_h)\in \mathcal{U}_h^N\times \mathcal{W}_h^{N+1}$ be a solution of the fully discrete scheme \eqref{eqn.dg.v}-\eqref{eqn.fem.p} with the same boundary conditions as in proposition \ref{prop:curl}. Then, for all $\sigma_p$, we have
		$$\mathbf{CURL}_{\sigma_p}\mathbf{v}^{n}_h=\mathbf{CURL}_{\sigma_p}\mathbf{v}^{0}_h.$$
	\end{proposition}
	\begin{proof}
		The proof is the same as   that of proposition \ref{prop:curl}.
	\end{proof}
	\begin{remark}
		In the numerical simulations, we have taken periodic boundary conditions.
	\end{remark}
}
\subsection{Discrete divergence}
{
	Using the same notations, and if $\mathbf{w}$ were regular enough, we would have
	$$\int_\Omega \psi_{\sigma_p}\; \text{div }\mathbf{w}\dx=-\int_\Omega \nabla\psi_{\sigma_p}\cdot \mathbf{w}\; \dx.$$
	This leads to the following definition of the divergence
	\begin{definition}
		Let be  $\mathbf{w}\in \mathcal{U}_h^{N}$. We define its divergence in $\mathcal{W}_h^{N+1}$  at the degrees of freedom $\sigma_p$ by
		$$\mathbb M\begin{pmatrix}\mathbf{DIV}_{\sigma_p}\mathbf{w}\end{pmatrix}_{\sigma_p}=
		\begin{pmatrix}\int_\Omega \nabla \psi_{\sigma_p}\cdot \mathbf{w}\dx \end{pmatrix}_{\sigma_p}.$$
	\end{definition}
	This allows to show the following propositions for the schemes on the Maxwell equations.
	\begin{proposition}
		Let $(\mathbf{\tilde E}_h,\mathbf{B}_h)\in \mathcal{W}_h^{N+1}\times \mathcal{U}_h^N$ be a solution of the semi-discrete scheme \eqref{eqn.maxwell.E}-\eqref{eqn.maxwell.B} with the boundary condition
		$\mathbf{E}_h \times \mathbf{n}=0$ on $\partial \Omega$. Then the discrete divergence of the magnetic field is invariant over time.
	\end{proposition}
	\begin{proof}
		If $\psi_p\in \mathbb P_{N+1}$ and is continuous across the faces of the elements $T$, we have because ${\nabla\psi_{\sigma_p}}_{|T}\in \mathcal{U}_h^{N}$,
		$$\int_\Omega \nabla\psi_p\cdot \dfrac{d\mathbf{B}_h}{dt}\dx=\sum_T\int_T\nabla\psi_p\cdot \dfrac{d\mathbf{B}_h}{dt}\dx
		=-\sum_T \int_T \nabla\psi_p\cdot \nabla\times \tilde{\mathbf{E}}_h\dx.$$
		Since for smooth vector fields $$\text{div }\big (\mathbf{w}\times \mathbf{k}\big )=\mathbf{w}\cdot \text{curl }\mathbf{k}-
		\mathbf{k}\cdot \text{curl }\mathbf{w},$$
		we get
		$$\int_T \nabla\psi_p\cdot \nabla\times \tilde{\mathbf{E}}_h\dx=\int_{\partial T} \big (\nabla \psi_p\times \mathbf{E}_h\big )\cdot \mathbf{n} \; d\gamma+\int_T \mathbf{E}_h\cdot \text{ curl}\big ( \nabla \psi_p)\dx=\int_{\partial T} \big (\nabla \psi_p\times \mathbf{E}_h\big )\cdot \mathbf{n} \; d\gamma.$$
		This leads to
		$$\int_\Omega \nabla\psi_p\cdot \dfrac{d\mathbf{B}_h}{dt}\dx=\sum_T \int_{\partial T}\big (\nabla \psi_p\times \mathbf{E}_h\big )\cdot \mathbf{n} \; d\gamma.$$
		The classical remark is that, since
		$$\big (\nabla \psi_p\times \mathbf{E}_h\big )\cdot \mathbf{n}=\det\big ( \nabla \psi_p, \mathbf{E}_h, \mathbf{n}\big )=
		-\mathbf{E}_h\cdot \big (\nabla \psi_p\times \mathbf{n}\big ),$$ since $\psi_p$ is continuous  across the faces of $T$ and $\mathbf{n}$ is normal to the face, that the jump of $\nabla \psi_p\times \mathbf{n}$ is null, see lemma \ref{zob}., and using the continuity of $\mathbf{E}_h$, 
		$$\sum_T \int_{\partial T}\big (\nabla \psi_p\times \mathbf{E}_h\big )\cdot \mathbf{n} \; d\gamma=-
		\int_{\partial \Omega} \mathbf{E}_h \nabla \psi_p\times \mathbf{n}\; d\gamma=\int_{\partial \Omega} \mathbf{E}_h \times \mathbf{n} \cdot \nabla \psi_p\; d\gamma.$$
		This shows that if $\mathbf{E}_h \times \mathbf{n}=0$ on $\partial \Omega$, 
		$$\int_\Omega \nabla\psi_p\cdot \dfrac{d\mathbf{B}_h}{dt}\dx=0$$ so that for all $\sigma_p$,
		$$ \dfrac{d}{dt}\mathbf{DIV}_{\sigma_p}\mathbf{B}_h=0. $$
	\end{proof}
	Using the same technique, we also get the following result on the fully discrete problem:
	\begin{proposition}
		Let $(\mathbf{\tilde E}_h,\mathbf{B}_h)\in \mathcal{W}_h^{N+1}\times \mathcal{U}_h^N$ be a solution of the discrete scheme \eqref{eqn.dg.B}-\eqref{eqn.fem.E} with the boundary condition
		$\mathbf{E}_h \times \mathbf{n}=0$ on $\partial \Omega$. Then the discrete divergence of the magnetic field is invariant over time.
	\end{proposition}
}

\subsection{Case of the Maxwell GLM equations}
Following the same reasoning, we get the following expression for the time derivative of the discrete divergence of $\B_h$
$$
\dfrac{d}{dt}\int_\Omega \nabla \psi_i\cdot \mathbf{B}_h \, \dx + \int_\Omega \nabla\psi_i\cdot \nabla \ph \, \dx=0
$$
and the relation 
$$
\dfrac{d}{dt}\int_\Omega \nabla \psi_i \times \mathbf{B}_h \, \dx + \int_\Omega \nabla\psi_i \times \nabla \Eh \, \dx=0
$$
for the time derivative of the discrete curl of $\Bh$. 
\begin{itemize}
	\item If we now assume Maxwell-compatible initial data $\ph=\qh=0$, see eqn. \eqref{eqn.mm.maxwellic}, we  immediately get that the discrete magnetic field remains divergence-free for all times thanks to the first of the above relations.
	\item In the case of acoustics-compatible initial data ($\Eh = 0$, $\qh=0$), see eqn. \eqref{eqn.mm.acousticsic}, the discrete vector field $\Bh$ remains curl-free for all times thanks to the second of the above relations.
\end{itemize}      

\begin{remark}
	With this formulation, it does not seem to be possible to say anything about the divergence and the curl of the electric field. The same remark also applies to the Maxwell equations. 
\end{remark}

\section{Numerical experiments}\label{sec:experiments}

In this section we carry out a series of numerical experiments to verify the discrete vector calculus identities and total energy conservation also in the computer code. In addition, we
perform some numerical convergence studies to show that the new DG scheme reaches its designed
order of accuracy. For basis polynomials of approximation degree $N$ we expect the scheme to 
be of order $N+1$ in space. The Crank-Nicolson time integrator used in this paper is only second order accurate in time and in the case of time-dependent problems, the time step size will be chosen small enough so that time discretization errors are sufficiently small compared to the spatial discretization error. In all test cases the boundary conditions are periodic in all directions.  
In all numerical experiments, apart from Section 4.1, the curl errors $\epsilon_{c}$ and the divergence errors $\epsilon_d$ of a generic vector field $\A$ are computed as follows: 
\begin{equation}
    \epsilon_c = \max_i  \left( \left| \int \limits_{\Omega} \nabla \psi_i \phi_j \times \A_j \, \dx \right| \right), 
    \qquad \textnormal{ and } \qquad 
    \epsilon_d = \max_i  \left( \left| \int \limits_{\Omega} \nabla \psi_i \phi_j \cdot \A_j \, \dx \right| \right), 
\end{equation}

{
\subsection{Basic sanity checks} 
The only aim of this section is the numerical verification of basic properties stated and proven in Proposition \ref{prop.td} and Proposition \ref{prop.nd} and as also outlined in Section \ref{sec.comp}, i.e. to show that the tangential components of the discrete primary gradient of a scalar field $\Zh \in \mathcal{W}_h^{N+1}$ and the normal component of the discrete primary curl of a vector field $\Ah \in \mathcal{W}_h^{N+1}$ are continuous across elements when projected into the DG space $\mathcal{U}_h^{N}$. We choose a computational domain $\Omega=[-\halb,+\halb]^d$, discretized on 902 triangles with $N_x=N_y=20$ elements along each edge of the computational domain in two space dimensions ($d=2$) and on  
6879 tetrahedra ($N_x=N_y=N_z=10$ elements along each edge) in three space dimensions ($d=3$). To make this test as difficult and as meaningful as possible, we define a \textit{random} discrete scalar field $\Zh \in \mathcal{W}_h^{N+1}$ and a \textit{random} discrete vector field $\Ah \in \mathcal{W}_h^{N+1}$, in the sense that the nodal degrees of freedom of the continuous finite element representation of $\Zh$ and $\Ah$ are chosen as pseudo-random numbers in Fortran in the range $[0,0.001]$, with seed value 1709 in the case of the Maxwell equations and seed value 2007 in the case of the equations of linear acoustics. The choice of a fixed seed makes our results reproducible when using the same computer code on the same CPU with the same compiler. The polynomial approximation degree of the DG scheme is chosen as $N=3$ in all cases and periodic boundary conditions are imposed on $\partial \Omega$. The definition of the unit normal vector on each element boundary is clear. In two space dimensions, the only tangent vector that needs to be considered at element boundaries points from the local point one of each edge to the local point two on each edge. In three space dimensions, the two tangent vectors considered for this test are those connecting the local points 1 and 2 and the local points 1 and 3 on each triangular face of all tetrahedral elements of the mesh. We then measure the maximum absolute value of the error in the tangential components of the gradient of the scalar field $\Zh$ ($\v_h = \nabla \Zh$, with $\v_h \in \mathcal{U}_h^{N}$) and of the normal component of the curl of the vector field $\Ah$ ($\B_h = \nabla \times \Ah$, with $\B_h \in \mathcal{U}_h^{N}$) computed in $(N+2)^{(d-1)}$ Gaussian quadrature points on all boundaries of all elements of the domain. Note that the Gaussian quadrature points used in this test (standard Gauss-Legendre points on the 1D edges in the two-dimensional case and the classical tensor-product formulas given by Stroud in \cite{stroud} for the triangular 2D faces in 3D) on the boundaries do in general not coincide with the nodal degrees of freedom of the FEM space, hence are totally unrelated and thus general enough and representative to compute the errors. The results are documented in Table \ref{tab.sanity} and clearly confirm that also the practical implementation of our scheme in actual computer code satisfies Propositions \ref{prop.td} and \ref{prop.nd}, since the reported errors are of the order of machine precision in all cases.  
\begin{table}  
	{
	\caption{Maximum absolute value of the errors in the tangential components of the gradient $\v_h = \nabla \Zh$ of the random scalar field $\Zh$ and in the normal component of the curl $\B_h = \nabla \times \Ah$ of the random vector field $\Ah$ in two and three space dimensions, respectively. The reported errors are the maximum error found in all Gaussian quadrature points of all element boundaries.} 	
	\begin{center} 
		\renewcommand{\arraystretch}{1.5}
		\begin{tabular}{ccc}  
			                       & two space dimensions & three space dimensions  \\ 
             \hline
			 max. error in the tangential components of $\v_h = \nabla \Zh$   &     1.72778458E-015   &  1.76247905E-015                      \\ 
			 max. error in the normal component of $\B_h = \nabla \times \Ah$ & 	1.49533164E-015   &  2.40085729E-015                     \\ 
			 \hline 
		\end{tabular}
	\end{center}
	\label{tab.sanity}
	}
\end{table}
}

\subsection{Equations of linear acoustics} 

\begin{figure}[!htbp]
	\begin{center}
		\begin{tabular}{cc} 
			\includegraphics[width=0.45\textwidth]{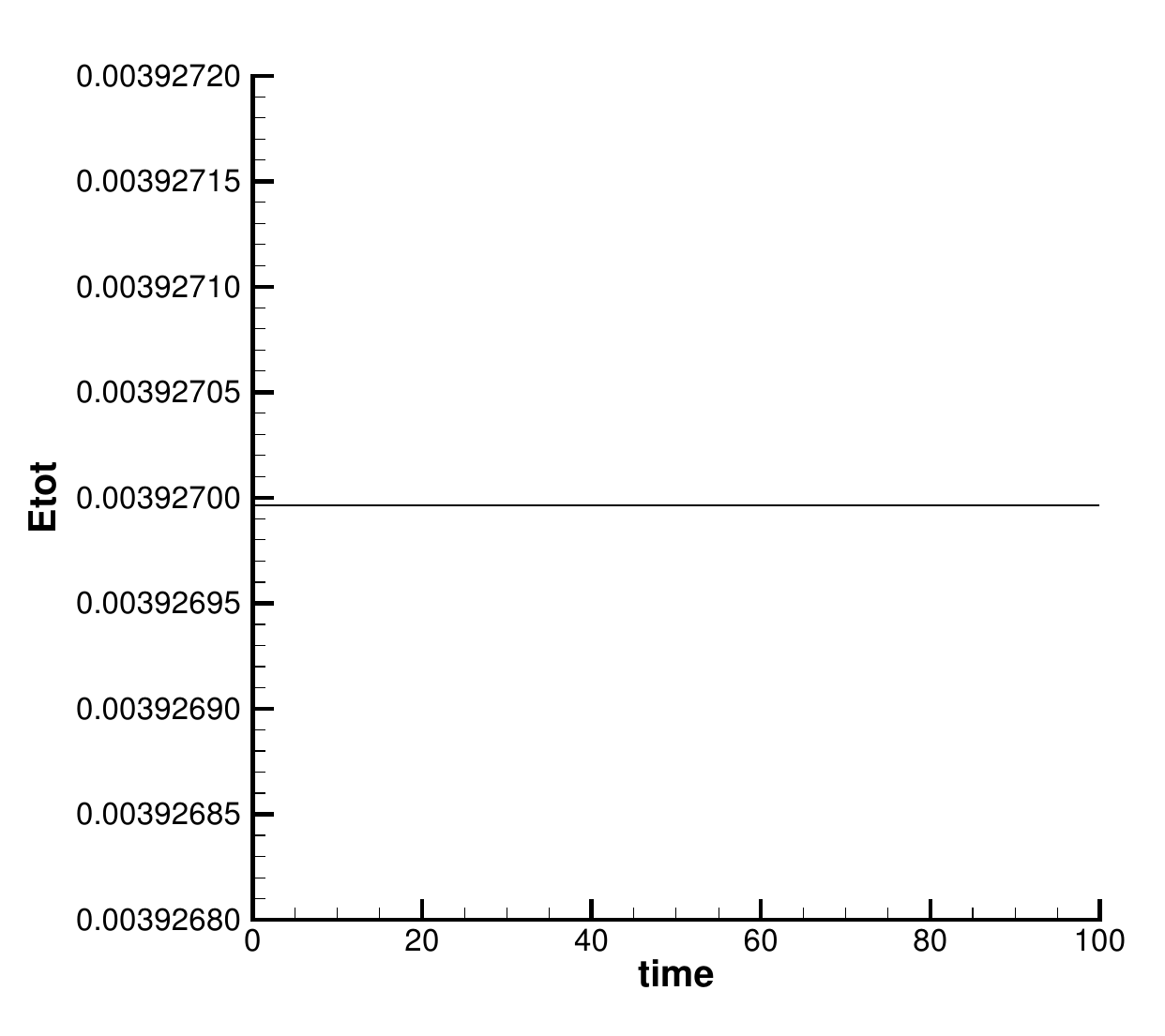}     & 
			\includegraphics[width=0.45\textwidth]{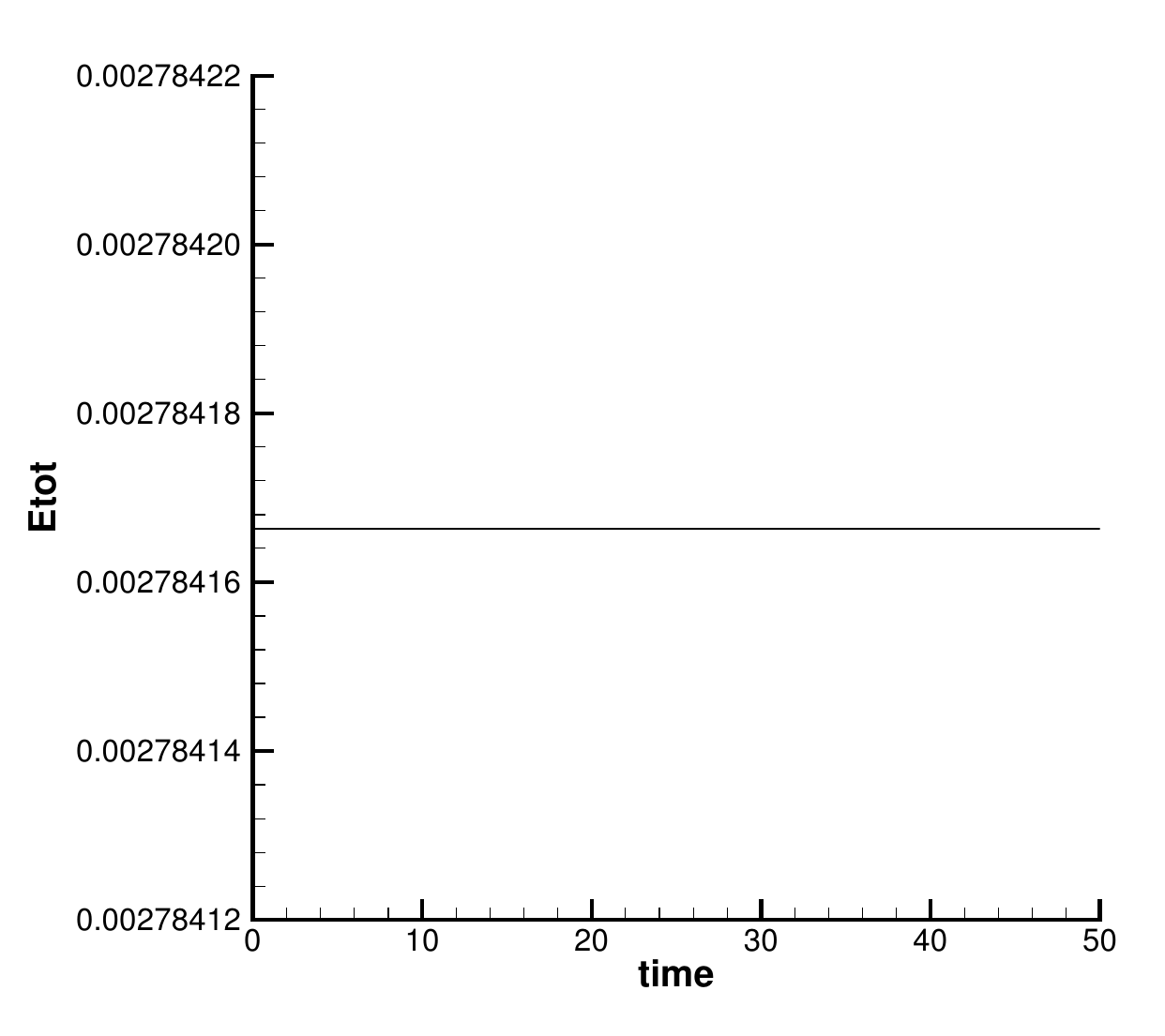}  \\  
			\includegraphics[width=0.45\textwidth]{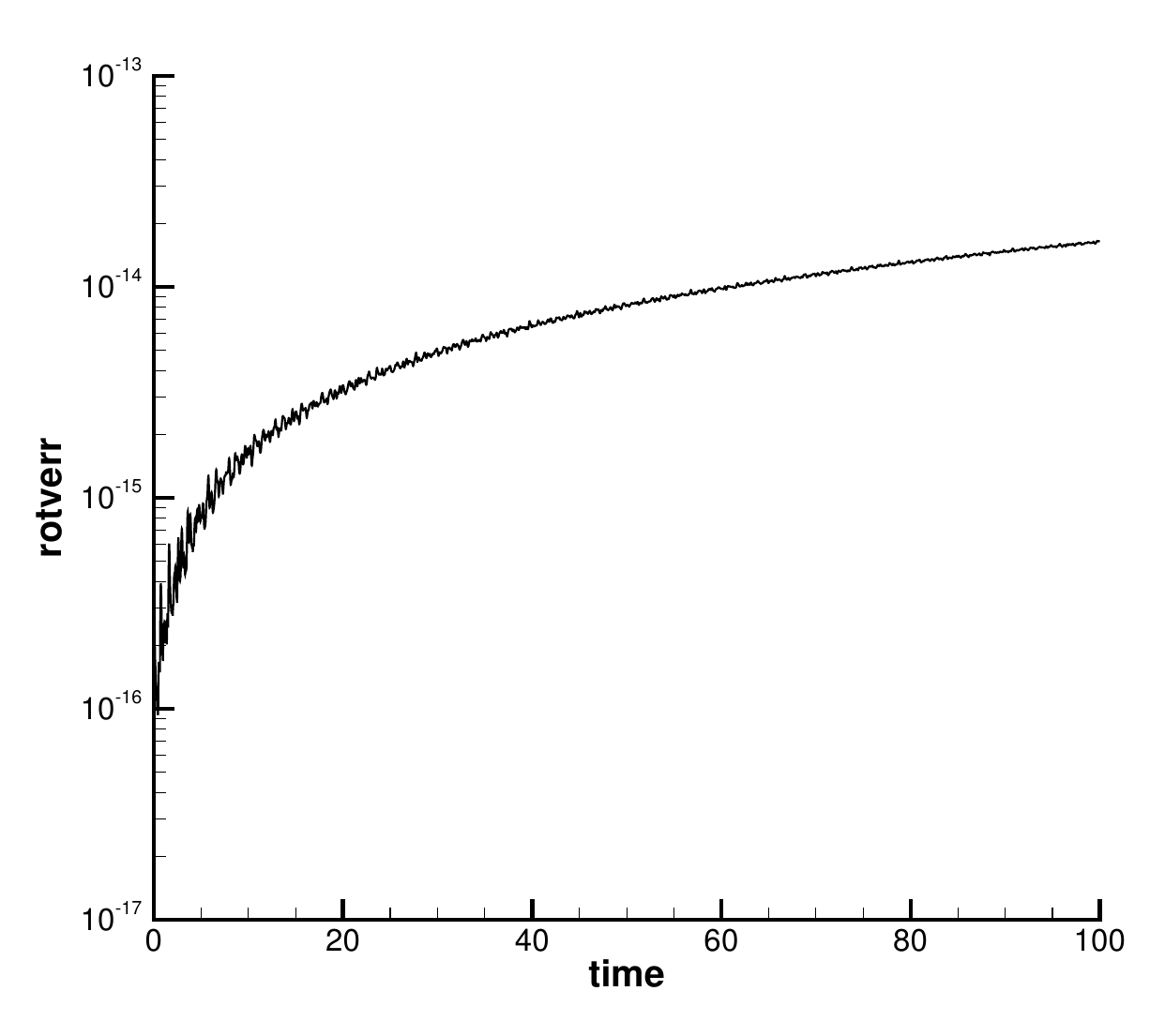}       & 
			\includegraphics[width=0.45\textwidth]{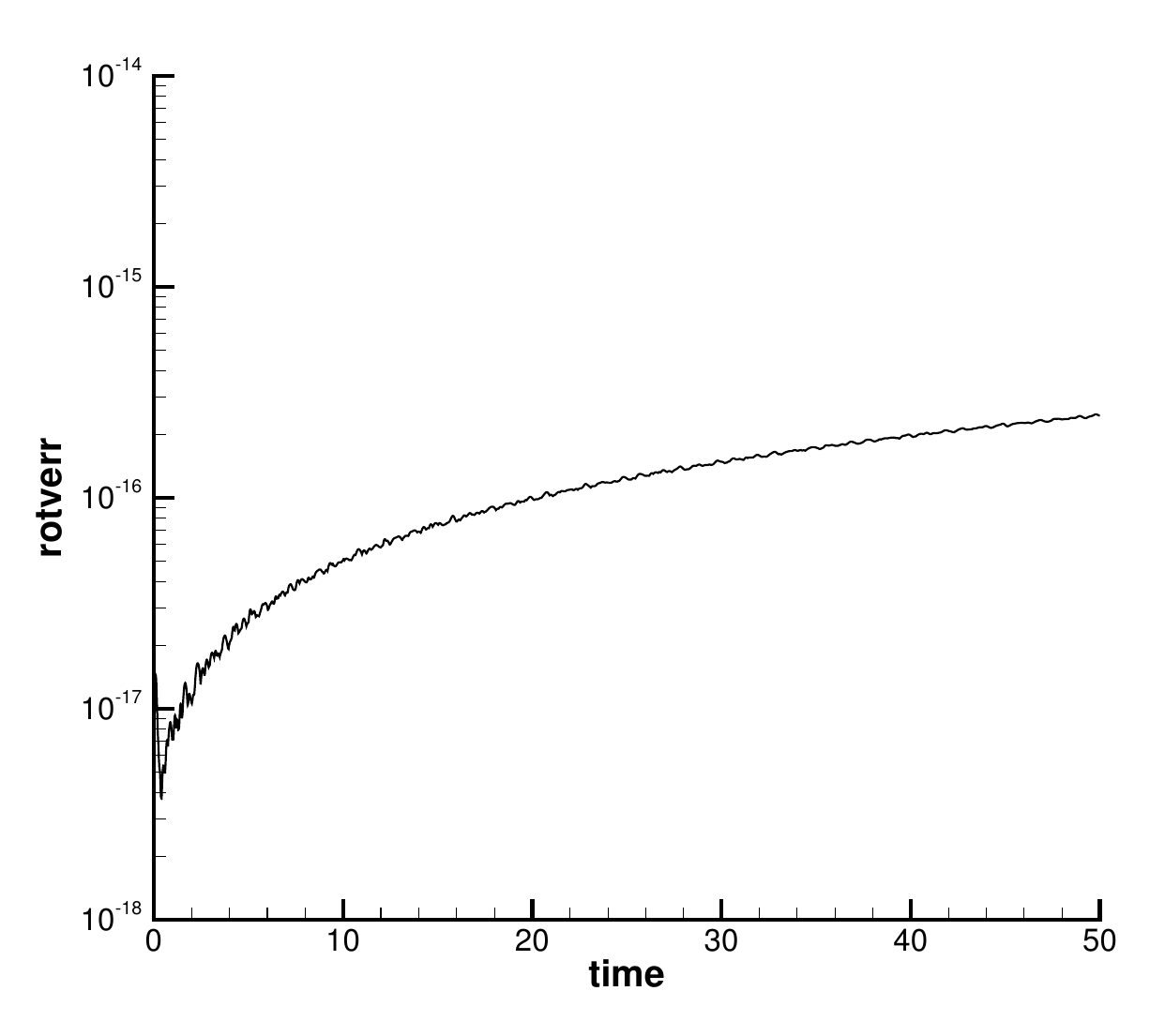}    \\  
			\includegraphics[width=0.45\textwidth]{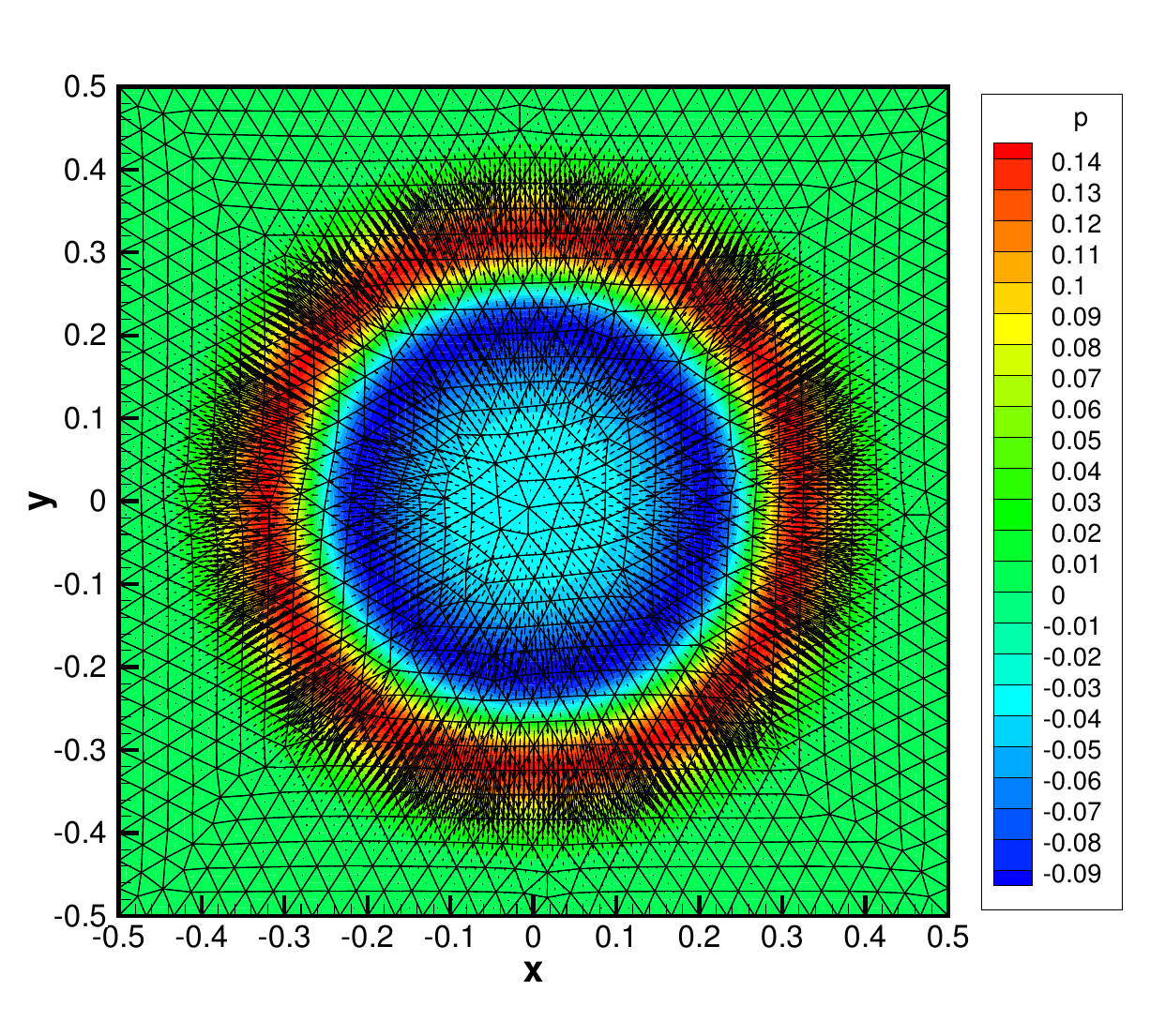}          & 
			\includegraphics[width=0.45\textwidth]{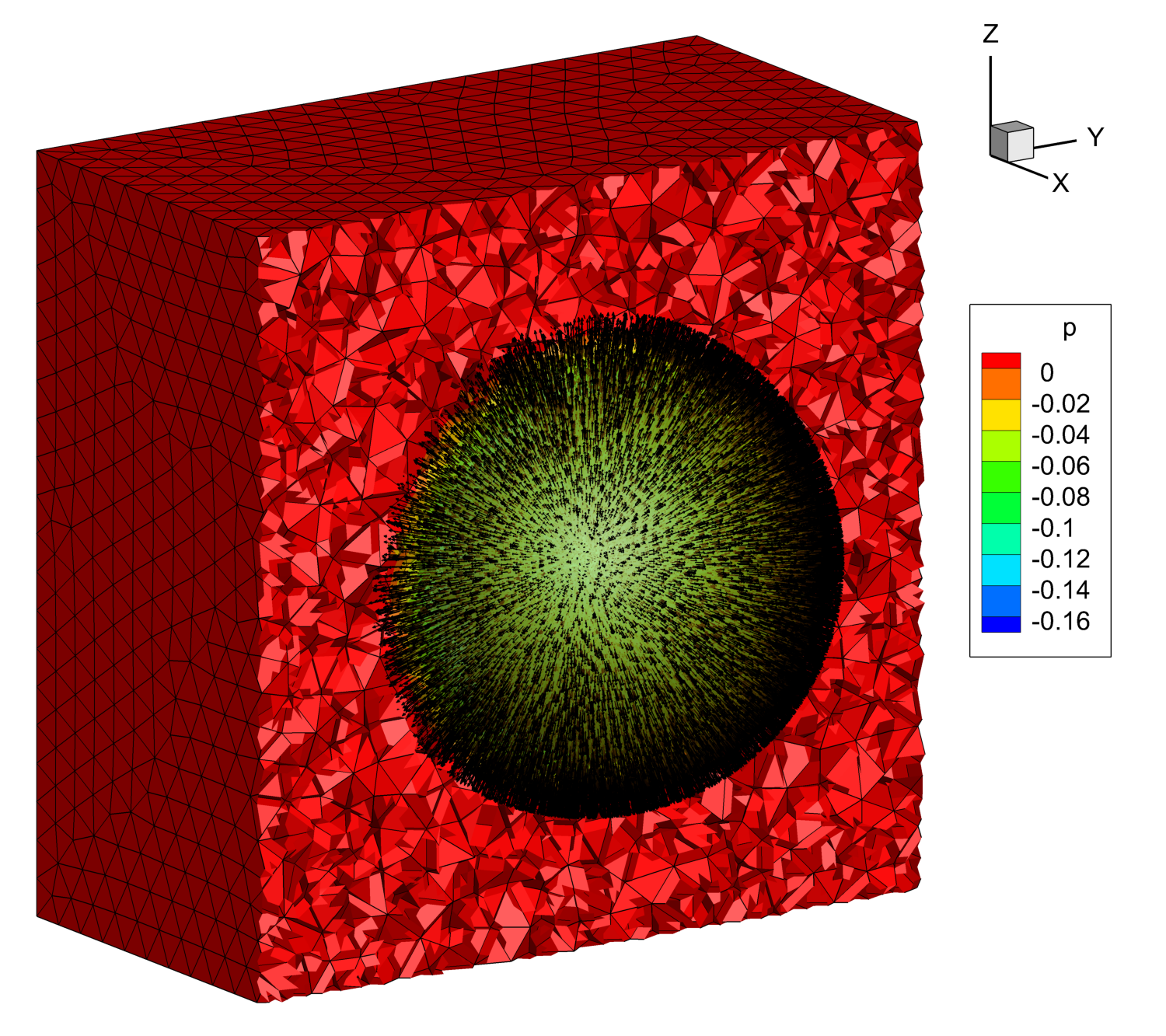}  
		\end{tabular} 
		\caption{{Numerical solution of the equations of linear acoustics \eqref{eqn.acoustics.v}-\eqref{eqn.acoustics.p} using the new structure-preserving semi-implicit DG scheme with $N=3$ in two and three space dimensions, respectively. Time evolution of the total energy in 2D (top left) and in 3D (top right), of the curl error of $\v$ in 2D (center left) and in 3D (center right). Mesh, contour colors and velocity vectors at time $t=3$ in the 2D case (bottom left).  Mesh, contour surfaces and velocity  vectors at time $t=3$ in the 3D case (bottom right).}} 
		\label{fig.acoustics}
	\end{center}
\end{figure}

\begin{figure}[!htbp]
	\begin{center}
		\begin{tabular}{cc} 
			\includegraphics[width=0.45\textwidth]{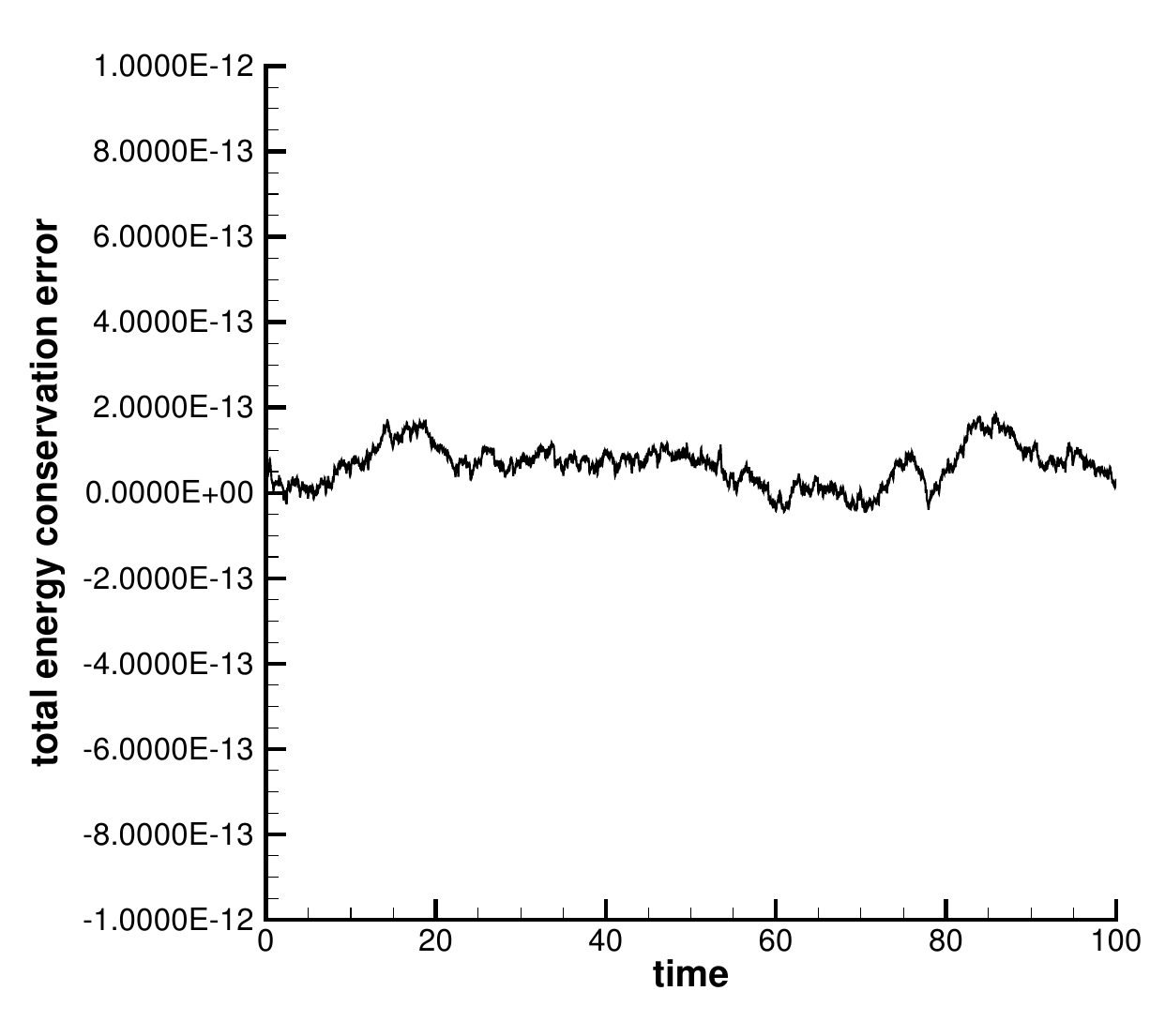}     & 
			\includegraphics[width=0.45\textwidth]{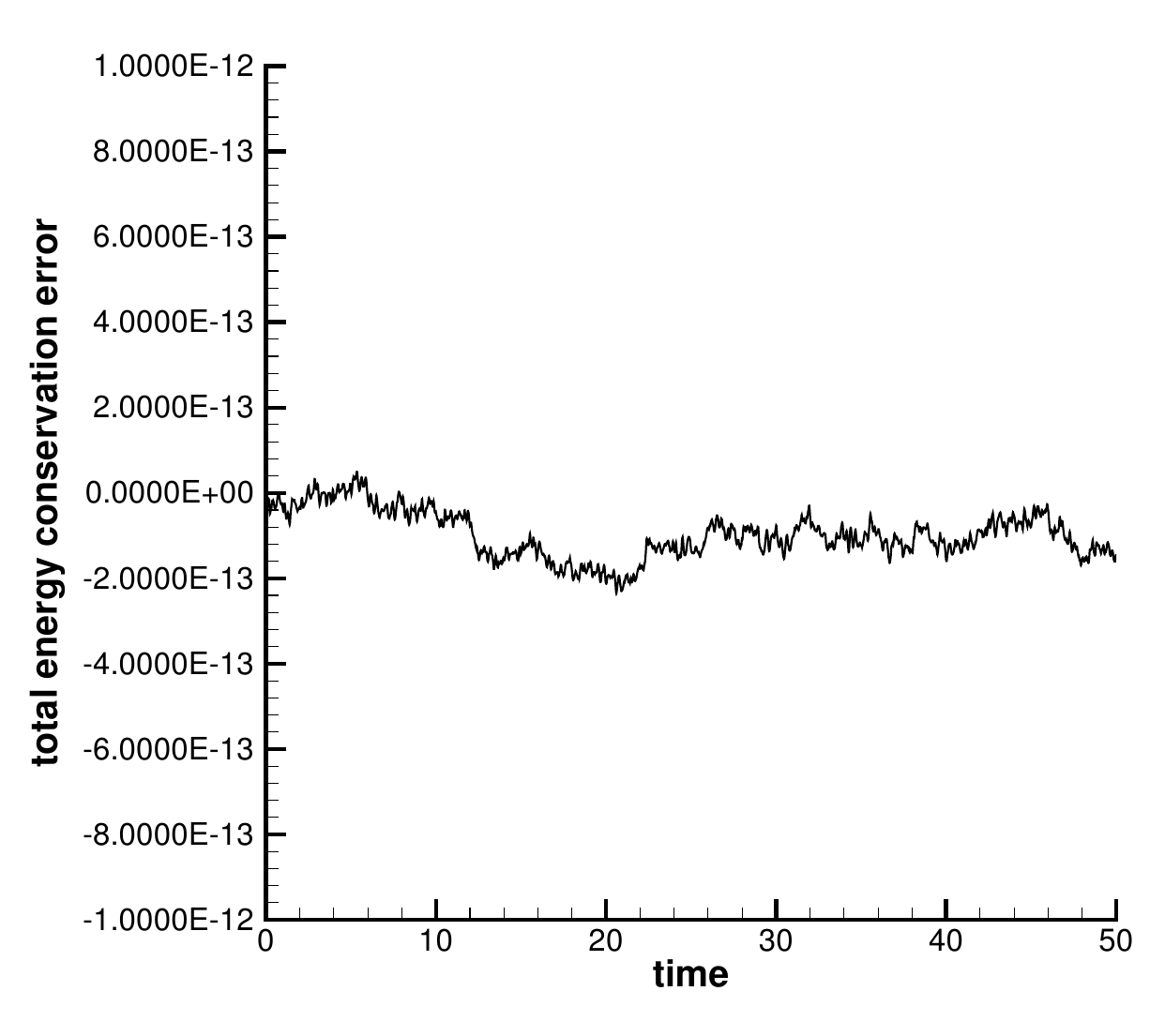}      
		\end{tabular} 
		\caption{{Time evolution of the total energy error in 2D (left) and in 3D (right).} } 
		\label{fig.acoustics.energyerror}
	\end{center}
\end{figure}

The first test concerns the equations of linear acoustics \eqref{eqn.acoustics.v}-\eqref{eqn.acoustics.p} in two and three space dimensions. We check numerically that the new schemes introduced in this paper conserve discrete total energy exactly and that they are able to preserve the curl-free property of the velocity field. 
{We recall that according to the ansatz \eqref{eqn.acoustics.ansatz} we assume $\vh \in \mathcal{U}^N_h$ and $\ph \in \mathcal{W}^{N+1}_h$. }
The initial condition of the test problem is given by 
\begin{equation}
  \v(\x, 0)  = 0, \qquad 
  p(\x, 0)  = p_0 \exp \left( -\halb \x^2/\sigma^2 \right).
	\label{eqn.ic.gauss.max.acoustic}
\end{equation}
The simulations are carried out in the domain $\Omega=[-\halb,+\halb]^d$ in $d$ space dimensions. In two space dimensions (2D, $d=2$) we run the simulation until a final time of $t=100$ and use an unstructured triangular mesh composed of a total number of 2028 triangles with $N_x=N_y=30$ elements along each edge of the computational domain. Instead, in three space dimensions (3D, $d=3$) the final time is set to $t=50$ and an unstructured mesh of $N_x=N_y=N_z=20$ elements along each edge of the domain is employed, leading to a total number of 59334 tetrahedra. In all cases we use a DG scheme of approximation degree $N=3$. For the initial data we set $p_0 = 1$. In the 2D test the half width of the Gaussian pressure pulse is chosen as $\sigma=0.05$, while it is set to $\sigma=0.1$ in 3D. 
The computational results obtained for both tests are depicted in Figure \ref{fig.acoustics}. 
{In addition, in Figure \ref{fig.acoustics.energyerror} we show the absolute total energy conservation error for the two and three-dimensional case, respectively. 
As expected, the method conserves total energy up to machine precision even over long times and the discrete velocity field remains curl-free up to machine precision, both in two and three space dimensions. } 

{
As a second test in this section we consider long-time wave propagation on coarse meshes, in order to check the dispersion properties of our schemes numerically. The initial data are given by
\begin{equation}
	p(\x,0) = v_1(\x,0) = \sin\left( 2 \pi \frac{x}{\lambda} \right),  \qquad v_2=v_3=0,
\end{equation}
with the wavelength $\lambda=0.25$. The computational domain is the rectangle $\Omega = [-0.5, +0.5] \times [-0.1, +0.1]$ and periodic boundary conditions are set everywhere. We test our new compatible DG scheme with polynomial degree $N=3$ using a very coarse mesh of characteristic mesh size $h=0.05$, corresponding to only 186 triangular elements. The discrete initial velocity field is computed as the compatible gradient of the scalar potential $Z= -\frac{\lambda}{2 \pi} \cos\left( 2 \pi \frac{x}{\lambda} \right)$ and the time step is fixed at $\Delta t = 10^{-3}$.
Simulations are run until a final time of $t=100$, which corresponds to 100 traveling periods of the signal in the given domain. The employed mesh and the pressure contours of the numerical solution at the final time, as well as a comparison between exact solution and numerical solution at time $t=100$ on a 1D cut through the solution at $y=0$ are depicted in Figure \ref{fig.acwave}. As our method is exactly energy conservative, one can observe a small phase shift in the signal, which corresponds to the numerical dispersion error of the method and which is actually expected. To reduce the dispersion error, higher order symplectic time integrators should be used, but this is clearly out of the scope of the present work and is left to future research. The curl errors remain of the order of machine precision for all times, as expected. 
    \begin{figure}[!htbp]
    	\begin{center}
    		\begin{tabular}{c} 
    			\includegraphics[trim=10 10 10 10,clip,width=0.55\textwidth]{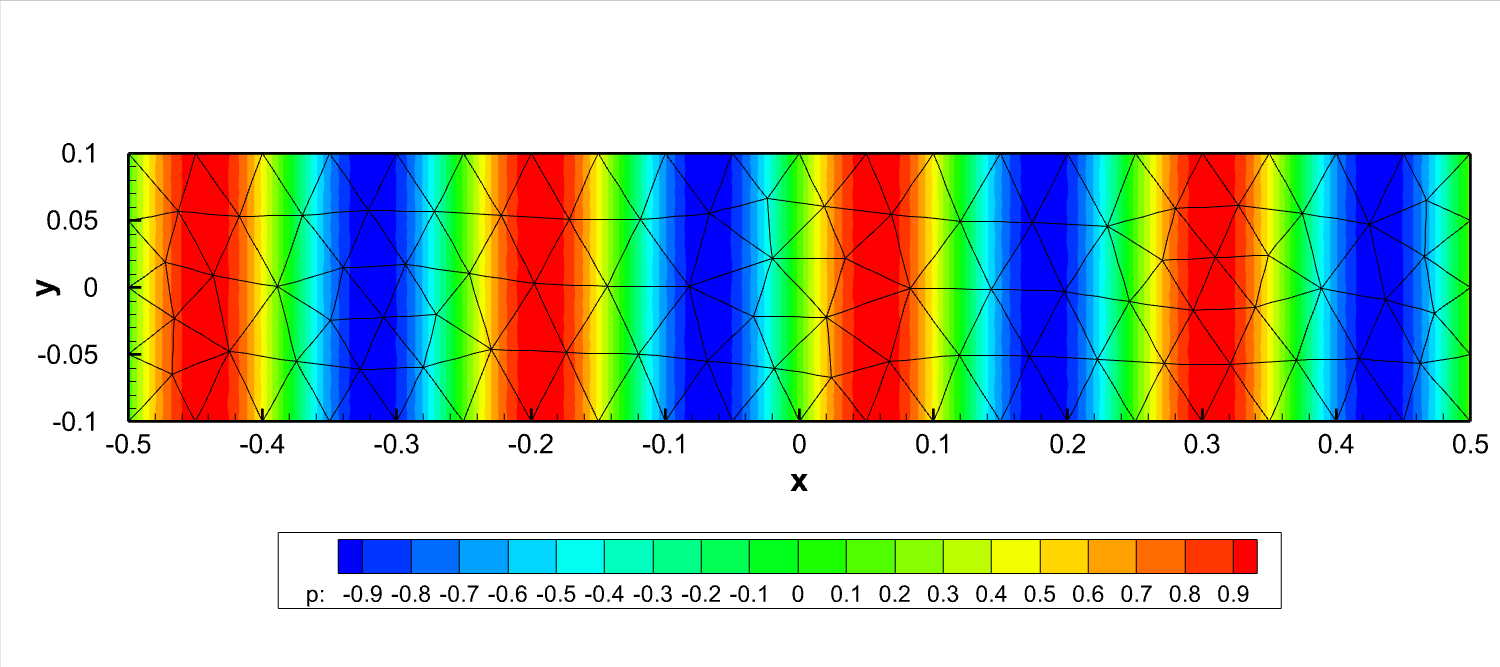}     \\ 
    			\includegraphics[width=0.55\textwidth]{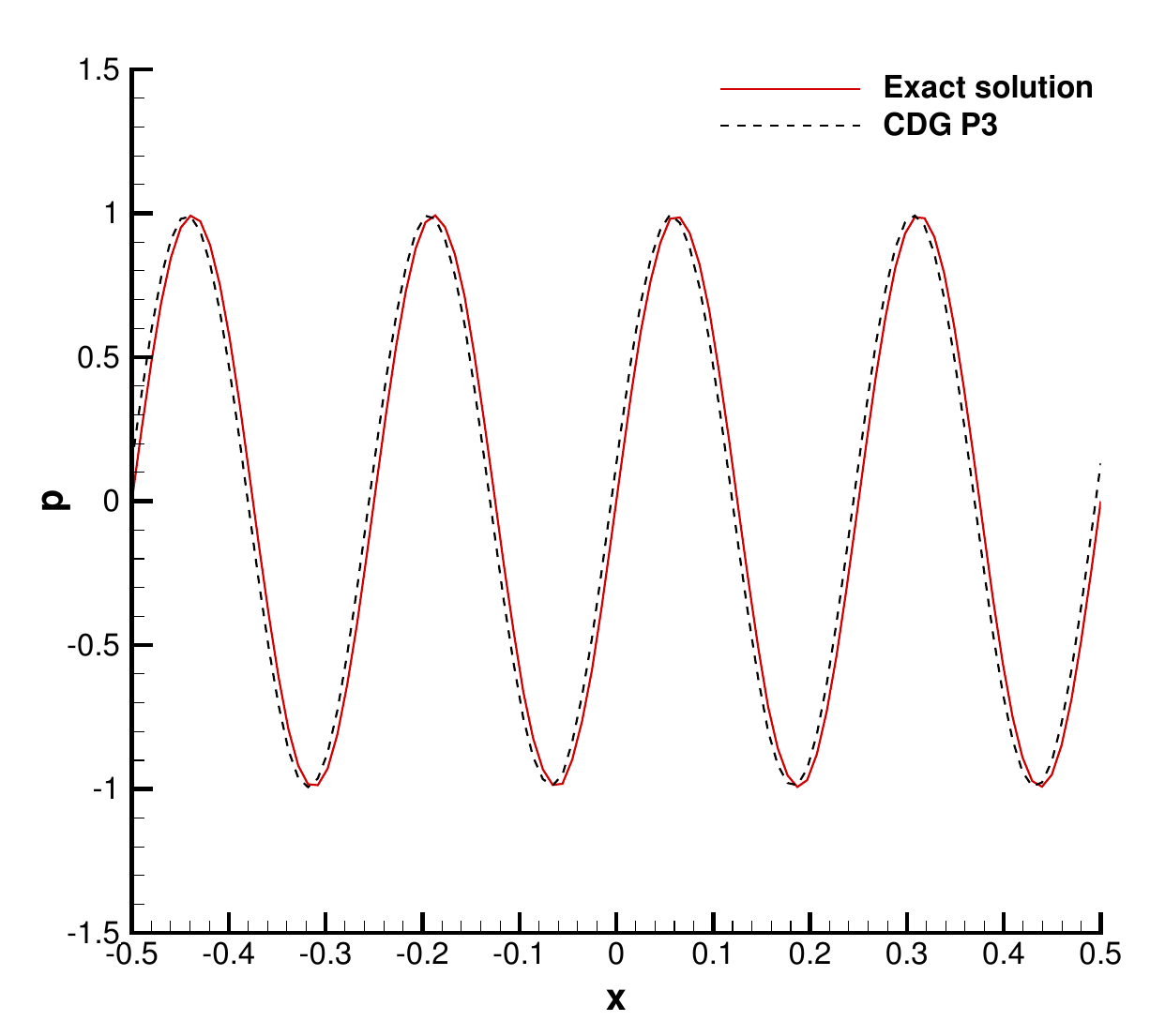}   \\    
    			\includegraphics[width=0.55\textwidth]{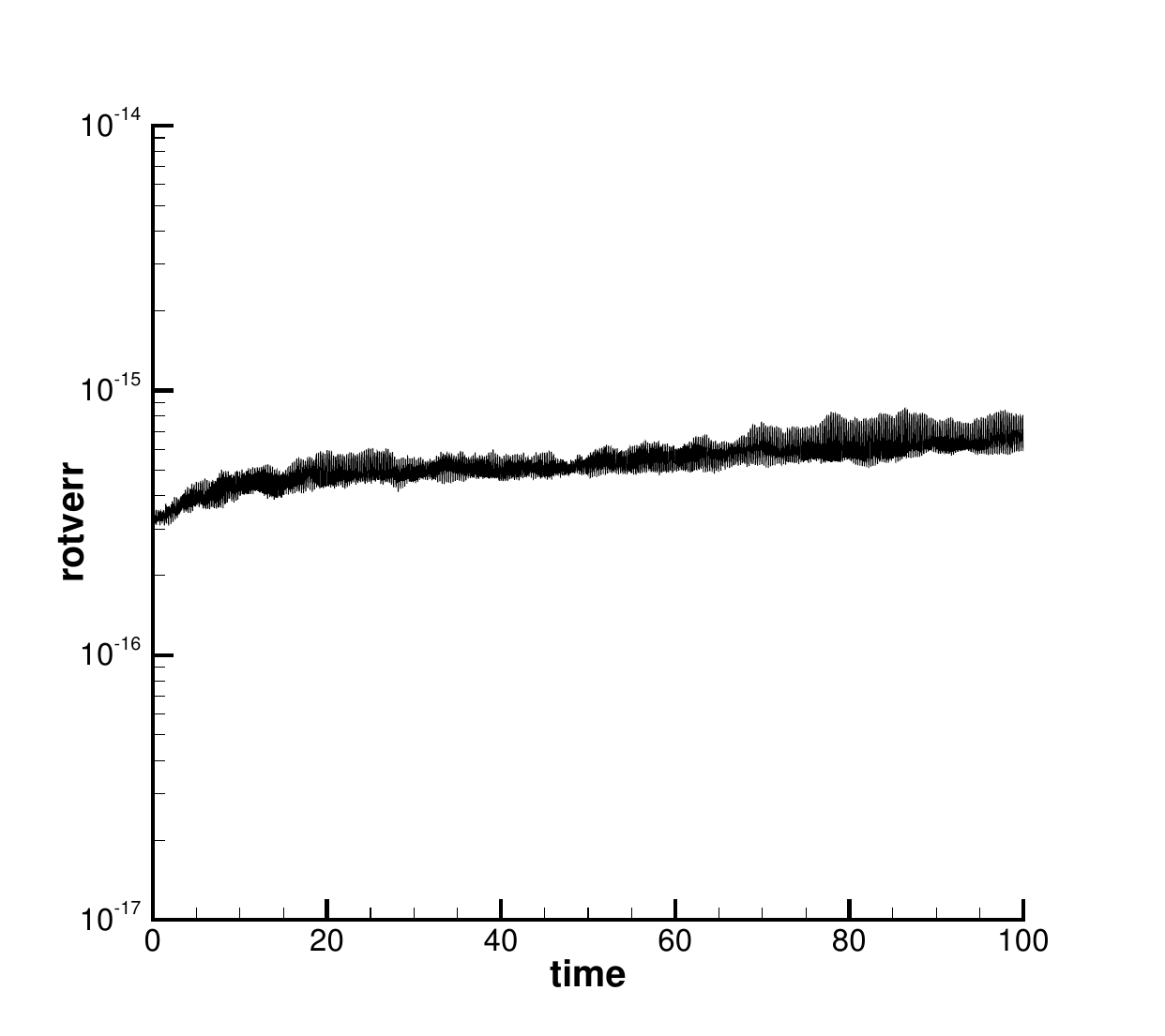}      
    		\end{tabular} 
    		\caption{{Computational mesh and pressure contour colors of the numerical solution at $t=100$ (top). Comparison between exact and numerical solution obtained with the new compatible DG scheme (CDG) at time $t=100$ (center). A slight phase shift in the signal is visible, related to the dispersion error of the numerical scheme. Time series of the curl error of the velocity field (bottom). } } 
    		\label{fig.acwave}
    	\end{center}
    \end{figure}
}

\subsection{Vacuum Maxwell equations} 

In this section we apply the new structure-preserving DG schemes proposed in this paper to the vacuum Maxwell equations \eqref{eqn.maxwell.B}-\eqref{eqn.maxwell.E} in two and three space dimensions. 
{According to the ansatz \eqref{eqn.maxwell.ansatz} we assume $\Bh \in \mathcal{U}^N_h$ and $\Eh \in \mathcal{W}^{N+1}_h$. }
We check numerically that they conserve total energy exactly and that they are able to preserve an exactly divergence-free magnetic field for all times. 

\begin{figure}[!htbp]
	\begin{center}
		\begin{tabular}{cc} 
			\includegraphics[width=0.45\textwidth]{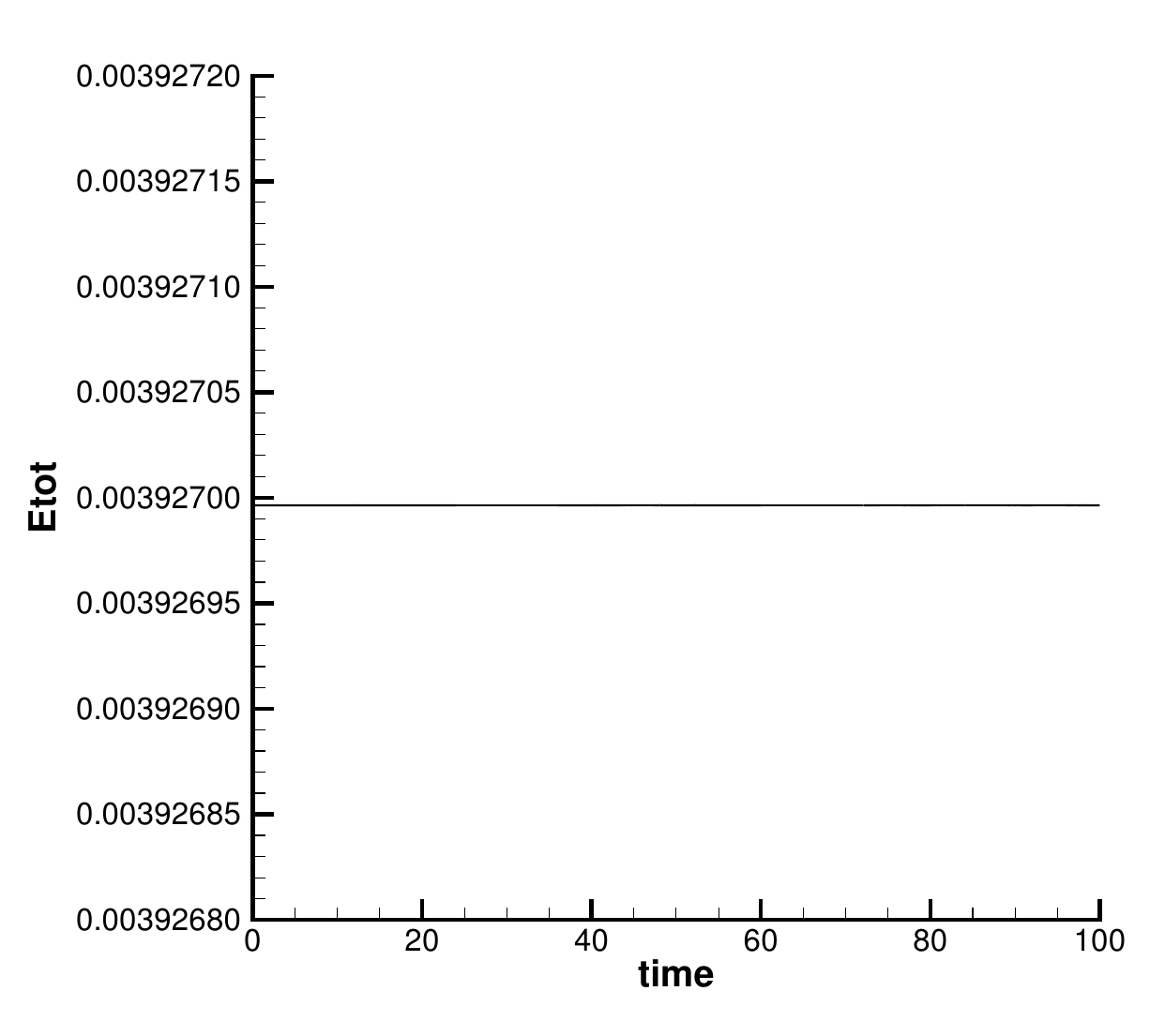}     & 
			\includegraphics[width=0.45\textwidth]{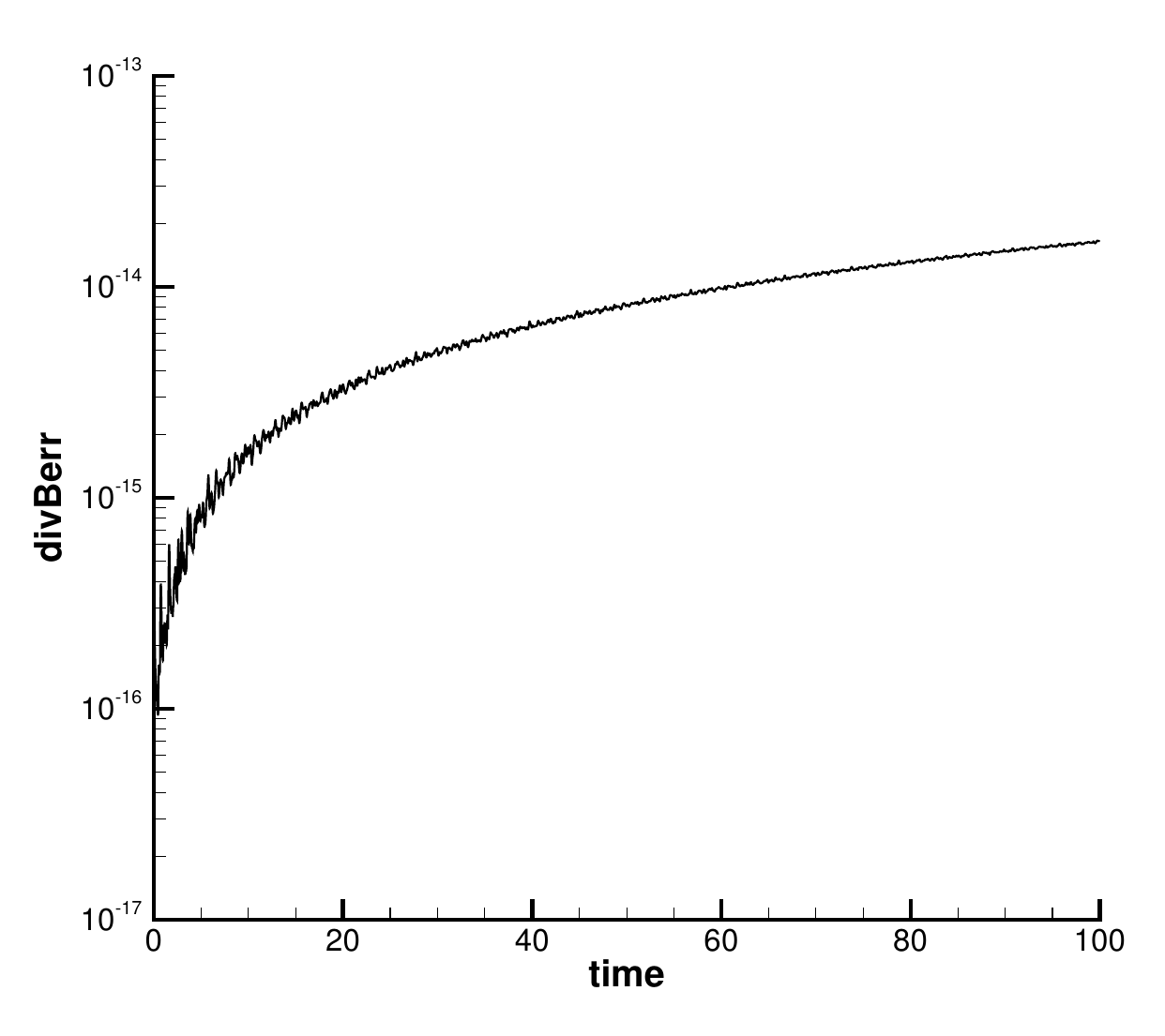}  \\  
			\includegraphics[width=0.45\textwidth]{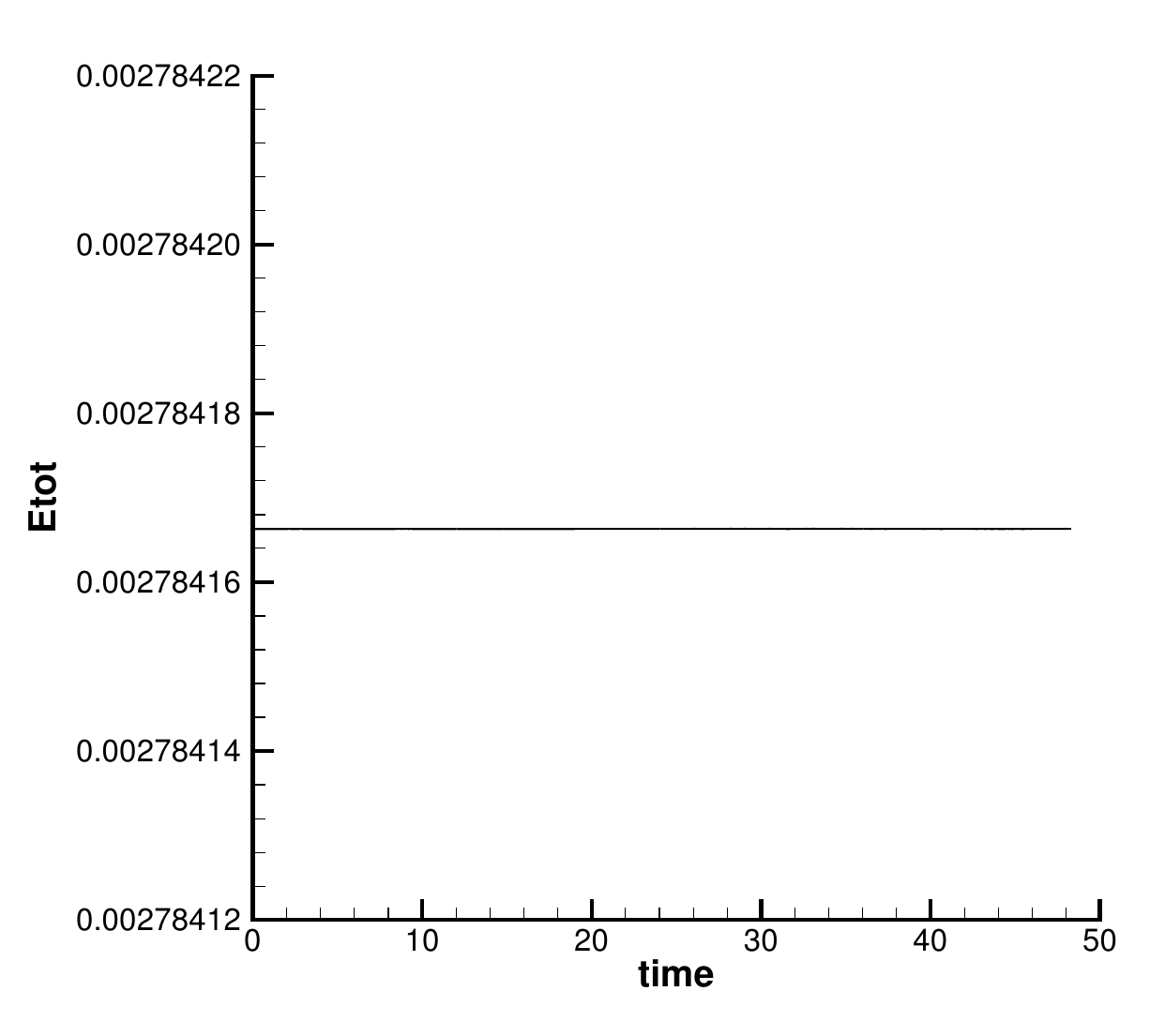}     & 
			\includegraphics[width=0.45\textwidth]{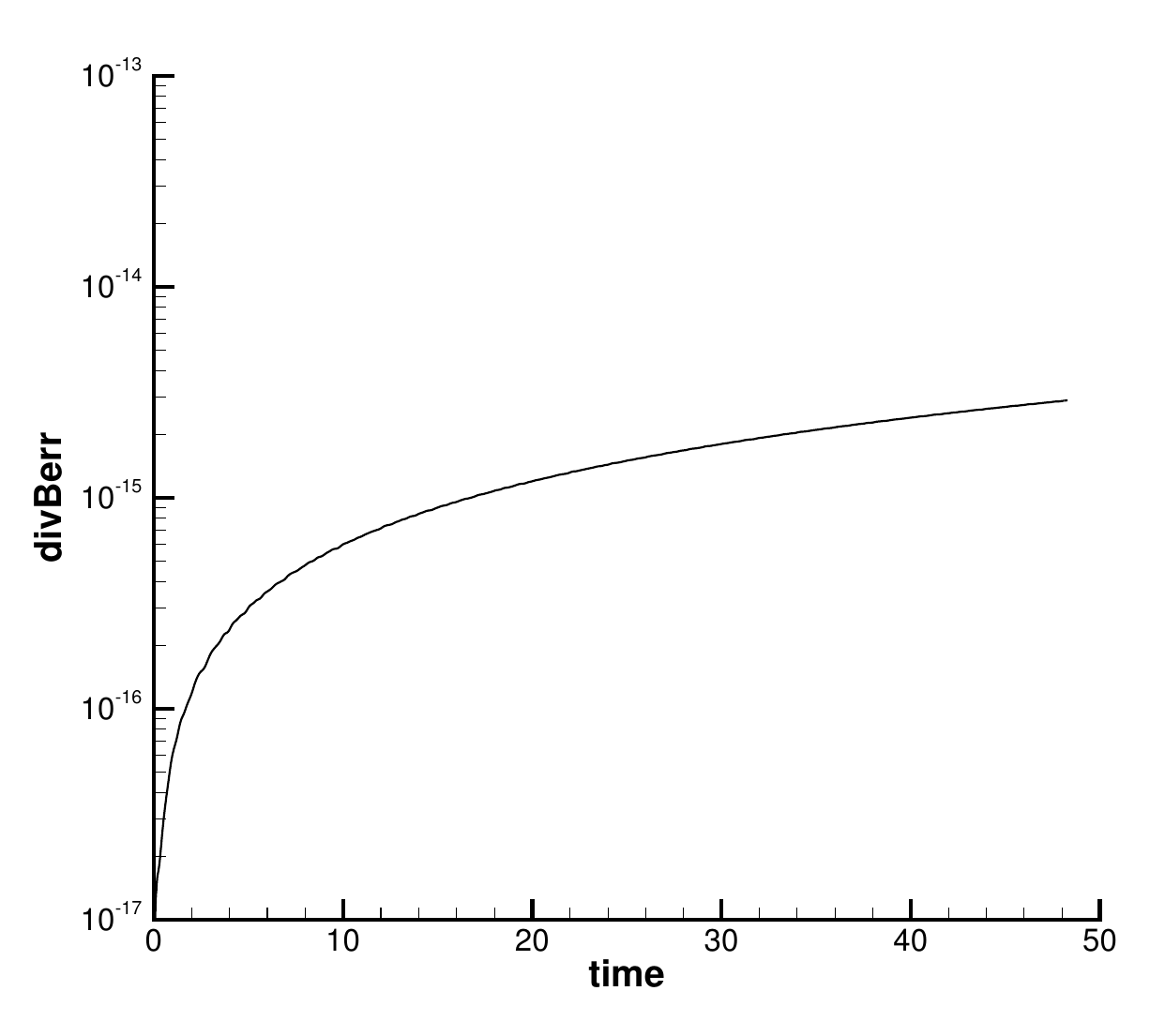}  \\  
			\includegraphics[width=0.45\textwidth]{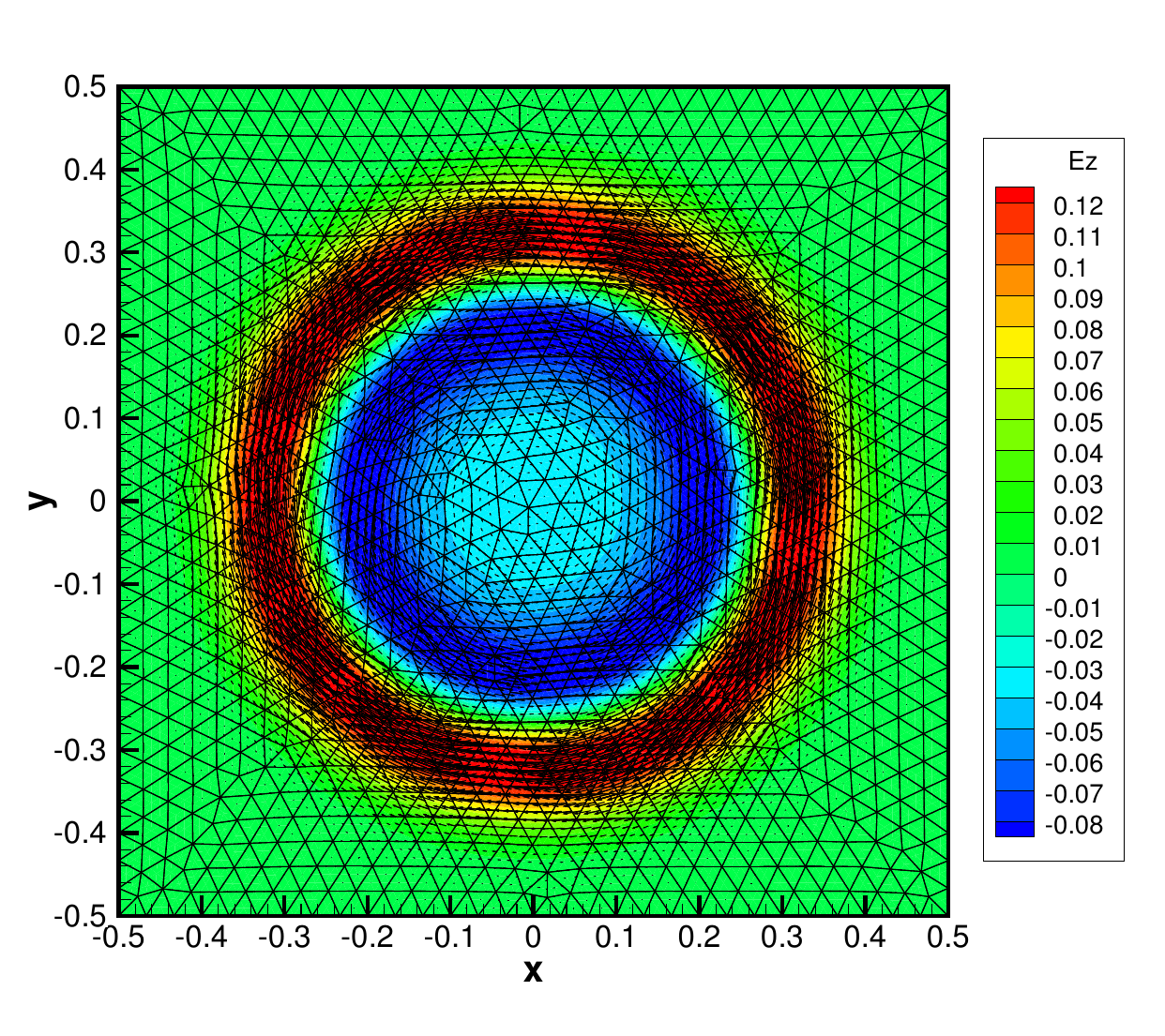}     & 
			\includegraphics[width=0.45\textwidth]{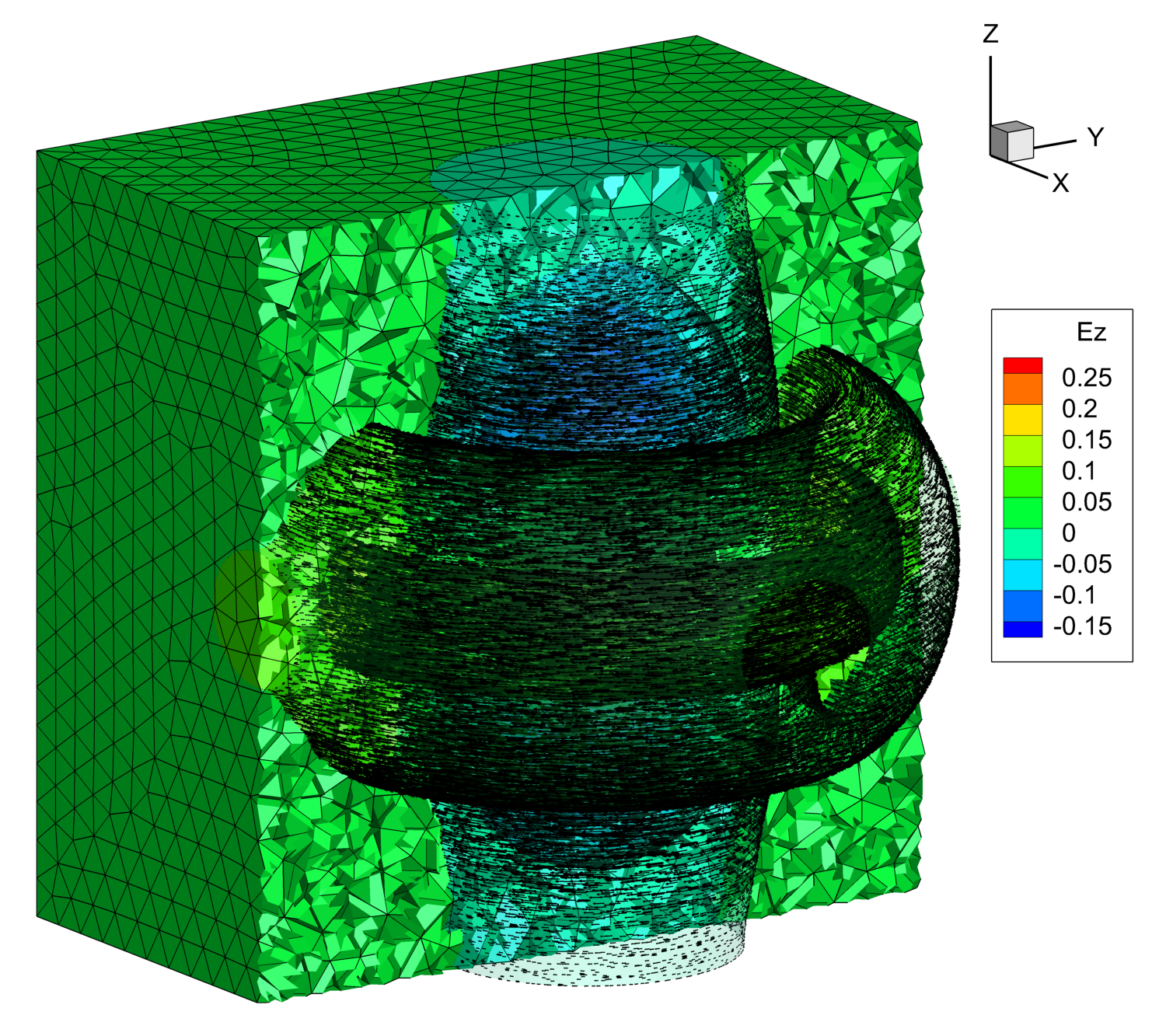}  
		\end{tabular} 
		\caption{Numerical solution of the vacuum Maxwell equations \eqref{eqn.maxwell.B}-\eqref{eqn.maxwell.E} until a final time of $t=100$ using the new structure-preserving semi-implicit DG scheme with $N=3$ in two and three space dimensions. Time evolution of the total energy in 2D (top left) and of the divergence error of $\B$ in 2D (top right). Time evolution of the total energy in 3D (center left) and of the divergence error of $\B$ in 3D (center right).  Mesh, contour colors and magnetic field vectors at time $t=3$ in the 2D case (bottom left).  Mesh, contour surfaces and magnetic field vectors at time $t=3$ in the 3D case (bottom right).} 
		\label{fig.maxwell}
	\end{center}
\end{figure}

To test the properties of the schemes, we solve the following test problem whose initial condition reads  
\begin{equation}
  \B(\x, 0)  = \mathbf{B}_0 \exp \left( -\halb \x^2/\sigma^2 \right), \qquad 
  \E(\x, 0)  = \mathbf{E}_0 \exp \left( -\halb \x^2/\sigma^2 \right).
	\label{eqn.ic.gauss.max.maxwell}
\end{equation}
As the the previous test, simulations are run in the domain $\Omega=[-\halb,+\halb]^d$ in $d$ space dimensions until a final time of $t=100$ in 2D and until $t=50$ in 3D, using a DG scheme of approximation degree $N=3$ in all cases. The employed computational meshes are as in the previous testcase. For the initial data we set $\B_0 = 0$ and $\E_0 = (0,0,1)$. In the 2D test the half width is chosen as $\sigma=0.05$, while it is set to $\sigma=0.1$ in 3D. 
The computational results obtained for both tests are depicted in Figure \ref{fig.maxwell}. As expected, the method conserves total energy perfectly even over long times and the involutions are satisfied up to machine precision, both in two and three space dimensions. 

\subsection{The Maxwell-GLM system} 

We first carry out a systematic numerical convergence study of the new structure-preserving scheme presented in this paper and applied to the Maxwell-GLM system  \eqref{eqn.maxmunz.B}-\eqref{eqn.maxmunz.psi}. 
{We recall that according to the ansatz \eqref{eqn.maxmunz.ansatz} we assume $\Bh \in \mathcal{U}^N_h$, $\qh \in \mathcal{U}^N_h$, $\Eh \in \mathcal{W}^{N+1}_h$ and  $\ph \in \mathcal{W}^{N+1}_h$. }
For this purpose we consider the following initial data of a planar wave travelling in the direction $\n = (1,1,0)$, see \cite{MaxwellGLM} for details. The initial condition 
is given by 
\begin{eqnarray}
	& \mathbf{B}(\x,0)  = \mathbf{B}_0 \sin\left( \pi (x_1-x_2) \right), \qquad & \varphi(\x,0) = \varphi_0 \sin \left(\pi (x_1-x_2)\right),	\\
	& \mathbf{E}(\x,0)  = \mathbf{E}_0 \sin\left( \pi (x_1-x_2) \right), \qquad & \psi(\x,0) = \psi_0  \sin\left(\pi (x_1-x_2)\right), 
\end{eqnarray}
with $\mathbf{B}_0=(0.25b, -0.25b, 1)^T $, $ \mathbf{E}_0=(1.5b, 0.5b, 0)^T $, $\varphi_0 = 0.25$, $\psi_0 = 0.5$ and $b=\sqrt{2}/2$.
The computational domain $\Omega=[-1,1]^2$ is discretized with a sequence of successively refined unstructured meshes with $N_x=N_y$ elements along each edge of $\Omega$. Simulations are run for one period, i.e. until a final time of $t=\sqrt{2}$, when the exact solution of the problem coincides again with the initial condition. Since this is a time-dependent test problem, to reduce the time discretization errors the CFL number has substantially been decreased for schemes of better than second order in space. Table \ref{tab.mm.conv} shows the $L^2$ errors of a large set of different field variables at the final time, depending on the chosen grid spacing and the polynomial approximation degree $N$. As one can observe, the method converges with the expected order of accuracy $N+1$ in space for all variables and at all polynomial approximation degrees $N$ under consideration. 

\begin{sidewaystable}  
		\caption{Numerical convergence results of the new structure-preserving DG scheme for the Maxwell-GLM system for polynomial approximation degrees $N \in [0,3]$. The $L^2$ error norms refer to the components $B_1$, $B_2$, $p$, $E_1$, $E_2$ and $q$ at the final time of $t=\sqrt{2}$. } 	
		\begin{center} 
				\renewcommand{\arraystretch}{1.1}
				\begin{tabular}{ccccccccccccc} 
                    \hline 
                    \multicolumn{13}{c}{$N=0$, CFL$=1$} \\ 
					\hline
					$N_x$ & ${L^2}(B_1)$ & $\mathcal{O}(B_1)$  & ${L^2}(B_2)$ & $\mathcal{O}(B_2)$  
					& ${L^2}(p)$ & $\mathcal{O}(p)$ & 
                    ${L^2}(E_1)$ & $\mathcal{O}(E_1)$  & ${L^2}(E_2)$ & $\mathcal{O}(E_2)$  
					& ${L^2}(q)$ & $\mathcal{O}(q)$ 
					\\ 
					\hline
					20	& 3.7615E-02 &      & 3.7620E-02 &       & 5.4471E-02 &       & 2.3110E-01 &      & 7.7033E-02 &       & 1.0644E-01 &       \\   
					40	& 1.9009E-02 & 1.0  & 1.9009E-02 &  1.0  & 2.7080E-02 & 1.0   & 1.1489E-01 & 1.0  & 3.8296E-02 &  1.0  & 5.3772E-02 & 1.0   \\  
					60	& 1.2773E-02 & 1.0  & 1.2774E-02 &  1.0  & 1.8118E-02 & 1.0   & 7.6869E-02 & 1.0  & 2.5623E-02 &  1.0  & 3.6127E-02 & 1.0   \\  
					80	& 9.2003E-03 & 1.1  & 9.1998E-03 &  1.1  & 1.3037E-02 & 1.1   & 5.5312E-02 & 1.1  & 1.8437E-02 &  1.1  & 2.6023E-02 & 1.1   \\ 
					120 & 6.1640E-03 & 1.0  & 6.1638E-03 &  1.0  & 8.7260E-03 & 1.0   & 3.7021E-02 & 1.0  & 1.2340E-02 &  1.0  & 1.7435E-02 & 1.0   \\  
                    \hline 
                    \multicolumn{13}{c}{$N=1$, CFL$=1$} \\ 
					\hline
					$N_x$ & ${L^2}(B_1)$ & $\mathcal{O}(B_1)$  & ${L^2}(B_2)$ & $\mathcal{O}(B_2)$  
					& ${L^2}(p)$ & $\mathcal{O}(p)$ & 
                    ${L^2}(E_1)$ & $\mathcal{O}(E_1)$  & ${L^2}(E_2)$ & $\mathcal{O}(E_2)$  
					& ${L^2}(q)$ & $\mathcal{O}(q)$ 
					\\ 
					\hline
					20	& 2.0402E-03 &      & 2.0399E-03 &       & 2.8843E-03 &       & 1.2237E-02 &      & 4.0790E-03 &       & 5.7708E-03 &       \\   
					40	& 4.6894E-04 & 2.1  & 4.6894E-04 &  2.1  & 6.6343E-04 & 2.1   & 2.8147E-03 & 2.1  & 9.3823E-04 &  2.1  & 1.3264E-03 & 2.1   \\  
					60	& 2.2859E-04 & 1.8  & 2.2859E-04 &  1.8  & 3.2328E-04 & 1.8   & 1.3716E-03 & 1.8  & 4.5719E-04 &  1.8  & 6.4654E-04 & 1.8   \\  
					80	& 1.1477E-04 & 2.4  & 1.1477E-04 &  2.4  & 1.6232E-04 & 2.4   & 6.8867E-04 & 2.4  & 2.2956E-04 &  2.4  & 3.2463E-04 & 2.4   \\ 
					120 & 5.0808E-05 & 2.0  & 5.0808E-05 &  2.0  & 7.1856E-05 & 2.0   & 3.0486E-04 & 2.0  & 1.0162E-04 &  2.0  & 1.4371E-04 & 2.0   \\  
                    \hline 
                    \multicolumn{13}{c}{$N=2$, CFL$=0.1$} \\ 
					\hline
					$N_x$ & ${L^2}(B_1)$ & $\mathcal{O}(B_1)$  & ${L^2}(B_2)$ & $\mathcal{O}(B_2)$  
					& ${L^2}(p)$ & $\mathcal{O}(p)$ & 
                    ${L^2}(E_1)$ & $\mathcal{O}(E_1)$  & ${L^2}(E_2)$ & $\mathcal{O}(E_2)$  
					& ${L^2}(q)$ & $\mathcal{O}(q)$ 
					\\ 
					\hline
					10	& 4.5152E-04 &      & 4.5152E-04 &       & 6.4454E-04 &       & 2.7346E-03 &      & 9.1152E-04 &       & 1.2789E-03 &       \\   
					20	& 5.4695E-05 & 3.0  & 5.4686E-05 &  3.0  & 7.7540E-05 & 3.1   & 3.2897E-04 & 3.1  & 1.0966E-04 &  3.1  & 1.5475E-04 & 3.0   \\  
					30	& 1.6362E-05 & 3.0  & 1.6362E-05 &  3.0  & 2.3171E-05 & 3.0   & 9.8306E-05 & 3.0  & 3.2769E-05 &  3.0  & 4.6285E-05 & 3.0   \\  
					40	& 7.0805E-06 & 2.9  & 7.0805E-06 &  2.9  & 1.0020E-05 & 2.9   & 4.2512E-05 & 2.9  & 1.4171E-05 &  2.9  & 2.0028E-05 & 2.9   \\ 
					50  & 3.7153E-06 & 2.9  & 3.7153E-06 &  2.9  & 5.2575E-06 & 2.9   & 2.2305E-05 & 2.9  & 7.4351E-06 &  2.9  & 1.0510E-05 & 2.9   \\  
                    \hline 
                    \multicolumn{13}{c}{$N=3$, CFL$=h$} \\ 
					\hline
					$N_x$ & ${L^2}(B_1)$ & $\mathcal{O}(B_1)$  & ${L^2}(B_2)$ & $\mathcal{O}(B_2)$  
					& ${L^2}(p)$ & $\mathcal{O}(p)$ & 
                    ${L^2}(E_1)$ & $\mathcal{O}(E_1)$  & ${L^2}(E_2)$ & $\mathcal{O}(E_2)$  
					& ${L^2}(q)$ & $\mathcal{O}(q)$ 
					\\ 
					\hline
					 5	& 3.7012E-04 &      & 3.7065E-04 &       & 5.2579E-04 &       & 2.2307E-03 &      & 7.4358E-04 &       & 1.0512E-03 &       \\   
					10	& 2.3266E-05 & 4.0  & 2.3263E-05 &  4.0  & 3.2970E-05 & 4.0   & 1.3988E-04 & 4.0  & 4.6627E-05 &  4.0  & 6.5916E-05 & 4.0   \\  
					15	& 4.5898E-06 & 4.0  & 4.5901E-06 &  4.0  & 6.4946E-06 & 4.0   & 2.7553E-05 & 4.0  & 9.1843E-06 &  4.0  & 1.2986E-05 & 4.0   \\  
					20	& 1.4010E-06 & 4.1  & 1.4010E-06 &  4.1  & 1.9853E-06 & 4.1   & 8.4201E-06 & 4.1  & 2.8067E-06 &  4.1  & 3.9681E-06 & 4.1   \\ 
					25  & 5.6858E-07 & 4.0  & 5.6772E-07 &  4.0  & 8.0137E-07 & 4.1   & 3.3969E-06 & 4.1  & 1.1323E-06 &  4.1  & 1.6001E-06 & 4.1   \\  
					\hline 					
				\end{tabular}
		\end{center}
		\label{tab.mm.conv}
\end{sidewaystable}

\begin{figure}[!htbp]
	\begin{center}
		\begin{tabular}{cc} 
			\includegraphics[width=0.45\textwidth]{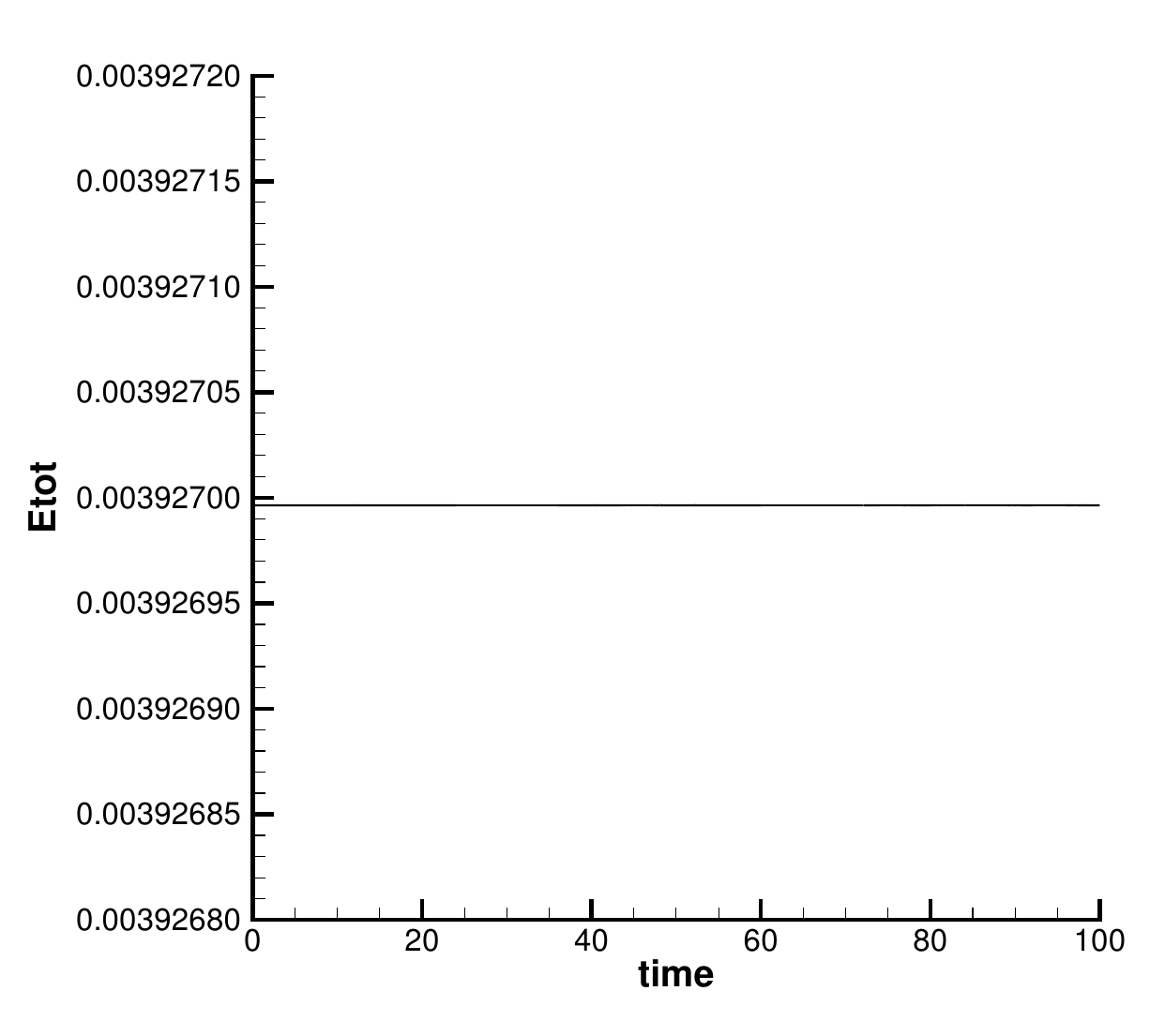}     & 
			\includegraphics[width=0.45\textwidth]{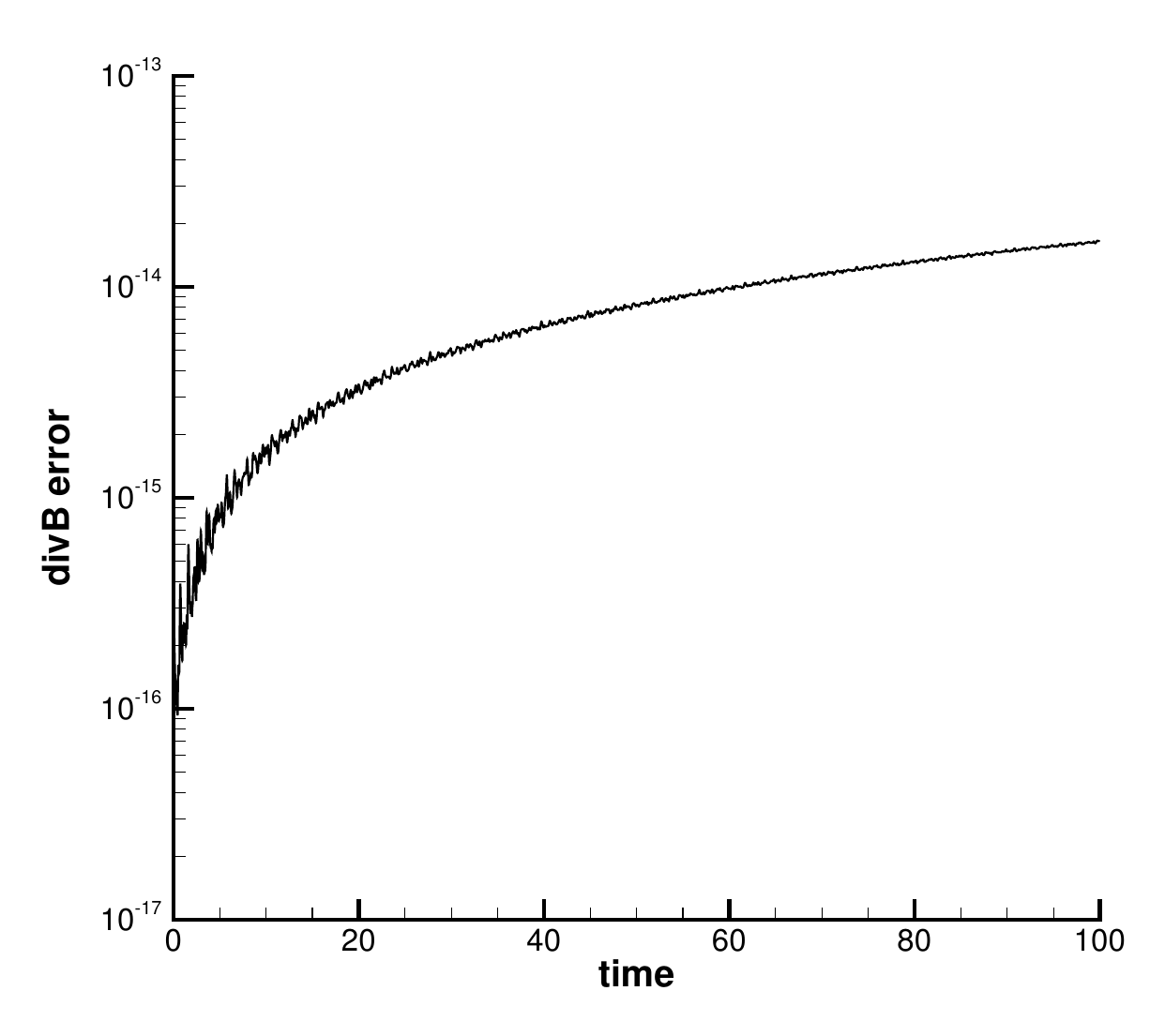}  \\  
			\includegraphics[width=0.45\textwidth]{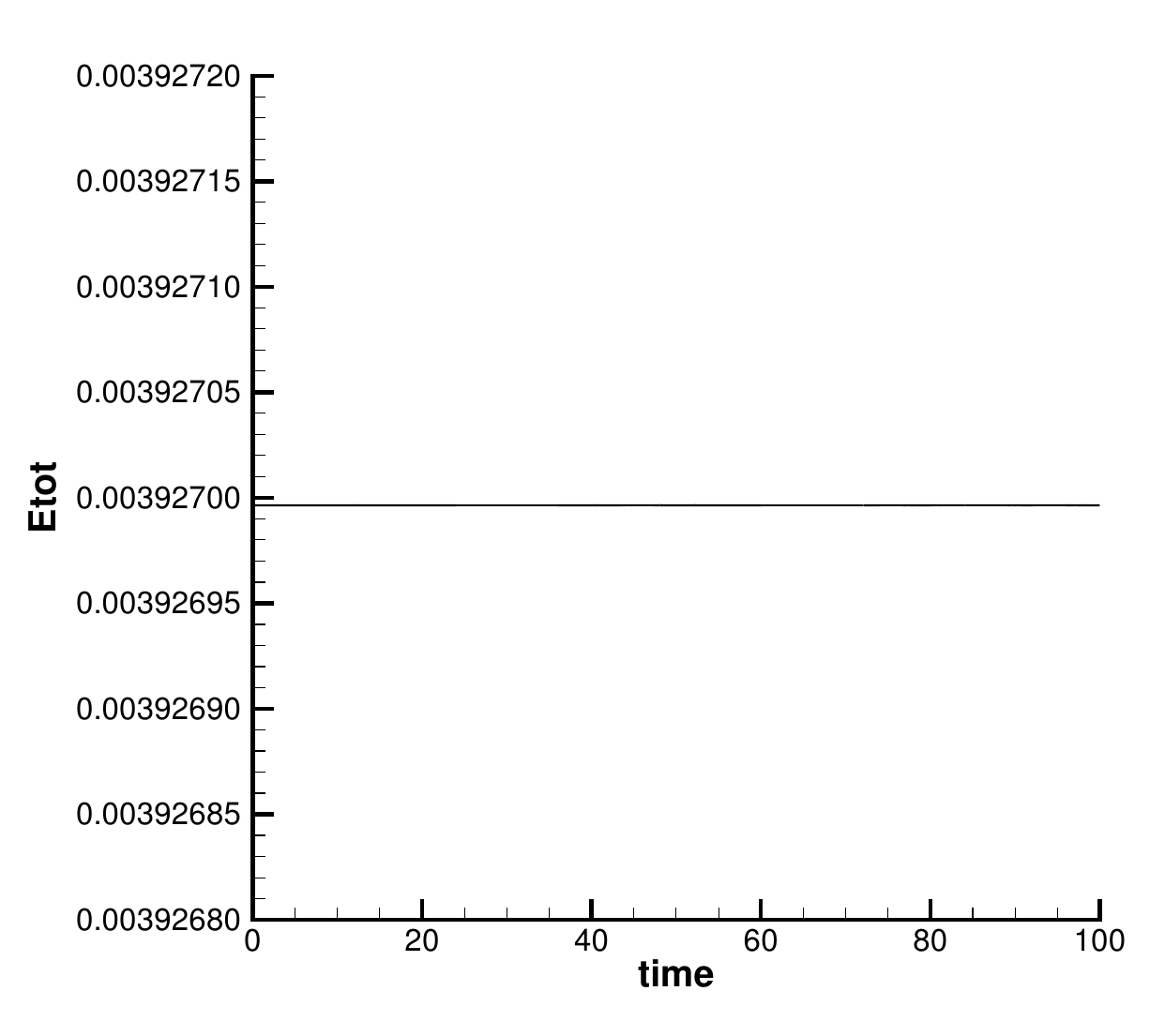}     & 
			\includegraphics[width=0.45\textwidth]{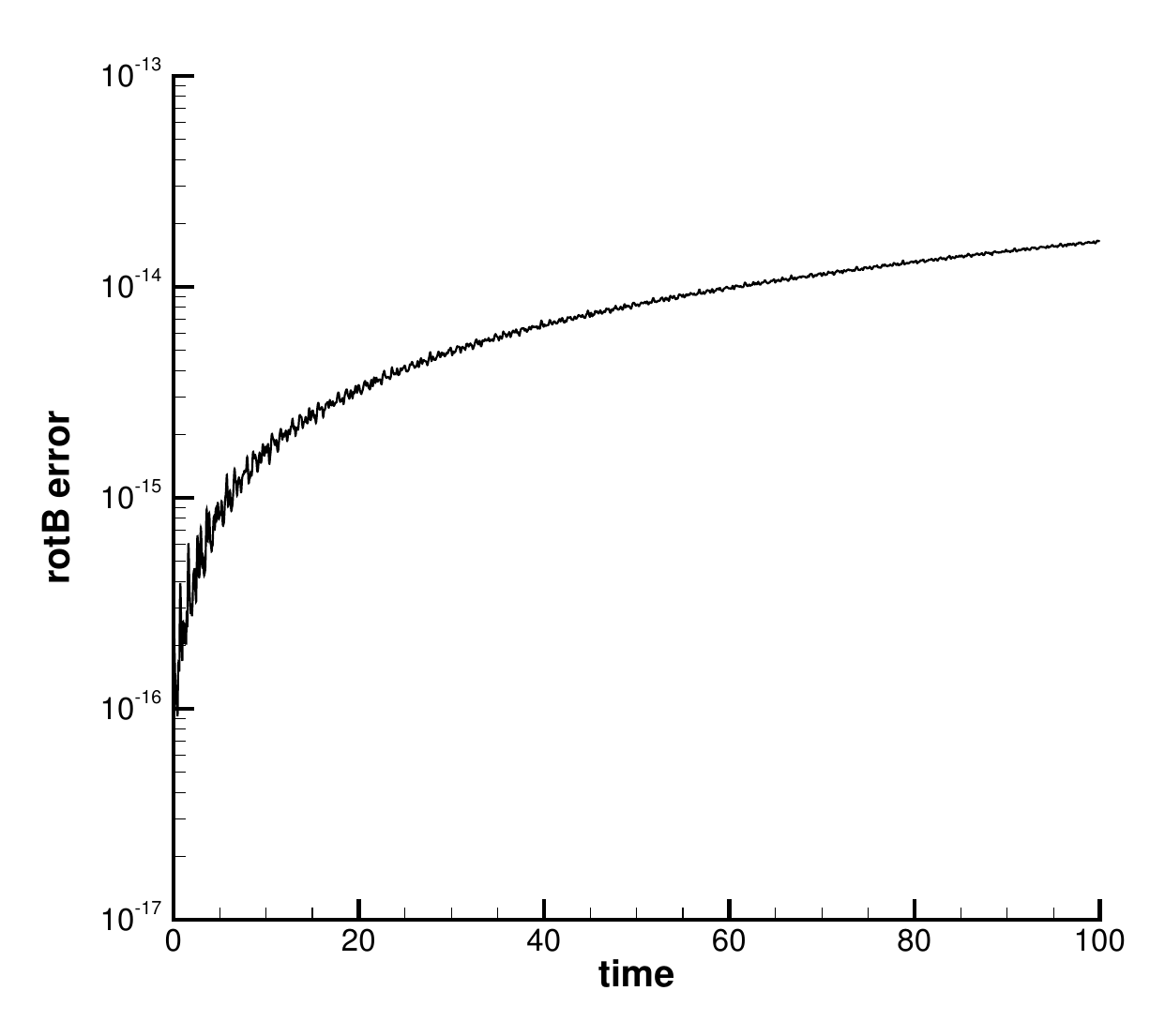}  
		\end{tabular} 
		\caption{Numerical solution of the Maxwell-GLM system \eqref{eqn.maxmunz.B}-\eqref{eqn.maxmunz.psi} until a final time of $t=100$ using the new structure-preserving semi-implicit DG scheme with $N=3$. Time evolution of the total energy (top left) and of the divergence error of $\B$ (top right) for test T1 (Maxwell-type initial data).
        Time evolution of the total energy (bottom left) and of the curl error of $\B$ (bottom right) for test T2 (acoustic-type initial data).} 
		\label{fig.mm}
	\end{center}
\end{figure}

We now verify via two numerical experiments that the proposed structure-preserving DG schemes conserve total energy and that they are able to satisfy the involutions for the vector field $\B$ exactly at the discrete level, depending on the chosen initial data. In the case 
$\phi(\x,0)=\psi(\x,0)=0$ one retrieves the vacuum Maxwell equations, hence $\B$  satisfies the involution $\nabla \cdot \B = 0$ if it was satisfied at the initial time, while for $\E=\psi=0$ the 
$\B$ field satisfies $\nabla \times \B = 0$ if it was initially curl-free. We emphasize that in the new numerical framework presented in this paper it is not necessary to adjust the approximation space of $\Bh$ to the particular involution that one wants to preserve, but the correct involution is automatically satisfied depending on the chosen initial data. This is a direct consequence of the compatibility of the chosen discrete nabla operators, primary and dual, and the two different piecewise polynomial approximation spaces, DG and FEM, as explained before. 
To test these properties of the schemes, we therefore solve the following family of test problems, given by the initial data 
\begin{eqnarray}
	& \B(\x, 0)  = \mathbf{B}_0 \exp \left( -\halb \x^2/\sigma^2 \right), \qquad & \varphi(\x, 0) = \varphi_0 \exp\left(-\halb \x^2/\sigma^2 \right),	\nonumber \\
	& \E(\x, 0)  = \mathbf{E}_0 \exp \left( -\halb \x^2/\sigma^2 \right), \qquad & \psi(\x, 0) = \psi_0 \exp\left(-\halb \x^2/\sigma^2 \right), 
	\label{eqn.ic.gauss}
\end{eqnarray}
Simulations are run in the domain $\Omega=[-\halb,+\halb]^2$ until a final time of $t=100$ using a DG scheme of approximation degree $N=3$ and taking an unstructured triangular mesh with $N_x=N_y=30$ elements along each edge of the domain. The time step is set to a constant value of $\Delta t = 10^{-2}$. We propose two different test cases, T1 and T2. In the first test (T1) we choose initial data that are compatible with the original vacuum Maxwell equations. We therefore set $\B_0 = 0$, $\E_0 = (0,0,1)$ and $\varphi_0 = \psi_0 = 0$. In the second test (T2) we choose $\B_0 = 0$, $\varphi_0=1$, $\psi_0 = 0$ and $\E_0 = 0$. In both tests the half width is chosen as $\sigma=0.05$. 
The computational results obtained for both tests are depicted in Figure \ref{fig.mm}. As expected, the method conserves total energy perfectly even over long times and the involutions are satisfied up to machine precision. 

\subsection{The incompressible Euler equations} 

{For the incompressible Euler equations we use the same collocation of ansatz spaces of the discrete velocity and pressure as in the case of linear acoustics, i.e. we use 
\eqref{eqn.acoustics.ansatz} with $\vh \in \mathcal{U}^N_h$ and $\ph \in \mathcal{W}^{N+1}_h$. }
We solve the two-dimensional Taylor-Green vortex, which is an exact stationary solution for the incompressible Euler equations and which reads 
\begin{equation}
v_1(\x,t) =  \sin(x)\cos(y), \qquad 
v_2(\x,t) = -\cos(x)\sin(y), \qquad 
p(\x,t) = p_0 + \frac{1}{4}(\cos(2x)+\cos(2y)).
\label{eq.tgv}
\end{equation}
The constant $p_0$ can be chosen arbitrarily. Here we set $p_0 = -\halb$ so that $p(\mathbf{0},t)=0$, which we also impose strongly in the linear algebraic system for the pressure in order to fix the constant inside the numerical scheme.  
{The discrete initial velocity field $\vh \in \mathcal{U}^N_h$ is computed as the discrete curl of the discrete vector potential $\tilde{\mathbf{A}}_h \in \mathcal{W}^{N+1}_h$, which is chosen as the continuous finite element representation of the continuous potential $\A = \sin(x) \sin(y)$.}
We run the problem on a sequence of successively refined unstructured triangular meshes in two space dimensions with a number of $N_x=N_y$ elements along each boundary edge of the square domain $\Omega = [0,2\pi]^2$. 
The simulations are run with schemes of nominal order one to six until a final time of $t=0.25$. The obtained
$L^2$ errors of the two velocity components and of the pressure are reported in Table
\ref{tab.euler.conv}. One can observe that the designed order of accuracy of $N+1$ has properly been reached for all schemes and for all variables.  

\begin{figure}[!htbp]
	\begin{center}
		\begin{tabular}{cc} 
			\includegraphics[width=0.45\textwidth]{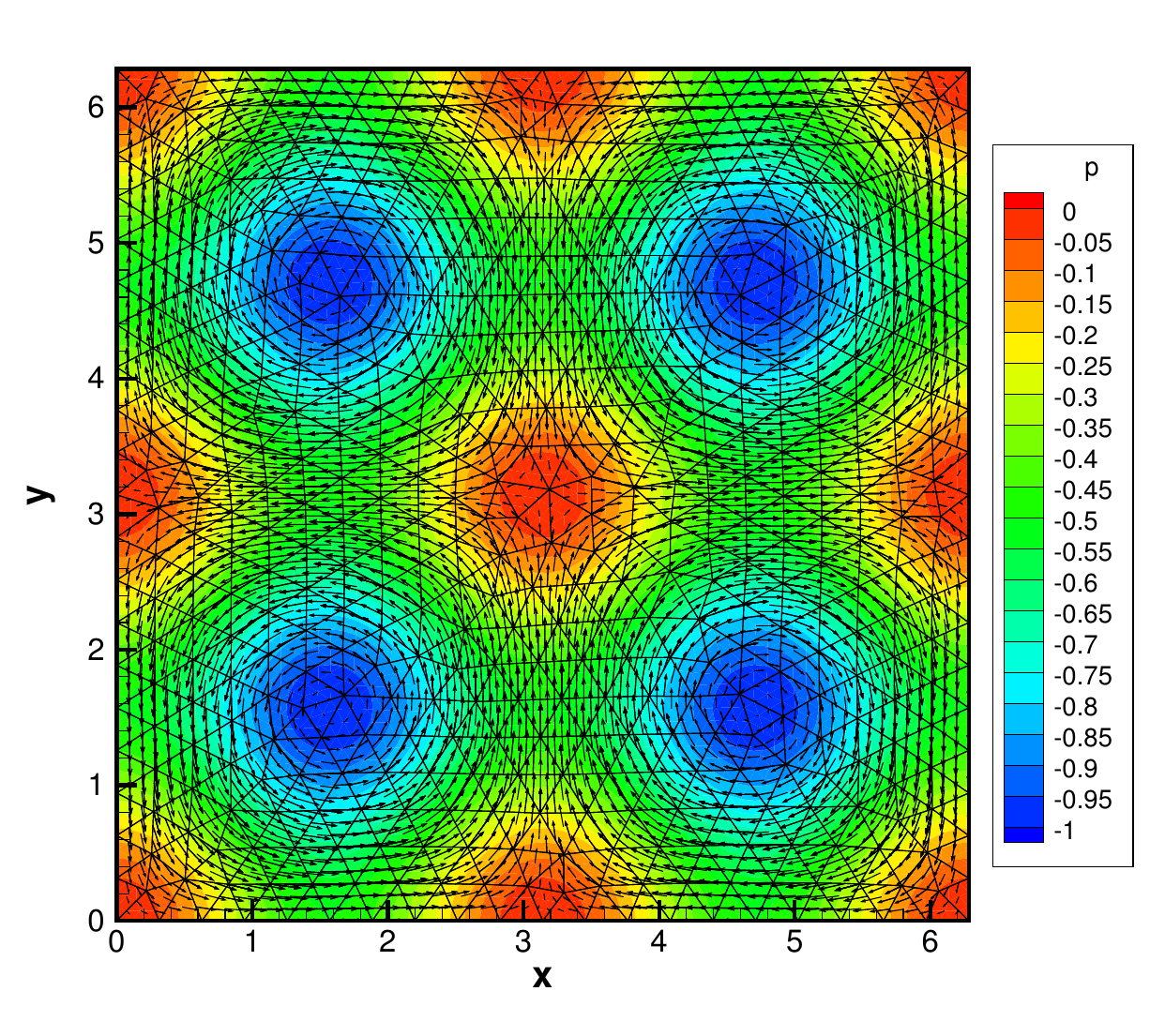}     & 
			\includegraphics[width=0.45\textwidth]{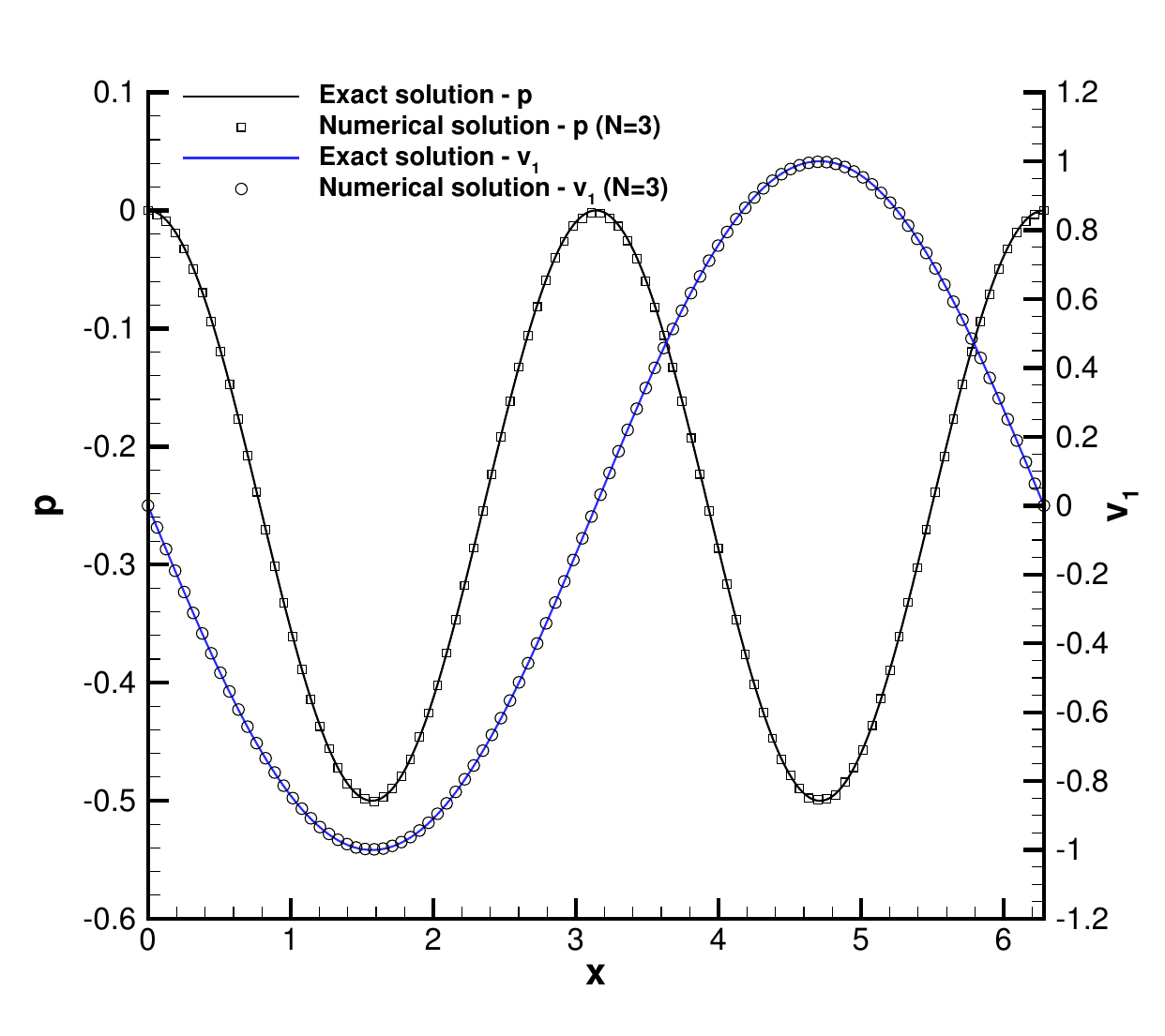}  \\  
			\includegraphics[width=0.45\textwidth]{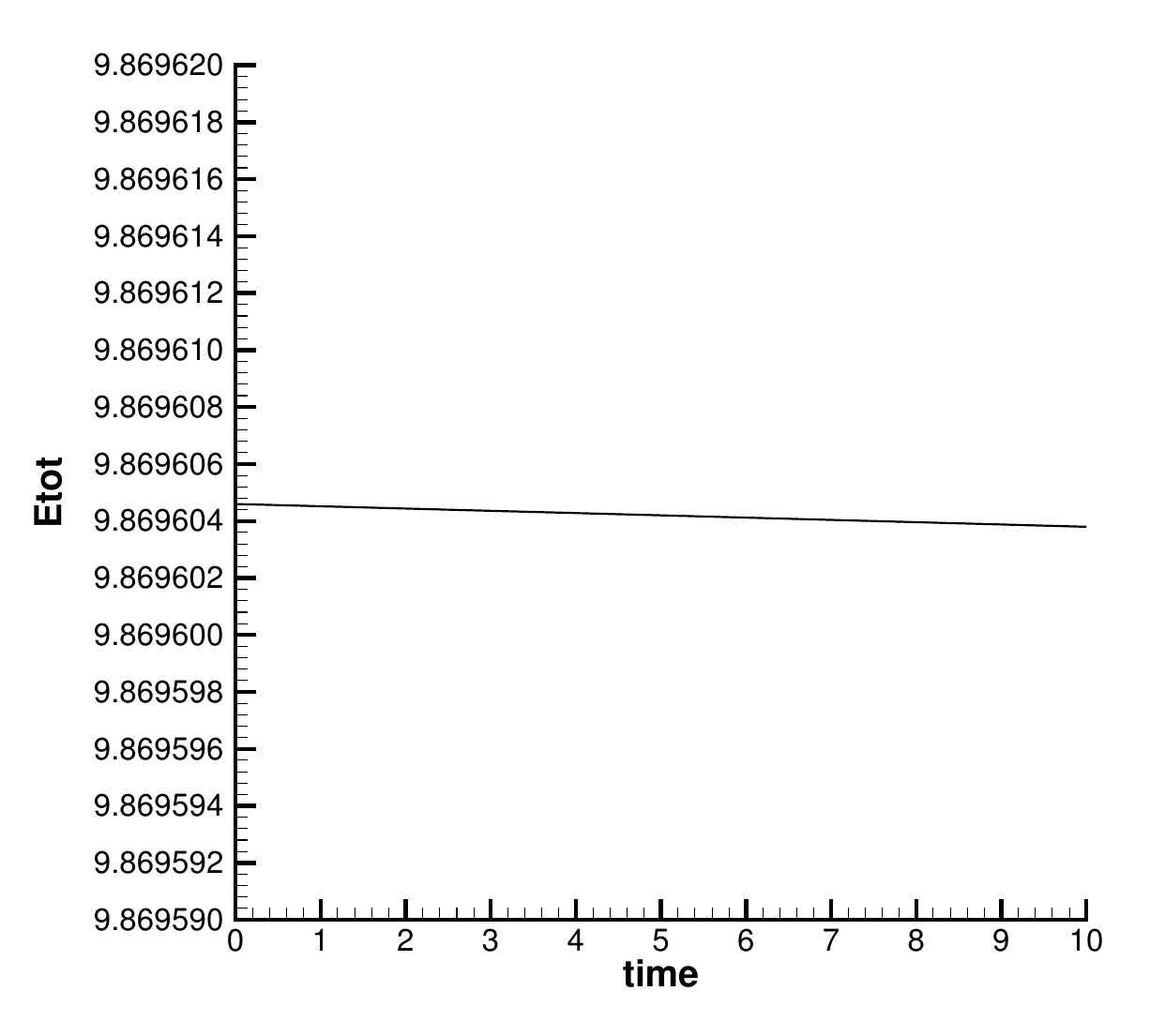}     & 
			\includegraphics[width=0.45\textwidth]{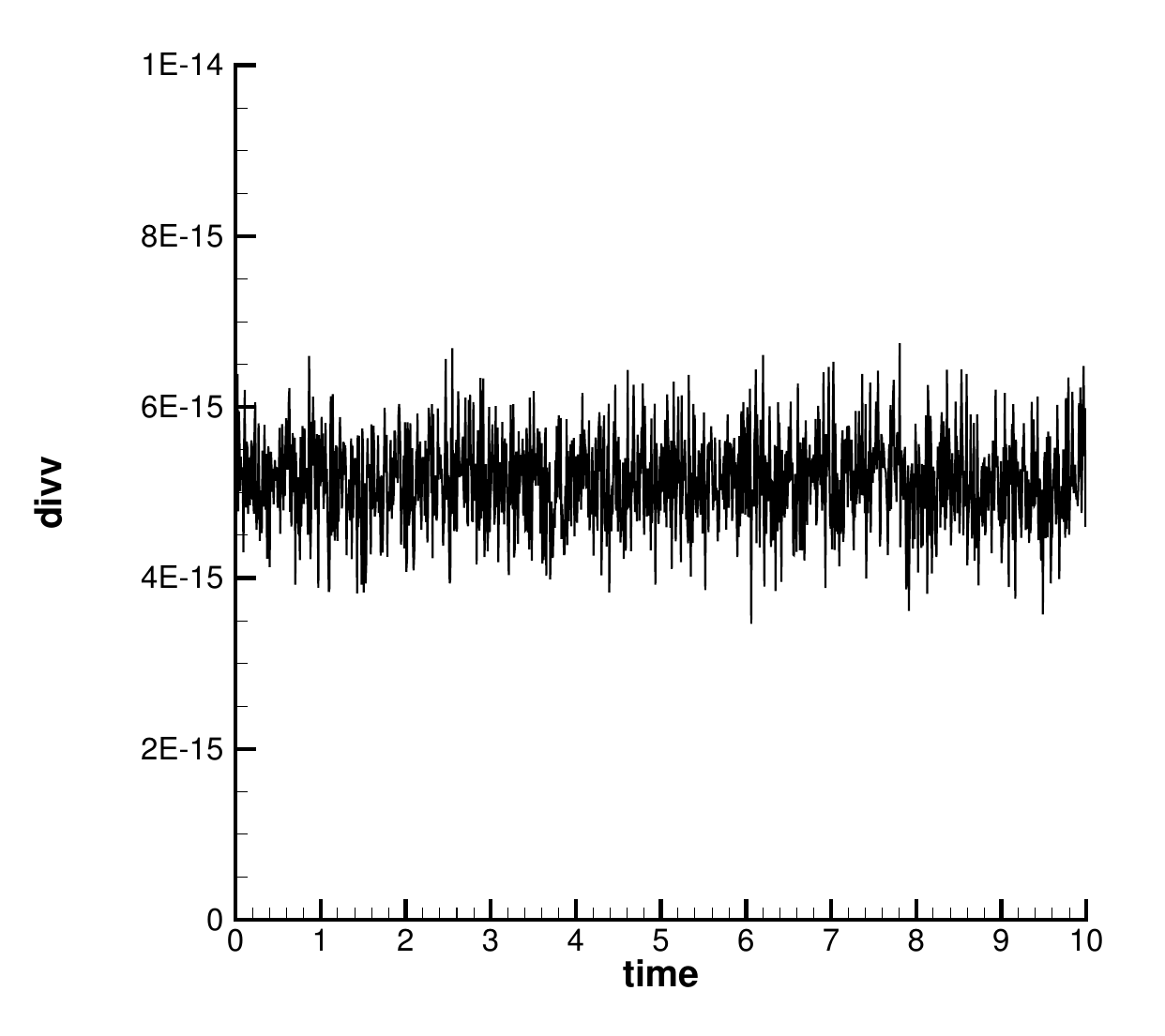}  
		\end{tabular} 
		\caption{Numerical solution of the Taylor-Green vortex at a final time of $t=10$ using the new structure-preserving semi-implicit DG scheme with $N=3$. Pressure contours and velocity vectors 
        (top left), 1D cut along the line $y=\pi$ for the pressure field and the velocity component $v_1$ and comparison with the exact solution (top right). Time evolution of the kinetic energy $E=\halb \mathbf{v}^2$ integrated over $\Omega$. Time evolution of the $L^{\infty}$ error of the velocity divergence (bottom right).  } 
		\label{fig.tgv}
	\end{center}
\end{figure}

We now run the same test again, but for large times, to verify that the discrete divergence of the velocity field remains zero up to machine precision and that the dissipation of the kinetic energy is small. Note that in the discretization of the nonlinear convective terms we use a dissipative numerical flux of the Ducros type.
The test is run with a DG scheme of approximation degree $N=3$ for the velocity field. As a consequence, the discrete pressure field is represented by continuous finite elements of degree $M=N+1=4$. Simulations are run until a final time of $t=10$ and the results are depicted in Figure \ref{fig.tgv}. 

\begin{table}  
		\caption{Numerical convergence results of the new structure-preserving DG scheme for the incompressible Euler equations and applied to the 2D Taylor-Green vortex problem for polynomial approximation degrees $N \in [0,5]$. The $L^2$ error norms refer to the velocity components $v_1$ and $v_2$ and the pressure $p$ at a final time of $t=0.25$. } 	
		\begin{center} 
				\renewcommand{\arraystretch}{1.1}
				\begin{tabular}{ccccccc} 
                    \hline 
                    \multicolumn{7}{c}{$N=0$} \\ 
					\hline
					$N_x = N_y$ & ${L^2}(v_1)$ & $\mathcal{O}(v_1)$  & ${L^2}(v_2)$ & $\mathcal{O}(v_2)$  
					& ${L^2}(p)$ & $\mathcal{O}(p)$  
					\\ 
					\hline
					20	& 1.6322E+00 &      & 1.6282E+00 &       & 1.2905E+00 &       \\ 
					40	& 8.0970E-01 & 1.0  & 8.1257E-01 &  1.0  & 6.1312E-01 & 1.1   \\ 
					60	& 5.3738E-01 & 1.0  & 5.4014E-01 &  1.0  & 4.0533E-01 & 1.0   \\ 
					80	& 4.0317E-01 & 1.0  & 4.0318E-01 &  1.0  & 3.0412E-01 & 1.0   \\ 
					100 & 3.2306E-01 & 1.0  & 3.2275E-01 &  1.0  & 2.4031E-01 & 1.1   \\ 
                    \hline 
                    \multicolumn{7}{c}{$N=1$} \\ 
					\hline
					$N_x = N_y$ & ${L^2}(v_1)$ & $\mathcal{O}(v_1)$  & ${L^2}(v_2)$ & $\mathcal{O}(v_2)$  
					& ${L^2}(p)$ & $\mathcal{O}(p)$  
					\\ 
					\hline
					20  & 1.0418E-01 &      & 1.0357E-01 &       & 8.6467E-02 &       \\ 
					40  & 3.0494E-02 & 1.8  & 2.9745E-02 &  1.8  & 2.1506E-02 & 2.0   \\ 
					60  & 1.5737E-02 & 1.6  & 1.4967E-02 &  1.7  & 9.5623E-03 & 2.0   \\ 
					80  & 9.2279E-03 & 1.9  & 8.9874E-03 &  1.8  & 5.3651E-03 & 2.0   \\ 
					100 & 6.0381E-03 & 1.9  & 6.0440E-03 &  1.8  & 3.4247E-03 & 2.0   \\ 
                    \hline 
                    \multicolumn{7}{c}{$N=2$} \\ 
					\hline
					$N_x = N_y$ & ${L^2}(v_1)$ & $\mathcal{O}(v_1)$  & ${L^2}(v_2)$ & $\mathcal{O}(v_2)$  
					& ${L^2}(p)$ & $\mathcal{O}(p)$  
					\\ 
					\hline
					20  & 3.5553E-03 &      & 3.5247E-03 &       & 3.7014E-03 &       \\ 
					40  & 4.0670E-04 & 3.1  & 4.0422E-04 &  3.1  & 4.2983E-04 & 3.1   \\ 
					60  & 1.1955E-04 & 3.0  & 1.1849E-04 &  3.0  & 1.2621E-04 & 3.0   \\ 
					80  & 4.9566E-05 & 3.1  & 4.9396E-05 &  3.0  & 5.2725E-05 & 3.0   \\ 
					100 & 2.5323E-05 & 3.0  & 2.5233E-05 &  3.0  & 2.6829E-05 & 3.0   \\ 
                    \hline 
                    \multicolumn{7}{c}{$N=3$} \\ 
					\hline
					$N_x = N_y$ & ${L^2}(v_1)$ & $\mathcal{O}(v_1)$  & ${L^2}(v_2)$ & $\mathcal{O}(v_2)$  
					& ${L^2}(p)$ & $\mathcal{O}(p)$  
					\\ 
					\hline
					10  & 1.5602E-03 &      & 1.6139E-03 &       & 2.3151E-03 &       \\ 
					20  & 1.0092E-04 & 4.0  & 1.0197E-04 &  4.0  & 1.2110E-04 & 4.3   \\ 
					30  & 2.0268E-05 & 4.0  & 2.0296E-05 &  4.0  & 2.3513E-05 & 4.0   \\ 
					40  & 6.4437E-06 & 4.0  & 6.4619E-06 &  4.0  & 7.3276E-06 & 4.1   \\ 
					80  & 4.3975E-07 & 3.9  & 4.3792E-07 &  3.9  & 4.5711E-07 & 4.0   \\ 
                    \hline 					
                    \multicolumn{7}{c}{$N=4$} \\ 
					\hline
					$N_x = N_y$ & ${L^2}(v_1)$ & $\mathcal{O}(v_1)$  & ${L^2}(v_2)$ & $\mathcal{O}(v_2)$  
					& ${L^2}(p)$ & $\mathcal{O}(p)$  
					\\ 
					\hline
					 5  & 2.4279E-03 &      & 2.4225E-03 &       & 8.5396E-03 &       \\ 
					10  & 8.7845E-05 & 4.8  & 8.8885E-05 &  4.8  & 1.1800E-04 & 6.2   \\ 
					20  & 2.2276E-06 & 5.3  & 2.1897E-06 &  5.3  & 3.5068E-06 & 5.1   \\ 
					30  & 2.7207E-07 & 5.2  & 2.6861E-07 &  5.2  & 4.4344E-07 & 5.1   \\ 
					40  & 6.1279E-08 & 5.2  & 6.0286E-08 &  5.2  & 1.0263E-07 & 5.1   \\ 
                    \hline 
                    \multicolumn{7}{c}{$N=5$} \\ 
					\hline
					$N_x = N_y$ & ${L^2}(v_1)$ & $\mathcal{O}(v_1)$  & ${L^2}(v_2)$ & $\mathcal{O}(v_2)$  
					& ${L^2}(p)$ & $\mathcal{O}(p)$  
					\\ 
					\hline
					5   & 3.5849E-04 &      & 3.8565E-04 &       & 2.4858E-03 &       \\ 
					10  & 5.7781E-06 & 6.0  & 5.9649E-06 &  6.0  & 2.0292E-05 & 6.9   \\ 
					15  & 4.1095E-07 & 6.5  & 4.1402E-07 &  6.6  & 8.8990E-07 & 7.7   \\ 
					20  & 6.7913E-08 & 6.3  & 6.5429E-08 &  6.4  & 1.1007E-07 & 7.3   \\ 
					30  & 5.1061E-09 & 6.4  & 5.0503E-09 &  6.3  & 7.4911E-09 & 6.6   \\ 
					\hline 					
				\end{tabular}
		\end{center}
		\label{tab.euler.conv}
\end{table}

\section{Conclusions}\label{sec:conclusions}

In this paper we have started to develop a new framework of discontinuous Galerkin finite element schemes on general unstructured simplex meshes that, in suitable combination with continuous Lagrange finite elements, achieve an \textit{exact discrete compatibility} with the two classical vector calculus identities $\nabla \times \nabla Z = 0$ and $\nabla \cdot \nabla \times \A = 0 $, for any generic scalar and vector field $Z$ and $\A$, respectively. 

The key idea behind our methods is the use of two suitable function spaces with their respective associated discrete nabla operators. The first approximation space $\mathcal{U}^N_h$ is of the discontinuous Galerkin type and is therefore based on piecewise polynomials of approximation degree $N$ that are allowed to jump across element interfaces. The second approximation space  $\mathcal{W}^{N+1}_h$ is the space of globally continuous Lagrange finite elements of degree $M=N+1$. 
The first discrete nabla operator, which we also call \textit{primary nabla operator}, is simply the exact derivative operator applied to the discrete solution represented in $\mathcal{W}^{N+1}_h$ by the continuous finite elements, projecting the result into the DG space $\mathcal{U}^N_h$. The corresponding \textit{discrete derivative} is therefore \textit{exact}, which is the first main key ingredient of our new method. In this paper we have also defined a suitable dual nabla operator, which uses the weak form of the equations and which allows to compute discrete gradients on the DG space $\mathcal{U}^N_h$, projecting the result back to the FEM space $\mathcal{W}^{N+1}_h$. 

The primary nabla operator being an exact derivative, it is obvious that the basic vector calculus identities must hold. More precisely, we were able to show a discrete Schwarz theorem, i.e. the symmetry of discrete second derivatives. As an immediate consequence, we obtain both discrete vector calculus identities. 

First, the discrete curl applied to a gradient computed via this primary discrete nabla operator is pointwise exactly zero inside each simplex element $T_k$. Furthermore, since the classical Lagrange finite elements are continuous, all tangential derivatives along the boundary $\partial T_k$ are necessarily continuous, since the FEM polynomials coincide in both elements adjacent to an element boundary. This provides at the same time trivially a locally and globally curl-free property of the resulting vector field in the case of linear acoustics. 

Second, the divergence of a curl computed via this primary discrete nabla operator vanishes locally pointwise everywhere inside each element. {By virtue of the same reasoning as above, also the tangential components of the gradient of a vector field are zero, and thus, after multiplication with the Levi-Civita tensor, the normal component of the magnetic field is continuous across element boundaries in the case of the Maxwell equations.}  

Or, to say it with the words of Leibniz: \textit{Omnibus ex nihilo ducendis sufficit unum}, since the discrete Schwarz theorem alone is enough to obtain as a consequence \textit{both} discrete vector calculus identities within one and the same numerical scheme and without needing to adapt the approximation spaces to the type of involution that one wants to preserve. 

The second key ingredient of our method is a new and universal recipe that allows to decide which variable needs to be defined in which approximation space. For this purpose we employ the framework of symmetric hyperbolic and thermodynamically compatible (SHTC) systems of Godunov and Romenski and coworkers. According to the SHTC framework the governing equations are always created in \textit{pairs}. For each single pair, the first type of equations is a mere consequence of the definitions via some abstract scalar or vector potentials. This type of PDE is usually the one endowed with the involutions that we want to preserve and therefore necessarily needs to be defined in the DG space $\mathcal{U}^N_h$. The second type of equations in each pair can be derived from an underlying variational principle via the classical Euler-Lagrange equations. The associated variables are thus naturally defined in the continuous FEM space $\mathcal{W}^{N+1}_h$.  

The proposed framework allows the construction of provably energy-preserving schemes in a very natural manner. 
Another interesting feature of the new schemes, at least from a practical point of view, is the possibility to make straightforward use of the Schur complement for the solution of the final linear algebraic systems, since the mass matrix of the DG scheme is block-diagonal or even diagonal for orthogonal DG basis functions and thus trivial to invert. We were able to show that the resulting algebraic systems are symmetric and positive definite and thus allow the use of efficient iterative solvers such as the matrix-free conjugate gradient method.  

We have applied our new numerical schemes to four different prototype systems: the equations of linear acoustics, the vacuum Maxwell equations and the Maxwell-GLM system. In all cases we have proven total energy conservation rigorously at the semi-discrete and at the fully-discrete level. We have also proven that the involutions are satisfied exactly at the discrete level. As a last system, which is outside the previous class of SHTC systems, we have applied our framework also to the incompressible Euler equations, just to show that the new approach is general enough to accommodate also naturally for more general elliptic-hyperbolic systems that are naturally outside the SHTC formalism, that it is indeed possible to combine our new DG method with existing standard DG schemes for nonlinear hyperbolic conservation laws {and that the new framework allows the straightforward construction of energy-stable schemes also for this type of PDE.} 

Future research will concern an extension to more general isoparametric elements, to the Euler equations of compressible gasdynamics and to Lagrangian hydrodynamics. We also plan to develop a semi-discrete and fully cell-centered version of our scheme that can be integrated in time via classical explicit Runge-Kutta methods. {In the future, we will also study the 
use of higher order symplectic time integrators for the scheme presented in this paper.}  

\section*{Acknowledgments}

	M.D. was financially supported by the Italian Ministry of Education, University 
	and Research (MIUR) in the framework of the PRIN 2022 project \textit{High order structure-preserving semi-implicit schemes for hyperbolic equations} and via the  Departments of Excellence  Initiative 2018--2027 attributed to DICAM of the University of Trento (grant L. 232/2016). 
	M.D. was also funded by the Fondazione Caritro via the SOPHOS project, the European Union Next Generation EU projects PNRR Spoke 7 CN HPC and PNRR Spoke 7 RESTART, as well as by the European Research Council (ERC) under the European Union's Horizon 2020 research and innovation programme, Grant agreement No. ERC-ADG-2021-101052956-BEYOND. 
	Views and opinions expressed are however those of the author(s) only and
	do not necessarily reflect those of the European Union or the European Research
	Council. Neither the European Union nor the granting authority can be held
	responsible for them.   
	{The authors would kindly like to thank the two anonymous referees for their constructive comments and remarks that helped to improve the quality and readability of this paper. }

    \vspace{3mm}
    
	M.D. is member of the Gruppo Nazionale per il Calcolo Scientifico dell'Istituto Nazionale di Alta Matematica (GNCS-INdAM) 
    and dedicates all his ideas and scientific contributions in this paper with great affection to Wilma Pilati, in deep thankfulness for the many pleasant, stimulating and very inspiring philosophical discussions. 

\vspace{5mm} 

%
%
%
\bibliographystyle{plain}
\bibliography{./biblio}


\end{document}